# NILPOTENT CONES AND THEIR REPRESENTATION THEORY

P. BROSNAN, G. PEARLSTEIN, C. ROBLES

*Dedicated to Steven Zucker in recognition of his contributions to Hodge theory, on the occasion of his 65th birthday.*

Abstract. We describe two approaches to classifying the possible monodromy cones $C$ arising from nilpotent orbits in Hodge theory. The first is based upon the observation that $C$ is contained in the open orbit of any interior point $N \in C$ under an associated Levi subgroup determined by the limit mixed Hodge structure. The possible relations between the interior of $C$ its faces are described in terms of signed Young diagrams.

The second approach is to understand the Tannakian category of nilpotent orbits via a category $D$ introduced by Deligne in a letter to Cattani and Kaplan. In analogy with Hodge theory, there is a functor from $D$ to a subcategory $\hat{D} \subset D$ of $\mathrm{SL}_2$-orbits. We prove that these fibers are, roughly speaking, algebraic. We also give a correction to a result [16] of K. Kato.

## 1. Introduction

The object of Steven Zucker's first published paper [27] was the study of normal functions arising from algebraic cycles and the Hodge conjecture. More precisely, let $X \subset \mathbb{P}^m$ be a smooth projective variety of dimension $2d$ and $\zeta$ be a primitive Hodge class of type $(d, d)$ on $X$. Let $Y$ be a smooth hyperplane section of $X$. Then, the long exact sequence for relative cohomology of the pair $(X, Y)$ gives

$$\cdots \to H^{2d-1}(X) \stackrel{i^*}{\hookrightarrow} H^{2d-1}(Y) \to H^{2d}(X, Y) \to H^{2d}(X) \to \cdots$$

*Date*: January 21, 2016.

2010 *Mathematics Subject Classification.* 14D07, 32G20, 17B08, 32S35.

Patrick Brosnan is partially supported by NSF grant DMS 1361159. Gregory Pearlstein is partially supported by NSF grant DMS 1361120. Colleen Robles is partially supported by NSF grants DMS 02468621 and 1361120.

where $i : Y \hookrightarrow X$ is inclusion. Letting $H^{2d-1}_{\rm van}(Y)$ denote the cokernel of $i^*$ and pulling back the above sequence along the morphism of Hodge structure $\mathbb{Z}(-d) \to H^{2d}(X)$ defined by $\zeta$ determines an extension

$$0 \to H^{2d-1}_{\rm van}(Y) \to E \to \mathbb{Z}(-d) \to 0$$

in the category of mixed Hodge structures. The set of all such extensions is the intermediate Jacobian $J(H^{2d-1}_{\rm van}(Y))$.

Applying the construction of the previous paragraph to the smooth fibers of a Lefschetz pencil of hyperplane sections of $X$ yields the prototypical example of a normal function $\nu_\zeta$. Moreover, in this context, there is a 1-1 correspondence $\zeta \leftrightarrow \nu_\zeta$ between normal functions and primitive, integral Hodge classes on $X$.

More precisely, Zucker proved [28] using $L^2$ methods that if $S$ is a curve with smooth completion $j : S \hookrightarrow \bar{S}$ and $\mathcal{V} \to S$ is a variation of Hodge structure of weight $k$ then $H^i(\bar{S}, j_*\mathcal{V})$ carries a functorial Hodge structure of weight $i+k$ (Theorem (7.12)). Furthermore (Theorem (9.2)), if $\mathcal{V}$ is pure of weight $2p-1$ then the cohomology classes of normal functions on $S$ surject onto the integral $(p,p)$ classes in $H^1(\bar{S}, j_*\mathcal{V})$.

One of the key tools in Zucker's proofs are W. Schmid's orbit theorems [23]. Roughly speaking, the nilpotent orbit theorem asserts that a variation of Hodge structure $\mathcal{V} \to S$ admits a local approximation near a point in the boundary of any normal crossing compactification $S \hookrightarrow \bar{S}$ by a nilpotent orbit

$$\theta(z_1, \dots, z_r) = \exp(\sum_j z_j N_j) F_\infty \tag{1.1}$$

determined by the local monodromy logarithms $N_1, \dots, N_r$ and limit Hodge filtration $F_\infty$ of $\mathcal{V}$. The SL$_2$-orbit theorem [23, 4] further asserts that $\theta$ can be approximated by an auxiliary nilpotent orbit which is governed by a representation $\rho$ of $SL_2(\mathbb{R})^r$. In this way, the norm of a flat multivalued section $\sigma$ of $\mathcal{V}$ is determined by it weights with respect to $\rho$. The SL$_2$-orbit theorem also implies that the limit Hodge filtration $F_\infty$ is part of a limit mixed Hodge structure $(F_\infty, W)$.

Moving beyond families of smooth projective varieties, Deligne conjectured that given a surjective, quasiprojective morphism $\bar{f} : \bar{X} \to \bar{S}$ there should exist a Zariski open subset $S \subset \bar{S}$ over which $\bar{f}$ restricts to give a variation of mixed Hodge structure [7] on the cohomology of the fibers. Furthermore, there should be a category of good variations of mixed Hodge structure which contains every variation of mixed Hodge structure of geometric origin, and has all of the salient features of the pure case.

In [26], Steenbrink and Zucker defined a category of admissible variations of mixed Hodge structure over a curve $S$, and proved that in this category one had limit mixed Hodge structures. Moreover, if $\mathcal{V} \to S$ is an admissible variation of graded-polarized mixed Hodge structure over $S$ and $j : S \hookrightarrow \bar{S}$ is a smooth completion of $S$ then the cohomology groups $H^i(\bar{S}, R^k j_*\mathcal{V})$ carry functorial mixed Hodge structures (cf. Theorem (4.1) [26]). In [15], Kashiwara defined a category of admissible variations of graded-polarized mixed Hodge structure in several variables using a curve test, and christened the associated nilpotent orbits *infinitesimal mixed Hodge modules.*

In particular, the category of infinitesimal mixed Hodge modules (IMHM) in a fixed number of variables is an abelian tensor category (4.3.3 and 5.2.6 [15]) which becomes a neutral Tannakian category when equipped with the functor $\omega$ which takes an IMHM to the underlying $\mathbb{R}$-vector space. The category of IMHM has a natural subcategory corresponding to nilpotent orbits with limit mixed Hodge structure which is split over $\mathbb{R}$. The Tannakian Galois group of the category of split orbits in one variable is described by the first two authors in [1].

In the case of nilpotent orbits of pure Hodge structures in one variable, a split orbit is the same thing as $\mathrm{SL}_2$-orbit: If $D$ is a period domain upon which the Lie group $G_{\mathbb{R}}$ acts transitively by automorphisms then a nilpotent orbit $\theta(z)$ with values in $D$ is an $\mathrm{SL}_2$-orbit if there exists a representation $\rho : \mathrm{SL}_2(\mathbb{R}) \to G_{\mathbb{R}}$ such that

$$\theta(g.\sqrt{-1}) = \rho(g).\theta(\sqrt{-1}) \tag{1.2}$$

for all $g \in \mathrm{SL}_2(\mathbb{R})$. A classification of such orbits may be deduced from (i) Lemma (6.24) of [23], and (ii) the classification of nilpotent $N \in \mathfrak{g}_{\mathbb{R}} = \mathrm{Aut}(V_{\mathbb{R}}, Q)$ (which is reviewed in §2.3). A full classification in the case of orbits into a Mumford–Tate domain is given in [21].

One of the intricacies of the $\mathrm{SL}_2$-orbit theorem [4] in several variables is the construction of the system

$$(\hat{N}_1, H_1, \hat{N}_1^+), \ldots, (\hat{N}_r, H_r, \hat{N}_r^+) \tag{1.3}$$

of commuting $\mathrm{SL}_2$-triples attached to the nilpotent orbit (1.1). In [8], Deligne gave a purely linear algebraic construction of (1.3) via an iterative construction which ensures that each $N_j$ is a sum of $\hat{N}_j$ and a collection of highest weight vectors for $(\hat{N}_j, H_j, \hat{N}_j^+)$. Published accounts of the resulting Deligne systems appears in [24, 14, 1].

The data of an IMHM consists of a set of commuting nilpotent endomorphisms $N_1, \dots, N_r$ together with Hodge and weight filtrations which satisfy a number of compatibility conditions. In [16], Kato observed that for each non-negative real number $a$, the substitution

$$N_j \mapsto \phi^a(N_j) = \sum_{k=0}^{j-1} \frac{a^k}{k!} N_{j-k}$$

defines a functor $\phi^a : \mathrm{IMHM} \to \mathrm{IMHM}$ such that $\phi^a \circ \phi^b = \phi^{a+b}$ and $\phi^0$ is the identity. Moreover, the functor $\phi$ extends to a category DH of Deligne–Hodge systems which contains IMHM as a subcategory. Kato further claimed that for any object $\theta$ of DH there exists an $a \geq 0$ such that $\phi^b\,\theta$ belongs to IMHM for all $b > a$.

As Kato explains in the introduction to [16], one of his motivations to study Deligne–Hodge systems was to have a framework to study degenerations of Hodge structure which are not polarizable. Such a framework could potentially be very useful in the study of degenerations of motives over non-archimedean local fields. However, in §6.2 we construct an explicit example of a two variable Deligne–Hodge system $\theta$ which violates Kato's assertion (i.e. $\phi^b\,\theta$ is never an IMHM). In §6.4 we show that Kato's claim is true for Deligne–Hodge systems which satisfy a suitable graded-polarization condition.

Accordingly, one can study IMHM in several variables by using the results of [21] to classify the possible several variable $\mathrm{SL}_2$-orbits with data (1.3), and then impose the representation theoretic conditions required to extend $\hat{N}_j$ to a candidate $N_j$. This second step can be done in the category DH. Application of $\phi^b$ for $b$ sufficiently large then produces the required IMHM. In §6.5 we show that the set of all Deligne systems with fixed data (1.3) forms an algebraic variety.

The approach outlined in the previous paragraph attempts to construct the monodromy cone

$$\mathcal{C} = \left\{ \sum_{j=1}^{r} a_j N_j \mid a_1, \dots, a_r > 0 \right\} \tag{1.4}$$

of an IMHM starting from the edges of the closure of $\mathcal{C}$. Alternatively, one can try an construct orbits starting from an element of the interior of $\mathcal{C}$. This second approach, *for nilpotent orbits of pure Hodge structure*, is the subject of §3–5.

A rough outline is as follows: Without loss of generality we can pass to the case where (1.1) is a nilpotent orbit with limit mixed Hodge structure split over $\mathbb{R}$. For any $N \in \mathcal{C}$, it then follows that $e^{zN}F_\infty$ is an $\mathrm{SL}_2$-orbit by the results of Cattani, Kaplan and Schmid [2, 4]. The results of [21] classify the possible pairs $(N, F_\infty)$.

To reverse this process, select an $\mathrm{SL}_2$-orbit $e^{zN}F_\infty$ with associated representation $\rho : \mathrm{SL}_2 \to G_\mathbb{R}$. Let $G_\mathbb{R}^{0,0}$ be the connected subgroup of $G_\mathbb{R}$ consisting of elements which preserve the limit mixed Hodge structure of $e^{zN}F_\infty$. Let $\mathcal{N}$ denote the orbit of $N$ under the adjoint action of $G_\mathbb{R}^{0,0}$. Then, by Lemma (3.5) and Corollary (3.6) it follows that if $\mathcal{C}$ is the cone (1.4) of a several variable nilpotent orbit with limit Hodge filtration $F_\infty$ and $N \in \mathcal{C}$ then $\mathcal{C} \subset \mathcal{N}$. Sections 4 and 5 implement this process for a number of examples related to period domains of weight 2.

Section 2 is of a different nature: The set $\mathrm{Nilp}(\mathfrak{g}_\mathbb{R})$ of nilpotent elements in a real semisimple (or reductive) Lie algebra is a classical, and very well understood, object of study in representation theory; for an excellent introduction see [6] and the references therein. In particular, there is a great deal of Hodge theoretic information that one can glean from the representation theorists' understanding of $\mathrm{Nilp}(\mathfrak{g}_\mathbb{R})$. In §2 we review the classification of nilpotent $N \in \mathfrak{g}_\mathbb{R}$ by signed Young diagrams. Particularly noteworthy here are (i) Đoković's Theorem 2.21 characterizing a partial order on the $\mathrm{Ad}(G_\mathbb{R})$ conjugacy classes of nilpotent elements $N \in \mathfrak{g}_\mathbb{R}$, and (ii) the description in §2.5 of the signed Young diagram associated to a polarized mixed Hodge structure $(F, W)$: as illustrated in §5.2 together these provide representation theoretic constraints on the degenerations associated with the faces of a nilpotent cone underlying a nilpotent orbit on a period domain $D$.

**Acknowledgments.** In 1995, Steven Zucker was one of the instructors at a workshop for graduate students and recent PhD's held in connection with the AMS Summer Institute in Algebraic Geometry held at UC Santa Cruz. This is where the first two authors became acquainted, and we thank Steven for all of the things we have learned from him. The third author would like to thank William McGovern, Mark Green and Phillip Griffiths for helpful discussions and correspondence in connection with this work.

## 2. Classification of nilpotent endomorphisms

The nilpotent elements $N$ is a classical Lie algebra $\mathfrak{g}_\mathbb{R}$ are classified by partially signed Young diagrams. This classification is up to the action of $\mathrm{Ad}(G_\mathbb{R})$; so what are really being classified are the $\mathrm{Ad}(G_\mathbb{R})$–orbits in $\mathrm{Nilp}(\mathfrak{g}_\mathbb{R})$, of which there are only finitely many.[1] Collingwood and McGovern's [6] is an excellent reference for the material in this section.

It is convenient to begin with the classification of nilpotent endomorphisms over $\mathbb{C}$ by (unsigned) Young diagrams; the basic idea is that (i) the Young diagram encodes the Jordan normal form of $N \in \mathfrak{g}_\mathbb{C}$, and (ii) the Jordan normal form determines $N$ up to the action of $\mathrm{Ad}(G_\mathbb{C})$. Before reviewing the classifications over $\mathbb{R}/\mathbb{C}$ we need to recall the notion of a standard triple.

### 2.1. Standard triples.

Let $\mathfrak{g}$ be a Lie algebra defined over $\Bbbk = \mathbb{R}$ or $\mathbb{C}$. A *standard triple* in $\mathfrak{g}$ is a set of three elements $\{N^+, Y, N\} \subset \mathfrak{g}$ such that

$$[Y, N^+] \;=\; 2\,N^+\,, \quad [N^+, N] \;=\; Y \quad \text{and} \quad [Y, N] \;=\; -2\,N\,.$$

Note that $\{N^+, Y, N\}$ span a *3–dimensional semisimple subalgebra* (TDS) of $\mathfrak{g}$ isomorphic to $\mathfrak{sl}(2, \Bbbk)$. We call $Y$ the *neutral element*, $N$ the *nilnegative element* and $N^+$ the *nilpositive element*, respectively, of the standard triple.

**Theorem 2.1** (Jacobson–Morosov)**.** *Every nilpotent $N \in \mathfrak{g}$ can be realized as the nilnegative of a standard triple.*

---

[1] In contrast, there are infinitely many $\mathrm{Ad}(G_\mathbb{R})$ semisimple orbits in $\mathfrak{g}_\mathbb{R}$.

*Example* 2.2. The matrices

$$\mathbf{n}^+ \;=\; \begin{pmatrix} 0 & 1 \\ 0 & 0 \end{pmatrix}\,, \quad \mathbf{y} \;=\; \begin{pmatrix} 1 & 0 \\ 0 & -1 \end{pmatrix} \quad \text{and} \quad \mathbf{n} \;=\; \begin{pmatrix} 0 & 0 \\ 1 & 0 \end{pmatrix} \tag{2.3}$$

form a standard triple in $\mathfrak{sl}(2,\mathbb{R})$; while the matrices

$$\bar{\mathbf{e}} \;=\; \tfrac{1}{2}\begin{pmatrix} \mathbf{i} & 1 \\ 1 & -\mathbf{i} \end{pmatrix}\,, \quad \mathbf{z} \;=\; \begin{pmatrix} 0 & \mathbf{i} \\ -\mathbf{i} & 0 \end{pmatrix} \quad \text{and} \quad \mathbf{e} = \tfrac{1}{2}\begin{pmatrix} -\mathbf{i} & 1 \\ 1 & \mathbf{i} \end{pmatrix} \tag{2.4}$$

form a standard triple in $\mathfrak{su}(1,1)$.

2.2. **Nilpotents over $\mathbb{C}$.** The classification of the nilpotent elements in classical, complex, semisimple Lie algebras (which is due to Gerstenhaber [11], and Springer and Steinberg [25]) is given by partitions or, equivalently, Young diagrams.

Let $V_\mathbb{C}$ be a $\mathbb{C}$–vector space of dimension $n$ and fix a nilpotent element $N \in \mathrm{End}(V_\mathbb{C})$. Let $\mathfrak{sl}_2\mathbb{C} \subset \mathrm{End}(V_\mathbb{C})$ be the TDS spanned by a standard triple containing $N$ as the nilnegative element (§2.1). Let

$$V_\mathbb{C} \;=\; \bigoplus_{\ell\geq 0} V(\ell) \tag{2.5}$$

be the $\mathfrak{sl}_2\mathbb{C}$–decomposition of $V_\mathbb{C}$; here each $V(\ell) \simeq (\mathrm{Sym}^\ell\mathbb{C}^2)^{m_\ell}$ is the direct sum of $m_\ell$ irreducible $\mathfrak{sl}_2\mathbb{C}$–modules of dimension $\ell+1$. In particular, $V(\ell)$ admits a basis of the form

$$\{N^a v_i \mid 1 \leq i \leq m_\ell\,,\ 0 \leq a \leq \ell\}\,.$$

Here $N^\ell v_i \neq 0$ and $N^{\ell+1} v_i = 0$. Each

$$\{N^a v_i \mid 0 \leq a \leq \ell\}$$

is an $N$*–string* of length $\ell+1$ and we think of $V(\ell)$ as spanned by $m_\ell$ of these $N$–strings. Let

$$P(\ell) \;:=\; \mathrm{span}_\mathbb{C}\{v_i \mid 1 \leq i \leq m_\ell\} \;\subset\; V(\ell)$$

be the subspace of highest weight vectors. (In Hodge–theoretic language, this is the vector space of $N$*–primitive vectors* in $V(\ell)$.)

Note that $\sum m_\ell(\ell+1) = n$. So we may associate to the nilpotent $N$ a partition

$$\mathbf{d} \;=\; [d_i] \;=\; [(\ell+1)^{m_\ell}]_{0\leq\ell\in\mathbb{Z}}$$

of $n$; here $(\ell+1)^{m_\ell}$ indicates that the part $d_i = \ell+1$ occurs $m_\ell$ times. The partition $\mathbf{d}$ is identified with the Young diagram whose $i$–th row contains $d_i$ boxes. For example, the partition $\mathbf{d} = (4,2,2,1)$ of $n = 9$ is identified with the Young diagram

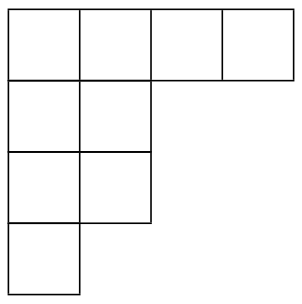

.

We think of each row of the Young diagram as representing an $N$–string. In this example, we have

| | | | |
|---|---|---|---|
| $u$ | $Nu$ | $N^2u$ | $N^3u$ |
| $v$ | $Nv$ | | |
| $w$ | $Nw$ | | |
| $x$ | | | |

and $V_{\mathbb{C}} = V(3) \oplus V(1) \oplus V(0)$ with $m_3 = 1 = m_0$ and $m_2 = 2$.

2.2.1. $G_{\mathbb{C}} = \mathrm{Aut}(V_{\mathbb{C}})$. The Jordan normal form for elements of $\mathfrak{g}_{\mathbb{C}} = \mathrm{End}(V_{\mathbb{C}})$ implies that two nilpotents $N_1, N_2 \in \mathfrak{g}_{\mathbb{C}}$ lie in the same $\mathrm{Ad}(G_{\mathbb{C}})$–orbit if and only if the corresponding partitions $\mathbf{d}_1$ and $\mathbf{d}_2$ are equal. That is, the $\mathrm{Aut}(V_{\mathbb{C}})$–orbits $\mathcal{N}$ in $\mathrm{Nilp}(\mathrm{End}(V_{\mathbb{C}}))$ are indexed by partitions of $n = \dim V$; equivalently, they are indexed by Young diagrams of size $n$.

*Example* 2.6. The nilpotent $\mathrm{Aut}(V_{\mathbb{C}})$–conjugacy classes $\mathcal{N} \subset \mathrm{End}(V_{\mathbb{C}})$ for $n = 5$ are indexed by

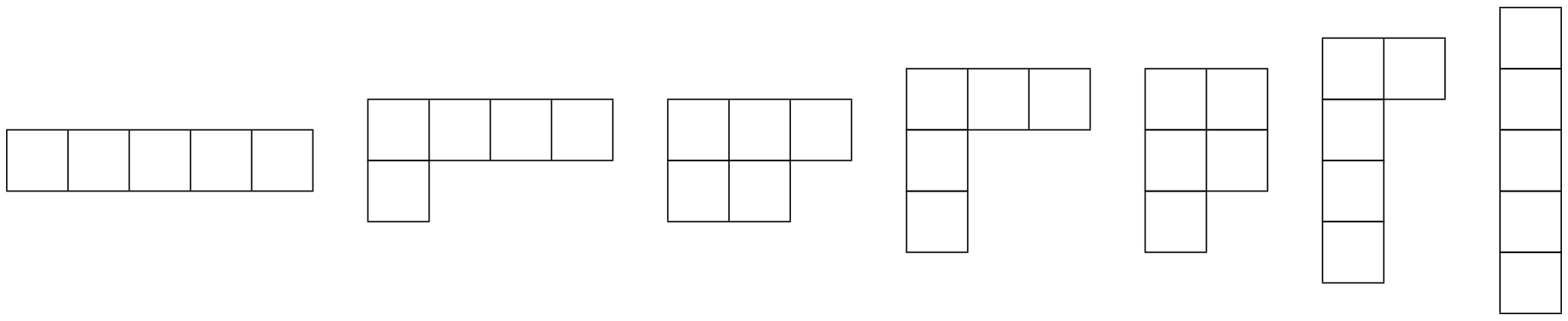

- Given a nilpotent $N \in \mathcal{N}$ in the conjugacy class indexed by the Young diagram 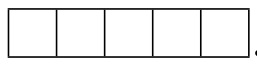, $V$ admits a basis of the form $\{v, Nv, \ldots, N^4v\}$ with $N^5v = 0$.
- Given a nilpotent $N \in \mathcal{N}$ in the conjugacy class indexed by the Young diagram 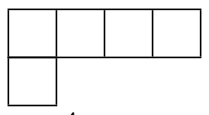, $V$ admits a basis of the form $\{u\,;\, v, Nv, \ldots, N^3v\}$ with $Nu = 0$ and $N^4 = 0$.

- And so on.

*Remark* 2.7. Note that the trivial nilpotent conjugacy class $\mathcal{N} = \{0\}$ is indexed by the vertical partition $[1^n]$.

2.2.2. $G_\mathbb{C} = \mathrm{Aut}(V_\mathbb{C}, Q)$. Fix $w \in \mathbb{Z}$ and let $Q$ be a nondegenerate bilinear form on $V_\mathbb{C}$ satisfying

$$Q(u,v) \ = \ (-1)^w Q(v,u)\,.$$

Set $G_\mathbb{C} = \mathrm{Aut}(V_\mathbb{C}, Q)$ and $\mathfrak{g}_\mathbb{C} = \mathrm{End}(V_\mathbb{C}, Q)$. Given a nonzero $N \in \mathrm{Nilp}(\mathfrak{g}_\mathbb{C})$, we may assume that the TDS is contained in $\mathfrak{g}_\mathbb{C}$ (cf. Jacobson and Morosov's Theorem 2.1). Then

$$Q_\ell(u,v) \ := \ Q(u, N^\ell v)$$

defines a non–degenerate bilinear form on $P(\ell)$.

- If $w+\ell$ is even, the $Q_\ell$ is symmetric.
- If $w+\ell$ is odd, then $Q_\ell$ is skew–symmetric. This implies that $m_\ell$ is even. So, if $w$ is even/odd, then the even/odd parts of $\mathbf{d}$ must occur with even multiplicity.

**Theorem 2.8** (Symplectic algebras ($w$ odd))**.** *Let $Q$ be a skew-symmetric bilinear form on a complex vector space $V_\mathbb{C}$, and set $G_\mathbb{C} = \mathrm{Aut}(V_\mathbb{C}, Q)$ with Lie algebra $\mathfrak{g}_\mathbb{C} = \mathrm{End}(V_\mathbb{C}, Q)$. Then the $\mathrm{Ad}(G_\mathbb{C})$–orbits in $\mathrm{Nilp}(\mathfrak{g}_\mathbb{C})$ are indexed by the partitions of $2m = \dim V_\mathbb{C}$ in which the odd parts occur with even multiplicity.*

*Example* 2.9. Suppose that $G_\mathbb{C} = \mathrm{Sp}(6,\mathbb{C})$. The nilpotent conjugacy classes in $\mathfrak{g}_\mathbb{C}$ are enumerated by the partitions

$$[6]\,,\ [4,2]\,,\ [4,1^2]\,,\ [3^2]\,,\ [2^3]\,,\ [2^2,1^2]\,,\ [2^2,1^4]\,,\ [1^6]\,.$$

The corresponding Young diagrams are

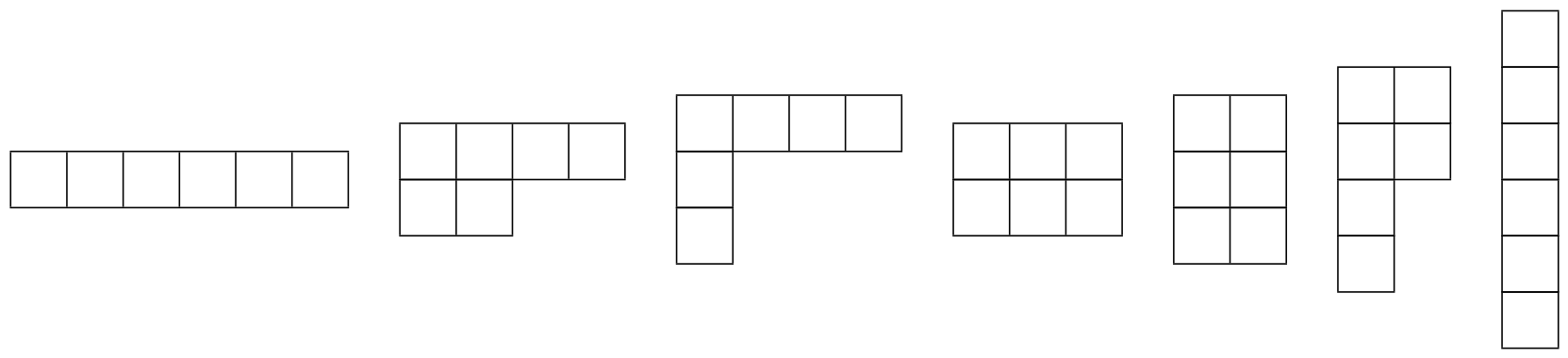

We say that a partition is *very even* if all parts $d_i$ are even and occur with even multiplicity.

**Theorem 2.10** (Orthogonal algebras ($w$ even)). *Let $Q$ be a symmetric bilinear form on a complex vector space $V_\mathbb{C}$, and set $\mathfrak{g}_\mathbb{C} = \mathrm{End}(V_\mathbb{C}, Q)$.*

(a) *Let $G_\mathbb{C} = \mathrm{Aut}(V_\mathbb{C}, Q)$. The $\mathrm{Ad}(G_\mathbb{C})$–orbits in $\mathrm{Nilp}(\mathfrak{g}_\mathbb{C})$ are indexed by the partitions of $n = \dim V_\mathbb{C}$ in which the even parts occur with even multiplicity.*

(b) *Let $G^\circ_\mathbb{C} = \mathrm{SO}(n, \mathbb{C}) \subset G_\mathbb{C}$. The $\mathrm{Ad}(G^\circ_\mathbb{C})$–orbits in $\mathrm{Nilp}(\mathfrak{g}_\mathbb{C})$ are indexed by partitions $\mathbf{d} = [d_i]$ of $n$ in which the even parts occur with even multiplicity, and with the caveat that a very even partition is associated with two distinct orbits.*

*Example* 2.11. Suppose that $G_\mathbb{C} = \mathrm{O}(7, \mathbb{C})$. The nilpotent conjugacy classes in $\mathfrak{g}_\mathbb{C} = \mathfrak{so}(7, \mathbb{C})$ are enumerated by the partitions

$$[7]\,,\ [5,1^2]\,,\ [3^2,1]\,,\ [3,2^2]\,,\ [3,1^4]\,,\ [2^2,1^3]\,,\ [1^7]\,.$$

The corresponding Young diagrams are

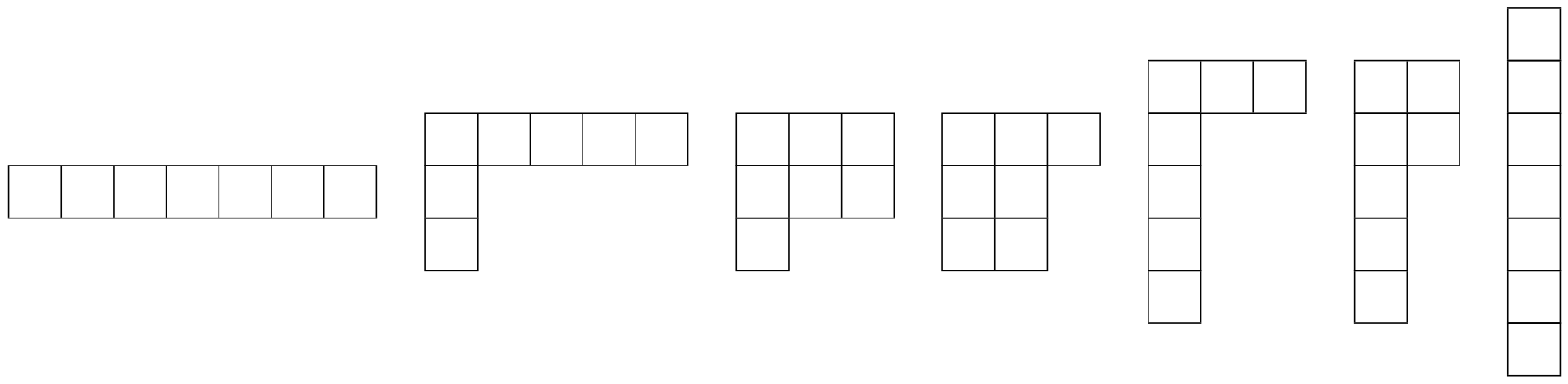

2.3. **Nilpotents over $\mathbb{R}$ (signed Young diagrams).** Let $V_\mathbb{R}$ be a real vector space of dimension $n$. Given $w \in \mathbb{Z}$, fix a nondegenerate bilinear form $Q : V_\mathbb{R} \times V_\mathbb{R} \to \mathbb{R}$ satisfying $Q(u, v) = (-1)^w Q(v, u)$. Set

$$G_\mathbb{R} \ = \ \mathrm{Aut}(V_\mathbb{R}, Q)\,.$$

The classification of $\mathrm{Ad}(G_\mathbb{R})$–conjugacy classes of nilpotent $N \in \mathfrak{g}_\mathbb{R}$ is due to Springer and Steinberg [25], and is given by (partially) signed Young diagrams. A *signed Young diagram* is a Young diagram in which the boxes of a fixed row are either labeled with alternating $\pm$ signs, or are left blank. For the real forms $G_\mathbb{R}$ under consideration, the blank rows occur with even multiplicity.

**Theorem 2.12** (Symplectic algebras ($w$ odd)). *Let $Q$ be a skew-symmetric bilinear form on a real vector space $V_\mathbb{R}$ of dimension $2m$, and set $G_\mathbb{R} = \mathrm{Aut}(V_\mathbb{R}, Q) \simeq \mathrm{Sp}(2m, \mathbb{R})$ with Lie algebra $\mathfrak{g}_\mathbb{R} = \mathrm{End}(V_\mathbb{R}, Q) \simeq \mathfrak{sp}(2m, \mathbb{R})$. Then the $\mathrm{Ad}(G_\mathbb{R})$–orbits in $\mathrm{Nilp}(\mathfrak{g}_\mathbb{R})$ are indexed by the signed Young diagrams of size $2m$ in which (i) the rows* *of even length are signed, and (ii) the rows of odd length are unsigned and occur with even multiplicity.*

*Example* 2.13. Suppose that $G_{\mathbb{R}} = \mathrm{Sp}(4, \mathbb{R})$. The nilpotent conjugacy classes in $\mathfrak{g}_{\mathbb{R}} = \mathfrak{sp}(4, \mathbb{R})$ are enumerated by the signed Young diagrams

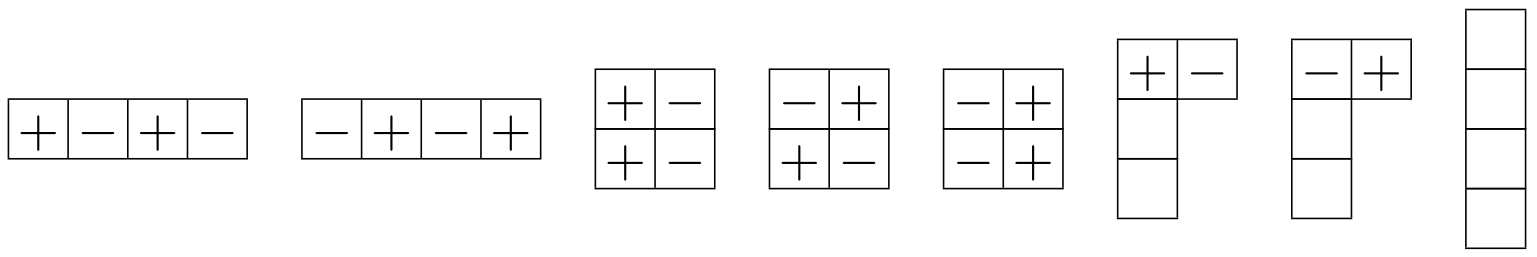

Given a signed Young diagram, let $b^{\pm}$ be the number of boxes labeled with a $\pm$, and let $2b^0$ be the number of unlabeled boxes. The *signature* of a signed Young digram is $\mathrm{sig}\, Y = (s^+, s^-)$ where $s^+ = b^+ + b^0$ and $s^- = b^- + b^0$.

*Example* 2.14. The signed Young diagrams of Theorem 2.12 are all of signature $(m, m)$.

**Theorem 2.15** (Orthogonal algebras ($w$ even))**.** *Let $Q$ be a symmetric bilinear form on a real vector space $V_{\mathbb{R}}$, and set $\mathfrak{g}_{\mathbb{R}} = \mathrm{End}(V_{\mathbb{R}}, Q)$ and $G_{\mathbb{R}} = \mathrm{Aut}(V_{\mathbb{R}}, Q) \simeq \mathrm{O}(a, b)$. The $\mathrm{Ad}(G_{\mathbb{R}})$–orbits in $\mathrm{Nilp}(\mathfrak{g}_{\mathbb{R}})$ are indexed by signed Young diagrams of size $n = \dim V_{\mathbb{R}}$ and signature $(a, b)$ in which (i) the rows of odd length are signed, and (ii) the rows of even length are unsigned and occur with even multiplicity.*

*Remark* 2.16. For the analog of Theorem 2.10(b) with $G^{\circ}_{\mathbb{R}}$ the connected identity component of $\mathrm{SO}(a, b)$, see [6, Theorem 9.3.4]: In the case that the the unsigned Young diagram characterizing the $\mathrm{Ad}(G_{\mathbb{R}})$–orbit $\mathcal{N}$ of $N$ is very even, the orbit $\mathcal{N}$ decomposes into two $\mathrm{Ad}(G^{\circ}_{\mathbb{R}})$–orbits.

*Example* 2.17. Suppose that $G_{\mathbb{R}} = \mathrm{O}(3, 3)$. The nilpotent conjugacy classes in $\mathfrak{g}_{\mathbb{R}} = \mathfrak{so}(3, 3)$ are enumerated by the signed Young diagrams

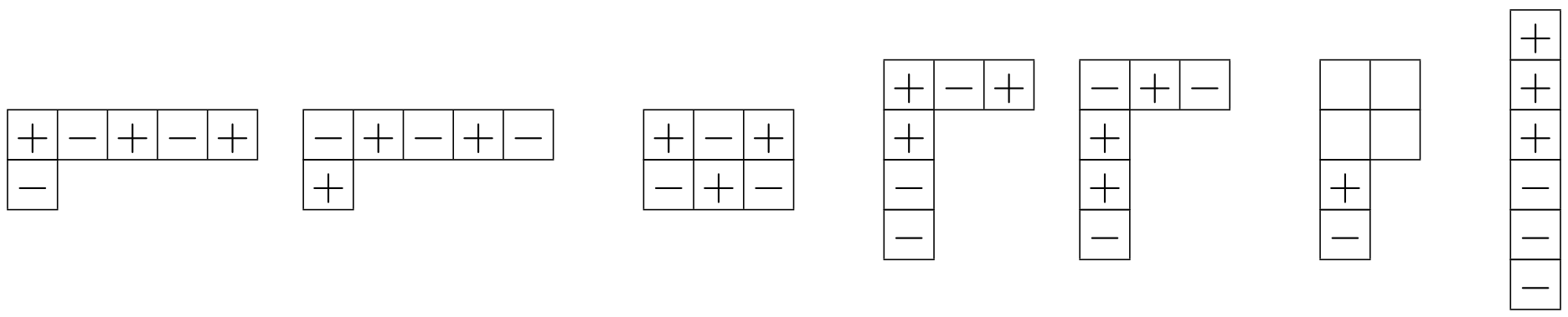

*Example* 2.18. Suppose that $G_\mathbb{R} = \mathrm{O}(4,2)$. The nilpotent conjugacy classes in $\mathfrak{g}_\mathbb{R} = \mathfrak{so}(4,2)$ are enumerated by the signed Young diagrams

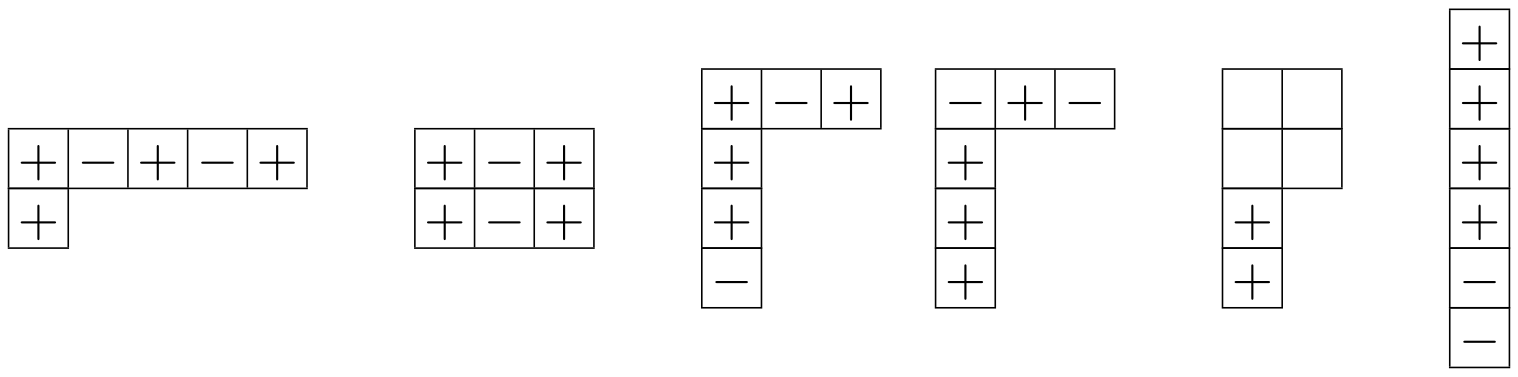

*Example* 2.19. Suppose that $G_\mathbb{R} = \mathrm{O}(5,1)$. The nilpotent conjugacy classes in $\mathfrak{g}_\mathbb{R} = \mathfrak{so}(5,1)$ are enumerated by the signed Young diagrams

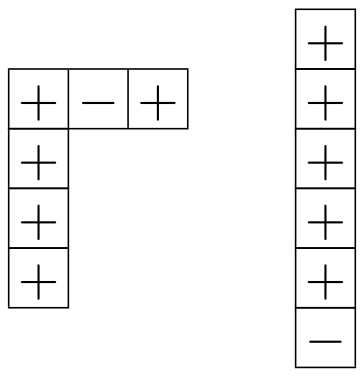

*Remark* 2.20. The Lie algebras $\mathfrak{g}_\mathbb{R}$ of *compact* Lie groups $G_\mathbb{R}$ contain no nilpotent elements other than the trivial $N = 0$.

2.4. **Partial order on conjugacy classes.** Given two nilpotent elements $N_1, N_2 \in \mathfrak{g}_\mathbb{R}$, let $\mathcal{N}_i = \mathrm{Ad}(G_\mathbb{R})N_i$ denote the associated conjugacy classes. We define a partial order on the set of conjugacy classes by

$$\mathcal{N}_1 \le \mathcal{N}_2 \quad \text{if} \quad \mathcal{N}_1 \subset \overline{\mathcal{N}}_2 \,.$$

Đoković's Theorem 2.21 characterizes the partial order in terms of the signed Young diagram classifying the conjugacy classes.

Given a signed Young diagram $Y$, let $Y'$ be the signed Young diagram obtained by removing the last (right-most) box from each row of $Y$. Inductively define $Y^{(k)}$ by $Y^{(0)} = Y$ and $Y^{(k+1)} = (Y^{(k)})'$. Given two signed Young diagrams $Y_1$ and $Y_2$ of signatures $s_1 = (s_1^+, s_1^-)$ and $s_2 = (s_2^+, s_2^-)$, respectively, we write $s_1 \le s_2$ if $s_1^+ \le s_2^+$ and $s_1^- \le s_2^-$. Then we put a partial order the signed Young diagrams by $Y_1 \le Y_2$ if $\mathrm{sig}\, Y_1^{(k)} \le \mathrm{sig}\, Y_2^{(k)}$ for all $k$. The following is [10, Theorem 5].

**Theorem 2.21** (Đoković). *Let $\mathcal{N}_1, \mathcal{N}_2 \subset \mathfrak{g}_\mathbb{R}$ be two nilpotent* $\mathrm{Ad}(G_\mathbb{R})$*–conjugacy classes, and let $Y_1$ and $Y_2$ be the associated signed Young diagrams* (§2.3). *Then $\mathcal{N}_1 \subset \overline{\mathcal{N}}_2$ if and only if $Y_1 \le Y_2$.*

*Remark* 2.22. The analogous partial order on $\mathrm{Ad}(G_{\mathbb{C}})$–conjugacy classes $\mathcal{N}_{\mathbb{C}} \subset \mathrm{Nilp}(\mathfrak{g}_{\mathbb{C}})$ was characterized by Gerstenhaber [11, 12].

*Example* 2.23. If we write $\mathcal{N}' \to \mathcal{N}$ to indicate $\mathcal{N}' < \mathcal{N}$, then the partial order on the nilpotent conjugacy classes of Example 2.13 is given by with the remaining relations given by transitivity.

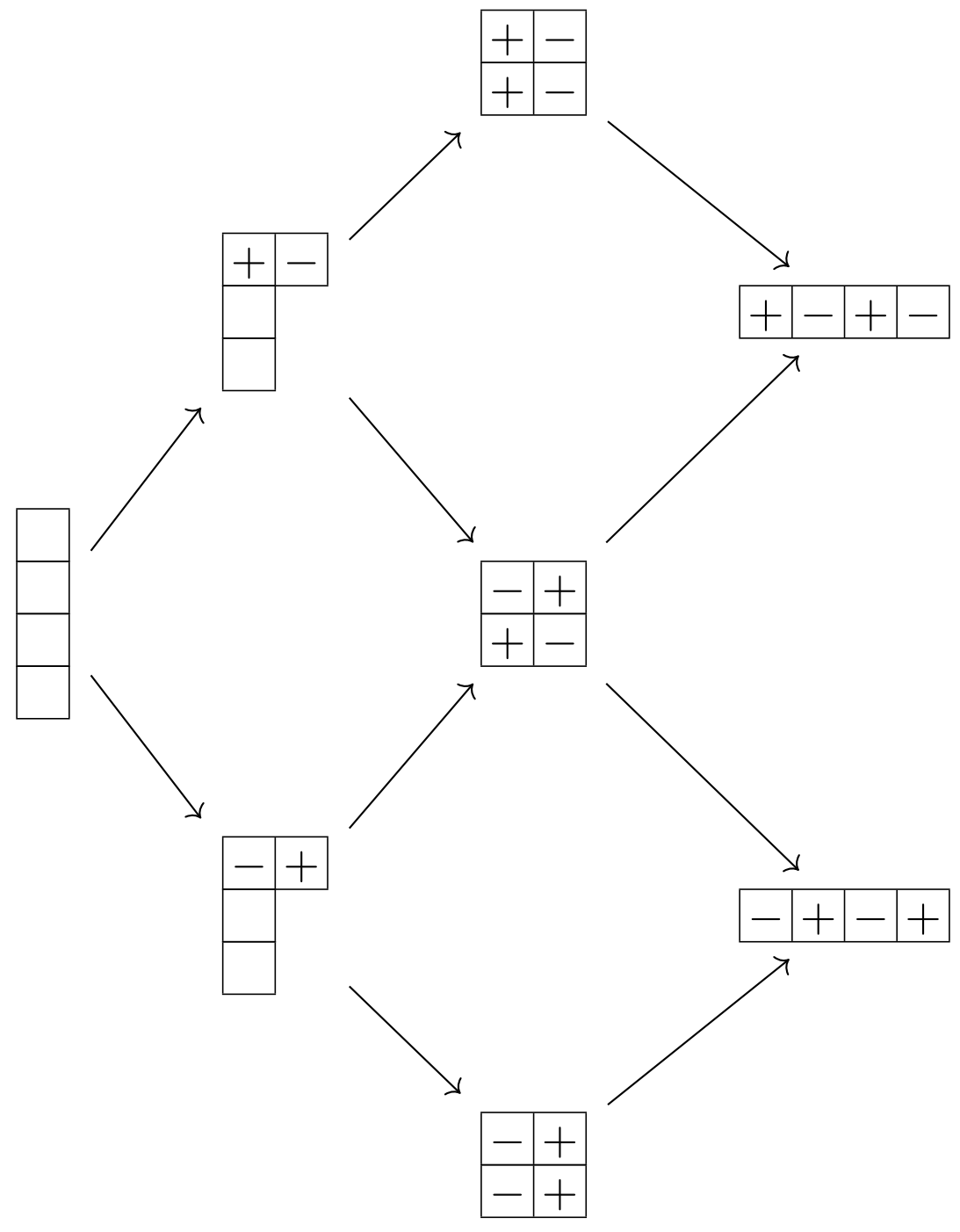

*Example* 2.24. Likewise the partial order on the nilpotent conjugacy classes of Example 2.17 is given by with the remaining relations given by transitivity.

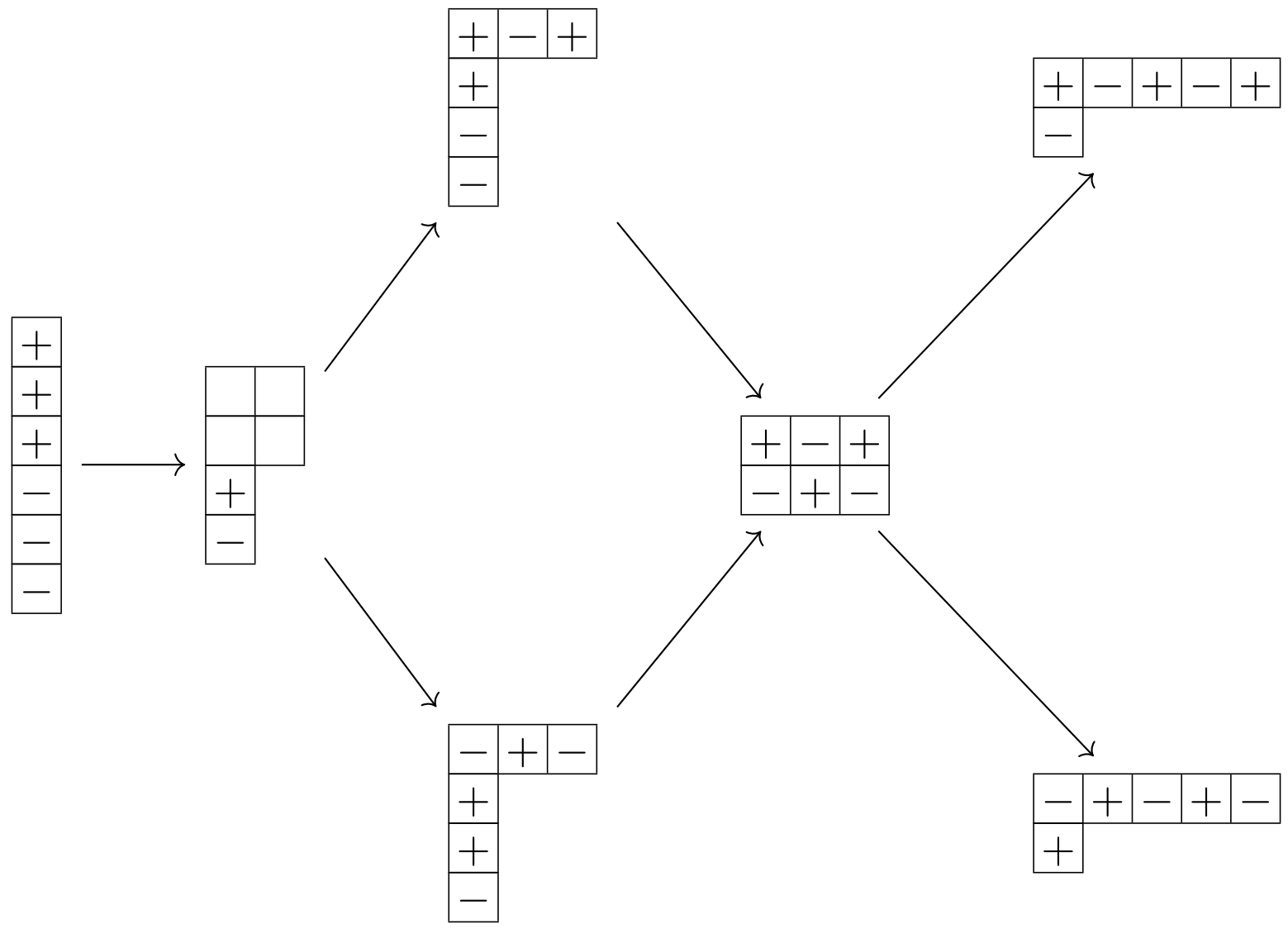

*Remark* 2.25. In the event that $\mathcal{N}$ decomposes into two $\mathrm{Ad}(G^{\circ}_{\mathbb{R}})$–orbits $\mathcal{N}_1$ and $\mathcal{N}_2$ (cf. Remark 2.16), the orbits are not comparable; that is, $\mathcal{N}_i \not\subset \overline{\mathcal{N}}_j$, for $i \neq j$.

2.5. **Polarized mixed Hodge structures and signed Young diagrams.** Let $D = G_{\mathbb{R}}/G^0_{\mathbb{R}}$ be a *period domain* parameterizing weight $w$, $Q$–polarized Hodge structures, so that

$$G_{\mathbb{R}} \;=\; \mathrm{Aut}(V_{\mathbb{R}}, Q)\,.$$

Let $(F^{\bullet}, N)$ be an $\mathbb{R}$–split nilpotent orbit on $D$, and let

$$\mathcal{N} \;:=\; \mathrm{Ad}(G_{\mathbb{R}})\cdot N$$

be the corresponding conjugacy class.

Let $W_{\bullet}(N)$ be the monodromy filtration of $N$; then $(F^{\bullet}, W_{\bullet}(N))$ is a polarized mixed Hodge structure. Let

$$V_{\mathbb{C}} \;=\; \bigoplus V^{p,q}$$

be the associated Deligne bigrading. Without loss of generality we assume that the $(F^{\bullet}, W_{\bullet}(N))$ is $\mathbb{R}$–split. The associated *Hodge diamond* is the configuration of points in the $pq$–plane for which $V^{p,q} \neq 0$. In this section we explain how to construct the signed Young diagram indexing $\mathcal{N}$ from the Hodge diamond. (This, along with Đoković's Theorem 2.21, will give constraints on the degenerations associated with the faces of a nilpotent cone $\sigma \ni N$ underlying a nilpotent orbit on $D$. See §5 for an illustration.)

Fix $p, q$ and define $\ell$ by $w + \ell = p + q$. Suppose $\ell \geq 0$ and let

$$P^{p,q} \;:=\; \ker\{N^\ell : V^{p,q} \to V^{p-\ell,q-\ell}\}$$

be the $N$–*primitive subspace*. By our hypothesis that the polarized mixed Hodge structure is $\mathbb{R}$–split, we have $\overline{P^{p,q}} = P^{q,p}$. Moreover,

$$Q_\ell(\cdot,\cdot) \;:=\; Q(\cdot, N^\ell \cdot)$$

is a nondegenerate bilinear form on $(P^{p,q} + P^{q,p}) \cap V_\mathbb{R}$ satisfying the symmetry

$$Q_\ell(u,v) \;=\; (-1)^{w+\ell} Q_\ell(v,u)\,.$$

First suppose that $p = q$. Then $P^{p,q}$ is real and admits a basis of $Q_\ell$–orthogonal real vectors. Given one such basis vector $v \in V_\mathbb{R}$,

$$v\,,\ Nv\,,\ \cdots\,,\ N^\ell v$$

is an $N$–string, and the polarization conditions assert

$$0 \;<\; Q(v, N^\ell v) \;=\; Q_\ell(v,v)\,.$$

In the lexicon of Đoković's [10], the "isomorphism class" of this $N$–string is the *rank* $\ell + 1$ *gene*

$$g^+(\ell+1)\,, \quad \begin{cases} \text{if } w \text{ is even and } \ell \equiv 0 \text{ mod } 4, \text{ or} \\ \text{if } w \text{ is odd and } \ell \equiv 1 \text{ mod } 4; \end{cases}$$

$$g^-(\ell+1)\,, \quad \begin{cases} \text{if } w \text{ is even and } \ell \equiv 2 \text{ mod } 4, \text{ or} \\ \text{if } w \text{ is odd and } \ell \equiv 3 \text{ mod } 4. \end{cases}$$

Pictorially the *positive* gene $g^+(\ell+1)$ is identified with a row of $\ell+1$ boxes, labeled with alternating signs and beginning with $+$; that is, $g^+(\ell+1)$ is visualized as $\boxed{+}\boxed{-}\boxed{+}\boxed{-}\cdots$. The *negative* gene $g^-(\ell+1)$ is depicted as $\boxed{-}\boxed{+}\boxed{-}\boxed{+}\cdots$.

Next suppose that $p \neq q$. Fix $\xi = u + \mathbf{i}v \in P^{p,q}$, with $u, v \in V_{\mathbb{R}}$. The polarization conditions assert that $Q_\ell(\xi, \xi) = 0$; equivalently,

$$Q_\ell(u,u) \;=\; Q_\ell(v,v)\,,$$

and

$$Q_\ell(u,v) \;=\; 0 \quad \text{if } w+\ell \text{ is even.}$$

The polarization conditions also impose $\mathbf{i}^{p-q} Q_\ell(\xi, \bar\xi) > 0$ for all $0 \neq \xi$. Equivalently we have the following:

$$\begin{array}{lllll}
\text{If } p-q \equiv 0 \bmod 4, & \text{then } 0 < Q_\ell(u,u)\,; & & \\
\text{if } p-q \equiv 2 \bmod 4, & \text{then } 0 > Q_\ell(u,u)\,; & & \\
\text{if } p-q \equiv 1 \bmod 4, & \text{then } 0 = Q_\ell(u,u) & \text{and} & 0 < Q_\ell(u,v)\,; \\
\text{if } p-q \equiv 3 \bmod 4, & \text{then } 0 = Q_\ell(u,u) & \text{and} & 0 > Q_\ell(u,v)\,.
\end{array}$$

(Note that $p-q$ is even if and only if $p+q = w+\ell$ is even.) Again, in the language of [10], a *chromosome* is a formal linear combination of genes with non-negative integral coefficients. (A chromosome is just a signed Young diagram.) The two $N$–strings $\{u, \ldots, N^\ell u\}$ and $\{v, \ldots, N^\ell v\}$ in $V_{\mathbb{R}}$ correspond to the chromosome

$$\begin{array}{ll}
2g(\ell+1)\,, & \text{if } w+\ell \text{ is odd;} \\
2g^+(\ell+1)\,, & \left\{ \begin{array}{l} \text{if } w \text{ is even and } \ell, p-q \equiv 0 \bmod 4, \text{ or} \\ \text{if } w \text{ is even and } \ell, p-q \not\equiv 0 \bmod 4, \text{ or} \\ \text{if } w \text{ is odd, } \ell \equiv 1 \bmod 4 \text{ and } p-q \equiv 0 \bmod 4, \text{ or} \\ \text{if } w \text{ is odd, } \ell \equiv 3 \bmod 4 \text{ and } p-q \not\equiv 0 \bmod 4; \end{array} \right. \\
2g^-(\ell+1)\,, & \text{otherwise.}
\end{array}$$

The *unpolarized* gene $g(\ell+1)$ is indicated by a row of $\ell+1$ boxes, without labels.

*Definition* 2.26. The *partially signed Young diagram* $Y(F^\bullet, N)$ (or chromosome) *associated with the* $\mathbb{R}$*–split polarized mixed Hodge structure* $(F^\bullet, N)$ is the union of genes obtained from the $N$–string decomposition of the standard representation $V_{\mathbb{R}}$.

Note that $G_{\mathbb{R}}$ acts on $\mathbb{R}$–split polarized mixed Hodge structures by $g \cdot (F^\bullet, N) = (gF^\bullet, \mathrm{Ad}_g N)$, for $g \in G_{\mathbb{R}}$. It follows from the classification results of §2.3 that $Y(F^\bullet, N)$ depends only on the $G_{\mathbb{R}}$–conjugacy class of $(F^\bullet, N)$.

*Example* 2.27 (Period domain for $\mathbf{h} = (3,3,3)$). This example was studied by Cattani and Kaplan in [5, §4]. We have $G_{\mathbb{R}} = \mathrm{O}(3,6)$, and there are five conjugacy classes of $\mathbb{R}$–split PMHS.

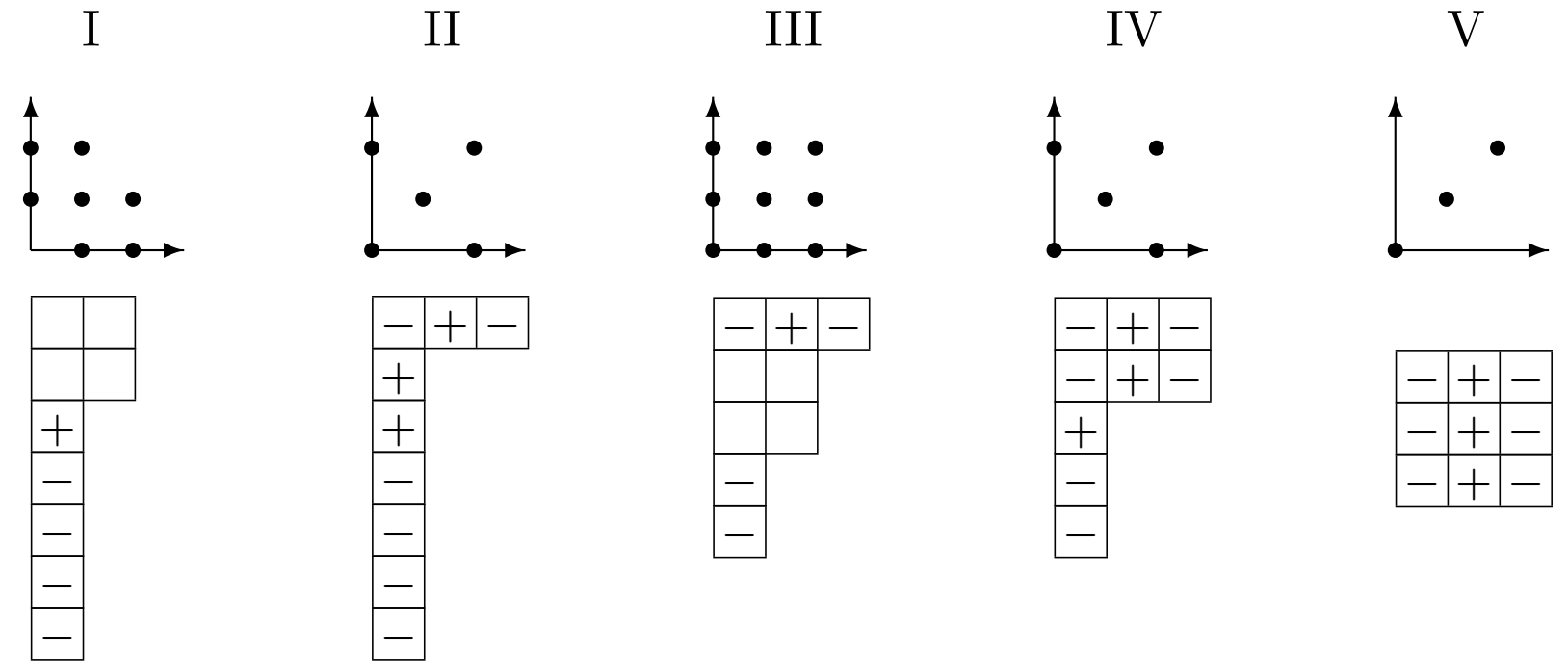

Đoković's Theorem 2.21 yields

$$\mathcal{N}_{\rm I} \ < \ \mathcal{N}_{\rm II} \ < \ \mathcal{N}_{\rm III} \ < \ \mathcal{N}_{\rm IV} \ < \ \mathcal{N}_{\rm IV} \,.$$

*Example* 2.28 (Period domain for $\mathbf{h} = (1,1,1,1,1,1)$). We have $G_{\mathbb{R}} = \mathrm{Sp}(3,\mathbb{R})$ there are seven (conjugacy classes of) $\mathbb{R}$–split PMHS.

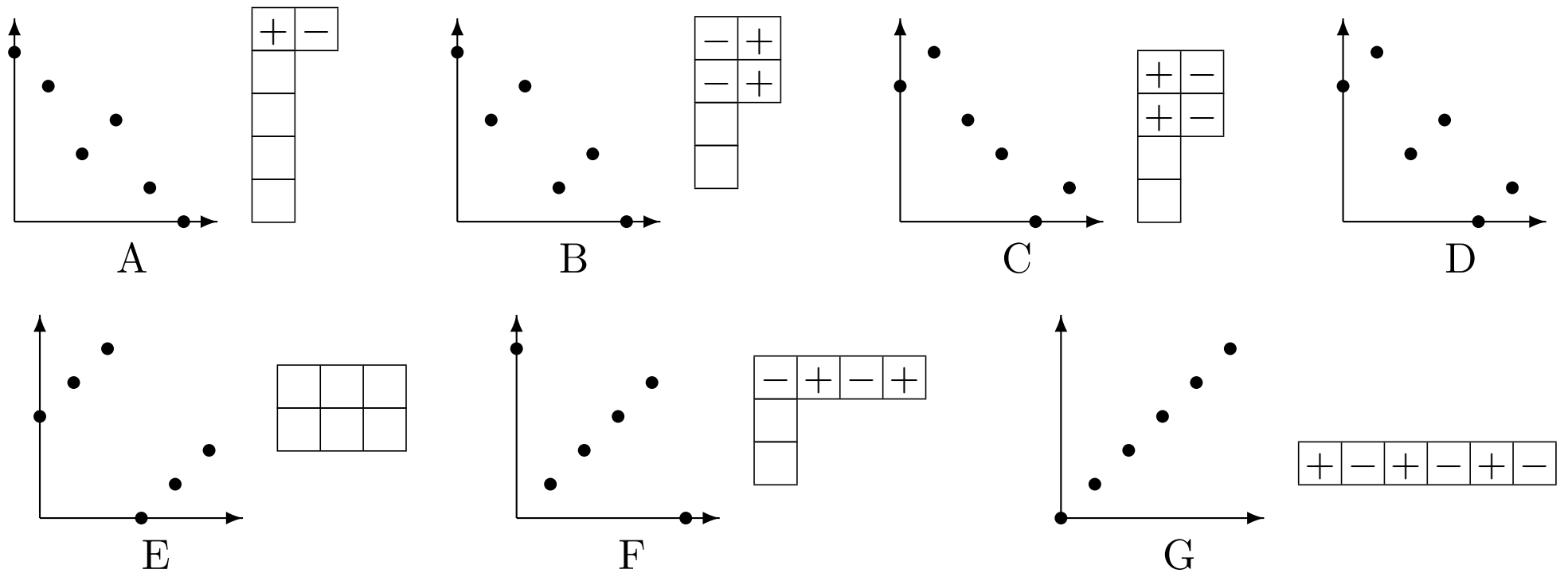

From Theorem 2.21 we find

$$\mathcal{N}_{\rm A} < \mathcal{N}_{\rm C}\,,\ \mathcal{N}_{\rm F}\,; \quad \mathcal{N}_{\rm B} < \mathcal{N}_{\rm E}\,,\ \mathcal{N}_{\rm F}\,; \quad \mathcal{N}_{\rm C} < \mathcal{N}_{\rm D}\,,\ \mathcal{N}_{\rm E}\,; \quad \mathcal{N}_{\rm D}\,,\ \mathcal{N}_{\rm E}\,,\ \mathcal{N}_{\rm F} < \mathcal{N}_{\rm G}\,.$$

where the relations obtained by the transitivity of the partial order are omitted.

## 3. Nilpotent cones

The goal of this section is to describe one approach to identifying the nilpotent cones that underlying nilpotent orbits on a period domain (or, more generally, a Mumford–Tate domain [13]). We begin in §3.1 by reviewing the definition of nilpotent orbits; in §3.1 we outline the strategy. The approach will be worked out for weight two period domains in §4, and illustrated in the case of Hodge numbers $\mathbf{h} = (2, m, 2)$ in §5.

3.1. **Nilpotent orbits.** Let $V_\mathbb{R}$ be a real vector space of dimension $n$ with a $\mathbb{Q}$–structure defined by a lattice $V_\mathbb{Z} \subset V_\mathbb{R}$. Fix $w \in \mathbb{Z}$ and let $Q$ be a nondegenerate $(-1)^w$–symmetric bilinear form on $V_\mathbb{R}$ defined over $\mathbb{Q}$. Fix Hodge numbers $\mathbf{h} = \{h^{p,q} \mid p+q = w\}$. Let $D$ denote the period domain parameterizing $Q$–polarized Hodge structures on $V_\mathbb{Q}$ with Hodge numbers $\mathbf{h}$. Let $\check{D}$ denote the compact dual of $D$. A ($m$–variable) *nilpotent orbit* on $D$ consists of a pair $(F^\bullet; N_1, \ldots, N_m)$ such that $F^\bullet \in \check{D}$, the $N_i \in \mathfrak{g}_\mathbb{R}$ are commuting nilpotents and $N_i F^p \subset F^{p-1}$, and the holomorphic map $\psi : \mathbb{C}^m \to \check{D}$ defined by

$$\psi(z^1, \ldots, z^m) \ = \ \exp(z^i N_i) F^\bullet \tag{3.1}$$

has the property that $\psi(z) \in D$ for $\mathrm{I}(z^i) \gg 0$. The associated (open) *nilpotent cone* is

$$\sigma \ = \ \{t^i N_i \mid t^i > 0\}\,. \tag{3.2}$$

Recall that the monodromy filtration $W_\bullet(N)$ is independent of our choice of $N \in \sigma$; so $W_\bullet(\sigma)$ is well-defined. Given a nilpotent orbit we will assume, without loss of generality, that the polarized mixed Hodge structure $(F^\bullet, W_\bullet(\sigma))$ is $\mathbb{R}$–split. Let

$$V_\mathbb{C} \ = \ \oplus V^{p,q} \quad \text{and} \quad \mathfrak{g}_\mathbb{C} \ = \ \oplus \mathfrak{g}^{p,q}$$

denote the Deligne bigradings. Recall that

$$\sigma \ \subset \ \mathfrak{g}_\mathbb{R}^{-1,-1}\,.$$

3.2. **Identification of cones underlying nilpotent orbits.** Observe that

$$\mathfrak{m}_\mathbb{C} \ := \ \oplus_p \, \mathfrak{g}^{p,p} \ \subset \ \mathfrak{g}_\mathbb{C}$$

is the subalgebra of $\mathfrak{g}_\mathbb{C}$ preserving the subspaces

$$V_m \ := \ \bigoplus_{q-p=m} V^{p,q} \ \subset \ V_\mathbb{C}\,, \quad m \in \mathbb{Z}. \tag{3.3}$$

Visually, $\mathfrak{m}$ corresponds to the dots on the diagonal $p = q$ in the Hodge diamond. Because the PMHS is $\mathbb{R}$–split, the subalgebra $\mathfrak{m}$ are defined over $\mathbb{R}$. Let $M^0_\mathbb{R} \subset G_\mathbb{R}$ be the *connected* Lie subgroup with (Levi) Lie algebra

$$\mathfrak{m}^0_\mathbb{R} \;:=\; \mathfrak{g}^{0,0}_\mathbb{R} .$$

It will be convenient to note that

(3.4) *the real form* $\mathfrak{m}^0_\mathbb{R} = \mathfrak{g}^{0,0}_\mathbb{R}$ *is the subalgebra of* $\mathfrak{g}_\mathbb{R}$ *preserving the* $V^{p,q}$.

**Lemma 3.5.** *Let* $(F^\bullet, N)$ *be an* $\mathbb{R}$*–split nilpotent orbit on* $D = G_\mathbb{R}/G^0_\mathbb{R}$ *with Deligne bigrading* $\mathfrak{g}_\mathbb{C} = \oplus \mathfrak{g}^{p,q}$. *Let* $M^0_\mathbb{R} \subset M_\mathbb{R}$ *be as defined above.*

(a) *The orbit*

$$\mathcal{N}^0 \;:=\; \mathrm{Ad}(M^0_\mathbb{R}) \cdot N$$

*is open in* $\mathfrak{g}^{-1,-1}_\mathbb{R}$.[2]

(b) *Suppose that* $\mathcal{W}^\circ_N$ *is the connected component of*

$$\mathcal{W}_N \;:=\; \{N' \in \mathfrak{g}^{-1,-1}_\mathbb{R} \mid W(N) = W(N')\}$$

*containing* $N$. *Then* $\mathcal{W}^\circ_N = \mathcal{N}^0$.

The lemma is well-known; the proof is included in the appendix for completeness. The key points are that (i) $\mathcal{W}_N$ is preserved under the adjoint action of $M^0_\mathbb{R}$, and (ii) the orbit $\mathrm{Ad}(M^0_\mathbb{R}) \cdot N'$ is open in $\mathfrak{g}^{-1,-1}_\mathbb{R}$ for every $N' \in \mathcal{W}_N$. Thus $\mathcal{W}_N$ is a disjoint union of open $\mathrm{Ad}(M^0_\mathbb{R})$–orbits.

From Lemma 3.5(b) we obtain:

**Corollary 3.6.** *Assume the hypotheses of Lemma 3.5. If a nilpotent cone* $\sigma$ *containing* $N$ *underlies a nilpotent orbit, then*

$$\sigma \;\subset\; \mathcal{N}^0 \;=\; \mathrm{Ad}(M^0_\mathbb{R}) \cdot N \;\subset\; \mathfrak{g}^{-1,-1} \quad \text{and} \quad N_i \in \overline{\mathcal{N}}^0 .$$

*Remark* 3.7. Together Theorem 2.21 and Corollary 3.6 gives us representation theoretic constraints on the degenerations associated with the faces of a nilpotent cone underlying a nilpotent orbit. See §5.2 for an illustration.

[2]The superscript of 0 in $\mathcal{N}^0$ is meant to distinguish these nilpotent $\mathrm{Ad}(M^0_\mathbb{R})$–conjugacy classes from the $\mathrm{Ad}(G_\mathbb{R})$–conjugacy classes $\mathcal{N}$ of §2.3.

For the converse to Corollary 3.6, recall Cattani and Kaplan's [3, Theorem 2.3]

**Theorem 3.8** (Cattani–Kaplan)**.** *Fix $F^\bullet \in \check{D}$ and a nilpotent cone* (3.2) *with the properties that:*

(i) *$N_i F^p \subset F^{p-1}$ for every $i$;*

(ii) *$N_i^{k+1} = 0$, where $k$ is the level of the Hodge structures on $D$; and*

(iii) *the filtration $W_\bullet(N)$ does not depend on the choice of $N \in \sigma$.*

*Then $(F^\bullet; N)$ is a limiting mixed Hodge structure for* some *$N \in \sigma$, if and only if $(F^\bullet; N_1, \ldots, N_m)$ is an $m$–variable nilpotent orbit.*

This, along with Lemma 3.5(b), yields the converse to Corollary 3.6:

**Proposition 3.9.** *Given an $\mathbb{R}$–split PMHS $(F^\bullet, N)$ on $D$, let $M^0_{\mathbb{R}}$ and $\mathcal{N}^0$ be as defined above. If $\sigma \subset \mathcal{N}^0$ is a nilpotent cone, then $\sigma$ underlies a nilpotent orbit at $F^\bullet$.*

The upshot of this discussion is

*Remark* 3.10. The cones $\sigma$ underlying a multivariable nilpotent orbit on a domain $D$ may be identified as follows. Begin with an $\mathbb{R}$–split PMHS $(F^\bullet, W_\bullet(N))$ on $D$. The Deligne bigrading determines the diagonal Levi subgroup $M$. Any nilpotent cone $\sigma \subset \mathcal{N}^0$ will underlie a nilpotent orbit, and all such cones arise in this fashion. So to identify the nilpotent cones underlying a nilpotent orbit we must have a good enough/explicit enough geometric description of $\mathcal{N}^0$ to understand how the nilpotent cones can "fit" inside. So the strategy proceeds in three steps:

*Step 1:* Enumerate the $\mathbb{R}$–split PMHS $(F^\bullet, N)$ on $D$. For an arbitrary Mumford–Tate domain $D = G_{\mathbb{R}}/G^0_{\mathbb{R}}$, with $G_{\mathbb{R}}$ connected, these are given (up to the action of $G_{\mathbb{R}}$) by [21]. As Cattani has pointed out, if $D$ is a period domain and $G_{\mathbb{R}}$ is the full automorphism group $\mathrm{Aut}(V_{\mathbb{R}}, Q)$, then the $\mathbb{R}$–split PMHS are enumerated (again, up to the action of $G_{\mathbb{R}}$) by the Hodge diamonds.

*Step 2:* Determine the $\mathrm{Ad}(M^0_{\mathbb{R}})$–orbit $\mathcal{N}^0$ of $N \in \mathfrak{g}_{\mathbb{R}}^{-1,-1}$. This will require a good description of $\mathfrak{g}_{\mathbb{R}}^{-1,-1}$ as a $M^0_{\mathbb{R}}$–module. In the case of weight two period domains this description is given in §4.

Also, in the context of "understanding the cones," it should be noted that a theorem of Đoković's can help us understand, given a nilpotent on a period domain, constraints on the degenerations coming from the faces of the cone; this is discussed in §5.2.

3.3. **The CKS commuting** SL(2)**'s.** To a nilpotent cone $\sigma$ underlying a nilpotent orbit on a Hodge domain $D$, Cattani, Kaplan and Schmid [4] associate a set of commuting $\mathfrak{sl}(2)$'s contained in $\mathfrak{m}^{\mathrm{ss}}_{\mathbb{R}}$.

Given a Cartan decomposition $\mathfrak{m}^{\mathrm{ss}}_{\mathbb{R}} = \mathfrak{k} \oplus \mathfrak{k}^{\perp}$, the *real rank* of $\mathfrak{m}^{\mathrm{ss}}_{\mathbb{R}}$ is the dimension of a maximal subspace of $\mathfrak{k}^{\perp}$ consisting of commuting semisimple elements.

**Lemma 3.11** (Nilpotent cones versus commuting $\mathfrak{sl}(2)$'s)**.** *The number of (nontrivial) commuting $\mathfrak{sl}(2)$'s is bounded by the real rank of the semisimple factor $\mathfrak{m}^{\mathrm{ss}}_{\mathbb{R}} = [\mathfrak{m}_{\mathbb{R}}, \mathfrak{m}_{\mathbb{R}}]$.*

*Proof.* To see this, let $Y_1, \ldots, Y_s$ denote the neutral elements of the commuting $\mathfrak{sl}(2)$'s. They are linearly independent and so span an $s$–dimensional abelian subspace of $\mathfrak{m}_{\mathbb{R}}$ consisting of semisimple elements. It remains to show that we may choose a Cartan decomposition $\mathfrak{m}_{\mathbb{R}} = \mathfrak{k} \oplus \mathfrak{k}^{\perp}$ so that $Y_i \in \mathfrak{k}^{\perp}$; for then $s \leq \mathrm{rank}_{\mathbb{R}}\, \mathfrak{m}^{\mathrm{ss}}_{\mathbb{R}}$. To see this, observe that the polarization conditions on the limit mixed Hodge structure imply that $Q(Y_i, Y_i) < 0$. Since the polarization $Q$ is minus the Killing form, it follows that the desired Cartan decomposition exists. □

3.4. **Two special cases.** We finish §3 with discussions of two special cases that we will encounter in these notes.

3.4.1. *The Hermitian case.* We say that a simple factor $\mathfrak{m}' \subset \mathfrak{m}$ is *Hermitian* if any of the following equivalent conditions hold:[3]

(a) $\mathfrak{m}' = (\mathfrak{m}' \cap \mathfrak{g}^{1,1}) \oplus (\mathfrak{m}' \cap \mathfrak{g}^{0,0}) \oplus (\mathfrak{m}' \cap \mathfrak{g}^{-1,-1})$;

(b) $\mathfrak{m}' \cap \mathfrak{g}^{p,p} = 0$ when $|p| \geq 2$;

(c) $\mathfrak{m}' \cap \mathfrak{g}^{-1,-1}$ is abelian.

When $\mathfrak{m}$ is Hermitian there will exist nilpotent cones $N \in \sigma \subset \mathcal{N}^0$ underlying nilpotent orbits that are open in $\mathfrak{g}^{-1,-1}_{\mathfrak{g},\mathbb{R}}$.

---

[3]Equivalence requires that we assume the horizontal distribution is bracket–generating.

3.4.2. *The contact case.* We say that a simple factor $\mathfrak{m}' \subset \mathfrak{m}$ is *contact* if $\dim \mathfrak{m}' \cap \mathfrak{g}^{2,2} = 1$ and $\mathfrak{m}' \cap \mathfrak{g}^{p,p} = 0$ for all $|p| \geq 3$. In this case the maximal abelian subspaces of $\mathfrak{g}^{-1,-1}$ are the Lagrangian subspaces $\Lambda$ of a symplectic form $\nu$ on $\mathfrak{g}^{-1,-1}$ that is invariant under the reductive (Levi) subalgebra $\mathfrak{g}^{0,0}$. The symplectic form is defined (up to scale) by choosing a nonzero $z$ in the one–dimensional $\mathfrak{g}_{\mathbb{R}}^{-2,-2}$ and setting $[x, y] =: \nu(x, y)z$ for any $x, y \in \mathfrak{g}^{-1,-1}$.

## 4. Weight two period domains

Suppose that $D$ is a period domain parameterizing effective weight two Hodge structures (Hodge numbers $\mathbf{h} = (h^{2,0}, h^{1,1}, h^{0,2})$). Our goal in this section is to address Steps 2 and 3 in the strategy (outlined in Remark 3.10) to identify the nilpotent cones underlying nilpotent orbits on $D$. Section 4.1 provides the necessary descriptions of both $M^0_{\mathbb{R}}$ and $\mathfrak{g}_{\mathbb{R}}^{-1,-1}$ as an $M^0_{\mathbb{R}}$–representation to explicitly describe the orbits $\mathcal{N}^0$. Section 4.2 describes a decomposition of the orbits $\mathcal{N}^0$ into simpler objects, culminating in an explicit description of $\mathcal{N}^0$ by Proposition 4.17.

4.1. **The representation theory.** We assume given an $\mathbb{R}$–split nilpotent orbit $(F, N)$ on $D$. The Hodge diamond of the Deligne bigrading $V_{\mathbb{C}} = \oplus V^{p,q}$ is contained in

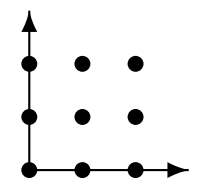

In particular,

$$V_0 \;=\; V^{0,0} \oplus V^{1,1} \oplus V^{2,2}\,, \quad V_1 \;=\; V^{1,0} \oplus V^{2,1} \quad \text{and} \quad V_2 = V^{2,0}\,.$$

Let

$$a \;:=\; \dim V^{2,2}\,, \quad a+b \;:=\; \dim V^{1,1}\,, \quad c \;:=\; \dim V^{1,2} \quad \text{and} \quad d \;:=\; \dim V^{0,2}\,.$$

Then

$$h^{2,0} \;=\; a+c+d \quad \text{and} \quad h^{1,1} \;=\; a+b+2c\,.$$

Let

$$M_{\mathbb{R}} \;=\; \{g \in G_{\mathbb{R}} \mid g(V_m) \subset V_m\,,\ \forall m\} \tag{4.1}$$

be the subgroup of $G_{\mathbb{R}}$ preserving the subspaces (3.3), and note that $\mathfrak{m}_{\mathbb{R}}$ is the Lie algebra of $M_{\mathbb{R}}$. In an abuse of notation, we let $G_{\mathbb{R}} \cap \mathrm{Aut}(V_m)$ denote the subgroup of elements $g \in G_{\mathbb{R}}$ that preserve $V_m$ (and therefore also preserve $V_{-m}$) and act trivially on $V_\ell$ for all $\ell \neq \pm m$.

**Proposition 4.2.** *Let $(F^\bullet, N)$ be an $\mathbb{R}$–split nilpotent orbit on a period domain parameterizing effective weight two Hodge structures. The diagonal subgroup is*

$$\begin{array}{rcl} M_{\mathbb{R}} & = & (G_{\mathbb{R}} \cap \mathrm{Aut}(V_0)) \ \times\ (G_{\mathbb{R}} \cap \mathrm{Aut}(V_1)) \ \times\ (G_{\mathbb{R}} \cap \mathrm{Aut}(V_2)) \, . \\ & \simeq & \mathrm{O}(a+b, 2a) \ \times\ \mathrm{Sp}(2c, \mathbb{R}) \ \times\ \mathrm{O}(d, \mathbb{R}) \, . \end{array} \tag{4.3}$$

*Remark* 4.4. For the groups $M_{\mathbb{R}}$ of Proposition 4.2, the subalgebra $\mathfrak{m}^{\mathrm{ss}}_{\mathbb{R}}$ has real rank $c + \min\{2a, a+b\}$, cf. [19].

*Proof.* From (4.1) and $\overline{V_m} = V_{-m}$ we see that

$$M_{\mathbb{R}} \ = \ G_{\mathbb{R}} \ \cap \ \{\mathrm{Aut}(V_0) \ \times \ \mathrm{Aut}(V_1) \ \times \ \mathrm{Aut}(V_2)\} \, .$$

From the facts that (i) $G_{\mathbb{R}} = \mathrm{Aut}(V_{\mathbb{R}}, Q)$ and (ii) the $V_m + V_{-m}$, $m = 0, 1, 2$, are all pairwise orthogonal we see that the first equation of (4.3) holds. It remains to identify the factors $G_{\mathbb{R}} \cap \mathrm{Aut}(V_m)$, $m = 0, 1, 2$. As $V_0$ is defined over $\mathbb{R}$, and $Q|_{V_{0,\mathbb{R}}}$ is nondegenerate of signature $(a+b, 2a)$, it follows directly that

$$G_{\mathbb{R}} \cap \mathrm{Aut}(V_0) \ = \ \mathrm{Aut}(V_{0,\mathbb{R}}, Q) \ \simeq \ \mathrm{O}(a+b, 2a) \, .$$

Next consider $G_{\mathbb{R}} \cap \mathrm{Aut}(V_2)$. From the polarization condition $0 < -Q(u, \overline{u})$ for all $0 \neq u \in V^{2,0} = V_2$, we see that we may pick a basis $\{z_s = x_s + \mathbf{i} y_s \mid s = 1, \ldots, d\}$ so that $Q(z_s, \bar{z}_t) = -\delta_{st}$ and $x_s, y_s \in V_{\mathbb{R}}$ for all $s, t$. Next observe that:

(i) From $-\delta_{st} = Q(z_s, \bar{z}_t)$ and $0 = Q(z_s, z_t)$ we may deduce that $-\frac{1}{2}\delta_{st} = Q(x_s, x_t) = Q(y_s, y_t)$ and $0 = Q(x_s, y_t)$ for all $s, t$.

(ii) Moreover, any element $g \in G_{\mathbb{R}}$ preserving $V_2$ must preserve both $\mathrm{span}_{\mathbb{R}}\{x_s \mid s = 1, \ldots, d\}$ and $\mathrm{span}_{\mathbb{R}}\{y_s \mid s = 1, \ldots, d\}$. In fact, the action of $g$ on $\mathrm{span}\{y_s\}_{s=1}^d$ is determined by the action of $g$ on $\mathrm{span}\{x_s\}_{s=1}^d$: if $g_s^r \in \mathbb{R}$ is defined by $g(x_s) = g_s^r x_r$, then $g(y_s) = g_s^r y_r$.

From these two observations we deduce that

$$G_{\mathbb{R}} \cap \mathrm{Aut}(V_2) \ = \ \mathrm{O}(d, \mathbb{R}) \, .$$

It remains to show that

$$G_{\mathbb{R}} \cap \mathrm{Aut}(V_1) \ = \ \mathrm{Sp}(2c, \mathbb{R}) \, .$$

The polarization condition $\mathbf{i}Q(u, N\bar{u})$ for all $0 \neq u \in V^{2,1}$ implies that we may pick a basis $\{w_s = u_s + \mathbf{i}v_s \mid s = 1, \ldots, c\}$ of $V^{2,1}$ so that $\mathbf{i}Q(w_s, N\overline{w}_t) = \delta_{st}$ and $u_s, v_s \in V_{\mathbb{R}}$ for all $s, t$.

(iii) Note that $\{\overline{w}_s\}_{s=1}^c$, $\{Nw_s\}_{s=1}^c$ and $\{N\overline{w}_s\}_{s=1}^c$ are bases of $V^{1,2}$, $V^{1,0}$ and $V^{0,1}$, respectively. Moreover, the fact that $Q(V^{p,q}, V^{r,s}) = 0$ if $(p,q) \neq (s,r)$ implies that

$$\begin{aligned}
0 &= Q(u_s, u_t) = Q(u_s, v_t) = Q(v_s, v_t)\,, \\
0 &= Q(Nu_s, Nu_t) = Q(Nu_s, Nv_t) = Q(Nv_s, Nv_t)\,, \\
0 &= Q(u_s, Nu_t) = Q(v_s, Nv_t)\,, \\
\tfrac{1}{2}\delta_{st} &= Q(u_s, Nv_t) = -Q(v_s, Nu_t)\,.
\end{aligned}$$

(iv) Any element $g \in G_{\mathbb{R}}$ preserving $V_1$ must preserve $\mathrm{span}_{\mathbb{R}}\{u_s, Nu_s \mid s = 1, \ldots, c\}$ and $\mathrm{span}_{\mathbb{R}}\{v_s, Nv_s \mid s = 1, \ldots, c\}$. Moreover, the action of $g$ on $\mathrm{span}\{v_s, N_s\}_{s=1}^c$ is determined by the action of $g$ on $\mathrm{span}\{u_s, Nu_s\}_{s=1}^c$: if $A_s^r, B_s^r, C_s^r, D_s^r \in \mathbb{R}$ are defined by $g(u_s) = A_s^r u_r + B_s^r Nu_r$ and $g(Nu_s) = C_s^r u_r + D_s^r Nu_r$, then $g(v_s) = A_s^r v_r + B_s^r Nv_r$ and $g(Nv_s) = C_s^r v_r + D_s^r Nv_r$.

From these two observations we see that $g \in G_{\mathbb{R}} \cap \mathrm{Aut}(V_1)$ if and only if the transformation represented by the $2c \times 2c$ matrix $\begin{pmatrix} A & C \\ B & D \end{pmatrix}$ preserves the skew form represented by the matrix $\begin{pmatrix} 0 & \mathbf{1}_c \\ -\mathbf{1}_c & 0 \end{pmatrix}$, where $\mathbf{1}_c$ is the $c \times c$ identity matrix. The claim now follows. □

*Remark* 4.5. It will be helpful to select a basis for $V_{0,\mathbb{R}}$. The polarization condition $0 < Q(u, N^2u)$ for all $0 \neq u \in V_{\mathbb{R}}^{2,2}$ implies that we may pick a basis $\{e_s \mid 1 \leq s \leq a\}$ of $V_{\mathbb{R}}^{2,2}$ so that $Q(e_s, N^2e_t) = \delta_{st}$. Note that $\{N^2e_s\}_{s=1}^a$ is a basis of $V_{\mathbb{R}}^{0,0}$, and that the $\{Ne_s\}_{s=1}^a \subset V_{\mathbb{R}}^{1,1}$ are linearly independent. Moreover, $Q(Ne_s, Ne_t) = -\delta_{st}$, and we may complete the $\{Ne_s\}$ to a basis $\{Ne_s, f_j \mid s = 1, \ldots, a\,,\ j = 1, \ldots, b\}$ so that $Q(Ne_s, f_j) = 0$ and $Q(f_j, f_k) = \delta_{jk}$ for all $s, j, k$. It then follows that, relative to the basis $\{e_s\,,\, Ne_s, f_j\,,\, N^2e_s\}$ of $V_{0,\mathbb{R}}$, the polarization $Q$ is given by

$$Q|_{V_{0,\mathbb{R}}} \;=\; \begin{pmatrix} 0 & 0 & \mathbf{1}_a \\ 0 & \mathbf{1}_{a,b} & 0 \\ \mathbf{1}_a & 0 & 0 \end{pmatrix}, \quad \text{where} \quad \mathbf{1}_{a,b} \;=\; \begin{pmatrix} -\mathbf{1}_a & 0 \\ 0 & \mathbf{1}_b \end{pmatrix}.$$

**Proposition 4.6.** *Let $(F^\bullet, N)$ be an $\mathbb{R}$–split nilpotent orbit on a period domain parameterizing effective weight two Hodge structures. The subgroup $M^0_\mathbb{R} \subset M_\mathbb{R}$ is the connected identity component of*

$$\begin{aligned} G^{0,0}_\mathbb{R} \;&:=\; \{g \in G_\mathbb{R} \mid g(V^{p,q}) \subset V^{p,q}\ \forall\ p,q\} \\ &=\; \big(G^{0,0}_\mathbb{R} \cap \mathrm{Aut}(V_0)\big) \;\times\; \big(G^{0,0}_\mathbb{R} \cap \mathrm{Aut}(V_1)\big) \;\times\; \big(G^{0,0}_\mathbb{R} \cap \mathrm{Aut}(V_2)\big) \qquad (4.7)\\ &\simeq\; \big(\mathrm{GL}(a,\mathbb{R}) \times \mathrm{O}(b,a)\big) \;\times\; \mathrm{GL}(c,\mathbb{R}) \;\times\; \mathrm{O}(d,\mathbb{R})\,. \end{aligned}$$

*Proof.* The proof is very like that of Proposition 4.2, so we will merely sketch the argument. To see that $M^0_\mathbb{R}$ is the connected identity component of $G^{0,0}_\mathbb{R}$, it suffices to observe that $G^{0,0}_\mathbb{R}$ is the maximal subgroup of $G_\mathbb{R}$ with Lie algebra $\mathfrak{m}^0_\mathbb{R} = \mathfrak{g}^{0,0}_\mathbb{R}$, and to recall that $M^0_\mathbb{R}$ is connected.

The subgroup $G^{0,0}_\mathbb{R} \subset M_\mathbb{R}$ is determined as follows: First, note that the argument establishing the first equality of (4.3) also yields the equality of (4.7). The factor $G_\mathbb{R} \cap \mathrm{Aut}(V_2) = \mathrm{O}(d,\mathbb{R})$ of $M_\mathbb{R}$ preserving $V_2 = V^{2,0}$ necessarily lies in $G^{0,0}_\mathbb{R}$. Next, from observations (iii) and (iv) in the proof of Proposition 4.2, we see that

$$G^{0,0}_\mathbb{R} \cap \mathrm{Aut}(V_1) \;\simeq\; \left\{ \begin{pmatrix} (D^t)^{-1} & 0 \\ 0 & D \end{pmatrix} \quad D \in \mathrm{GL}(c,\mathbb{R}) \right\}. \qquad (4.8)$$

Likewise, working with the basis of Remark 4.5 we see that

$$G^{0,0}_\mathbb{R} \cap \mathrm{Aut}(V_0) \;\simeq\; \left\{ \begin{pmatrix} E_1^{-1} & 0 & 0 \\ 0 & E_2 & 0 \\ 0 & 0 & E_1^t \end{pmatrix} \;\middle|\; \begin{array}{c} E_1 \in \mathrm{GL}(a,\mathbb{R}) \\ E_2 \in \mathrm{O}(b,a) \end{array} \right\}. \qquad (4.9)$$

□

Next we describe $\mathfrak{g}^{-1,-1}_\mathbb{R}$ as a $M^0_\mathbb{R}$–representation.

**Proposition 4.10.** *Let $(F^\bullet, N)$ be an $\mathbb{R}$–split nilpotent orbit on a period domain parameterizing effective weight two Hodge structures. As a $G^{0,0}_\mathbb{R}$–representation*

$$\mathfrak{g}^{-1,-1}_\mathbb{R} \;\simeq\; \mathrm{Sym}^2\,\mathbb{R}^c \;\oplus\; \mathrm{Hom}_\mathbb{R}(\mathbb{R}^a, \mathbb{R}^{a+b})\,. \qquad (4.11)$$

*More precisely:*

(i) *The factor* $\mathrm{O}(d,\mathbb{R}) \simeq G_{\mathbb{R}}^{0,0} \cap \mathrm{Aut}(V_2)$ *acts trivially on* $\mathfrak{g}_{\mathbb{R}}^{-1,-1}$.

(ii) *If* $g \in G_{\mathbb{R}}^{0,0} \cap (\mathrm{Aut}(V_1) \times \mathrm{Aut}(V_0)) \simeq \mathrm{GL}(c,\mathbb{R}) \times (\mathrm{GL}(a,\mathbb{R}) \times \mathrm{O}(b,a))$ *is represented by* $(D; E_1, E_2)$ *as in* (4.8) *and* (4.9)*, then the action of* $g$ *on* $(X,Y) \in \mathrm{Sym}^2\mathbb{R}^c \oplus \mathrm{Hom}_{\mathbb{R}}(\mathbb{R}^a, \mathbb{R}^{a+b})$ *is* $(X,Y) \mapsto (DXD^t, E_2 Y E_1)$.

*Proof.* The key observation is that $\mathfrak{g}_{\mathbb{R}}^{-1,-1}$ decomposes as

$$\mathfrak{g}_{\mathbb{R}}^{-1,-1} \;=\; \left(\mathfrak{g}_{\mathbb{R}}^{-1,-1} \cap \mathrm{Hom}(V_1)\right) \;\oplus\; \left(\mathfrak{g}_{\mathbb{R}}^{-1,-1} \cap \mathrm{Hom}(V_0)\right) .$$

From the characterization of $G_{\mathbb{R}} \cap \mathrm{Aut}(V_1)$ in the proof of Proposition 4.2 (items (iii) and (iv)) we see that the the first summand is

$$\mathfrak{g}_{\mathbb{R}}^{-1,-1} \cap \mathrm{Hom}(V_1) \;\simeq\; \left\{ \begin{pmatrix} 0 & 0 \\ X & 0 \end{pmatrix} \;\middle|\; X = X^t \text{ a symmetric } c \times c \text{ matrix} \right\} . \tag{4.12}$$

Likewise, relative to the basis of Remark 4.5, the second summand is

$$\mathfrak{g}_{\mathbb{R}}^{-1,-1} \cap \mathrm{Hom}(V_0) \;\simeq\; \left\{ \begin{pmatrix} 0 & 0 & 0 \\ Y & 0 & 0 \\ 0 & -Y^t \mathbf{1}_{a,b} & 0 \end{pmatrix} \;\middle|\; Y \text{ a } (a+b) \times a \text{ matrix} \right\} . \tag{4.13}$$

This establishes (4.11).

To complete the proof it remains to check that the adjoint action of $g \in G_{\mathbb{R}}^{0,0}$ on $\mathfrak{g}_{\mathbb{R}}^{-1,-1}$ is as described. Making use of the identifications (4.8) and (4.9), this is an exercise in matrix multiplication that we leave to the reader. □

4.2. **Orbit decomposition.** Assume the hypotheses of Proposition 4.10 and use the decomposition (4.11) to write $N = N_1 + N_0$ with

$$N_1 \in \mathrm{Sym}^2\,\mathbb{R}^c \quad \text{and} \quad N_0 \in \mathrm{Hom}_{\mathbb{R}}(\mathbb{R}^a, \mathbb{R}^{a+b}) .$$

Note that

- $N_0 \neq 0$ if and only if $V_{\mathbb{R}}^{2,2} \neq 0$ (equivalently, $a \neq 0$), and
- $N_1 \neq 0$ if and only if $V^{2,1} \neq 0$ (equivalently, $c \neq 0$).

Propositions 4.6 and 4.10 imply

$$\mathcal{N}^0 \;=\; \mathcal{N}_0^0 \times \mathcal{N}_1^0 \tag{4.14}$$

where

- $\mathcal{N}_0^0$ is the orbit of $N_0$ under $\mathrm{GL}(a,\mathbb{R}) \times \mathrm{O}(b,a)$, and
- $\mathcal{N}_1^0$ is the orbit of $N_1$ under $\mathrm{Sp}(2c,\mathbb{R})$.

Likewise

$$\sigma \;=\; \sigma_0 \times \sigma_1\,, \tag{4.15}$$

with $\sigma_i = \sigma \cap \mathcal{N}_i^0$. From Proposition 4.2 we see that the second and third factors, $\mathrm{Sp}(2c,\mathbb{R})$ and $\mathrm{O}(d,\mathbb{R})$, of $M_\mathbb{R}$, are always Hermitian (§3.4.1). The third factor we disregard as it acts trivially on $\mathfrak{g}_\mathbb{R}^{-1,-1}$, and we have

$$\max\dim_\mathbb{R} \sigma_1 \;=\; \tfrac{1}{2}c(c+1)\,. \tag{4.16}$$

The first factor $\mathfrak{so}(a+b,2a)$ is Hermitian if and only if $a = 0,1$. In this case $\max\dim_\mathbb{R}\sigma_0 = a(a+b)$. The first factor is contact (§3.4.2) if and only if $a = 2$, and in this case $\max\dim_\mathbb{R}\sigma_0 = (2+b)$. To summarize, with $\max\dim_\mathbb{R}\sigma$ denoting the maximal possible dimension of a nilpotent cone $\sigma \subset \mathcal{N}^0$ underlying a nilpotent orbit:

- If $a = 0,1$, then $\max\dim_\mathbb{R}\sigma = a(a+b) + \frac{1}{2}c(c+1)$.
- If $a = 2$, then $\max\dim_\mathbb{R}\sigma = (2+b) + \frac{1}{2}c(c+1)$.

For the cases $a = 0,1,2$, the maximal nilpotent cones $\sigma$ underlying a nilpotent orbit all have the same dimension. This will not be the case when $a > 2$. However, one may use [22] to identify the dimensions of the maximal cones.

It remains to describe the orbits $\mathcal{N}_0^0$ and $\mathcal{N}_1^0$.

**Proposition 4.17.** *Let $(F^\bullet, N)$ be an $\mathbb{R}$–split nilpotent orbit on a period domain parameterizing effective weight two Hodge structures. Let $N = N_0 + N_1$ be the decomposition given by* (4.11).

(a) *The orbit $\mathcal{N}_1^0 = \mathrm{GL}(c,\mathbb{R}) \cdot N_1$ is the set*

$$\mathcal{X} \;:=\; \left\{ X \in \mathfrak{g}_\mathbb{R}^{-1,-1} \cap \mathrm{End}V_1 \mid 0 < \mathbf{i}Q(u, X\bar{u}) \;\forall\; 0 \neq u \in V^{2,1} \right\}\,.$$

*Under the identification of* (4.12)*, this orbit is parameterized by the $c \times c$ symmetric matrices $X = DD^t$ with $D \in \mathrm{GL}(c,\mathbb{R})$.*

(b) *The orbit $\mathcal{N}_0^0 = (\mathrm{GL}(a,\mathbb{R}) \times \mathrm{O}(a,b)^\circ) \cdot N_0$ is the connected component $\mathcal{Y}^\circ$ of*

$$\mathcal{Y} \;:=\; \left\{ Y \in \mathfrak{g}_\mathbb{R}^{-1,-1} \cap \mathrm{End}V_0 \mid 0 < Q(u, Y^2 u) \;\forall\; 0 \neq u \in V_\mathbb{R}^{2,2} \right\}\,.$$

*containing $N_0$. Under the identification of* (4.13)*, $\mathcal{Y}$ is parameterized by the $(a+b)\times a$ matrices $Y = (\alpha E_1, \beta E_1)$ with $\alpha$ an $a\times a$ matrix and $\beta$ a $b\times a$ matrix such that $\alpha^t\alpha - \beta^t\beta = \mathbf{I}_a$, and $E_1 \in \mathrm{GL}(a,\mathbb{R})$.*

The proofs of $\mathcal{N}_1^0 = \mathcal{X}$ and $\mathcal{N}_0^0 = \mathcal{Y}^\circ$ are variations on the argument establishing Lemma 3.5, and will extend to more general situations in a fairly straightforward manner.

In the proof below, it will be helpful to keep in mind that

$$\begin{aligned}
\mathfrak{g}_{\mathbb{R}}^{-1,-1} \cap \mathrm{End}(V_1) &= \mathrm{End}(V_{\mathbb{R}}, Q) \cap \mathrm{Hom}(V^{2,1}, V^{1,0})\,, \\
\mathfrak{g}_{\mathbb{R}}^{-1,-1} \cap \mathrm{End}(V_0) &= \mathrm{End}(V_{\mathbb{R}}, Q) \cap \left(\mathrm{Hom}(V^{2,2}, V^{1,1}) \oplus \mathrm{Hom}(V^{1,1}, V^{0,0})\right)\,.
\end{aligned}$$

*Proof.* By Lemma 3.5, we know that $\mathcal{N}^0$ is open in $\mathfrak{g}_{\mathbb{R}}^{-1,-1}$. From Proposition 4.10 and (4.14), we see that this is equivalent to the two conditions that $\mathcal{N}_i^0$ is open in $\mathfrak{g}_{\mathbb{R}}^{-1,-1} \cap \mathrm{Hom}(V_i)$, $i = 0, 1$.

To establish $\mathcal{N}_1^0 = \mathcal{X}$, observe that $\mathcal{X}$ is open, convex and preserved under the action of $M_{\mathbb{R}}^0$. By definition $N_1 \in \mathcal{X}$. Since $\mathcal{X}$ is preserved under the action of $M_{\mathbb{R}}^0$, it is immediate that $\mathcal{N}_1^0 \subset \mathcal{X}$. It remains to show that equality holds. An argument analogous to that establishing Lemma 3.5 shows that $\mathcal{X}$ is a union of open $M_{\mathbb{R}}^0$–orbits. The equality $\mathcal{X} = \mathcal{N}_1^0$ then follows from the convexity of $\mathcal{X}$.

Under the identification (4.12), $N_1$ is represented by $X = \mathbf{1}_c$. (View this as indicating that $N_1$ gives us a specific isomorphism $V^{2,1} \simeq V^{1,0}$.) So the action of $g = (D; E_1, E_2) \in M_{\mathbb{R}}^0$ on $N_1$ is $\mathbf{1}_c \mapsto DD^t$, by Proposition 4.10.[4]

We briefly sketch the argument establishing $\mathcal{N}_0^0 = \mathcal{Y}^\circ$ which is very like that above for $\mathcal{N}_1^0 = \mathcal{X}$. Again we observe that $\mathcal{Y}$ is open and preserved under the action of $M_{\mathbb{R}}^0$ (but *not* convex[5]). As above $\mathcal{N}_0^0 \subset \mathcal{Y}$, and one may show that $\mathcal{Y}$ is a union of open $M_{\mathbb{R}}^0$–orbits.

Under the identification (4.13), $N_0$ is represented by $Y = (\mathbf{1}_a\ \mathbf{0}_{b,a})^t$, where $\mathbf{0}_{b,a}$ is the $b\times a$ zero matrix. Decompose $E_2 \in \mathrm{O}(b,a)$ as

$$E_2 = \begin{pmatrix} \alpha & * \\ \beta & * \end{pmatrix},$$

[4]The Cholesky decomposition yields a factorization of every element $X' \in \mathcal{X}$ of the form $X' = DD^t$.

[5]This is a feature of the non-classical case that $D$ is not Hermitian symmetric.

with $\alpha$ an $a \times a$ matrix and $\beta$ a $b \times a$ matrix. Then Proposition 4.10 asserts that the action of $g = (D; E_1, E_2) \in M^0_{\mathbb{R}}$ on $N_0$ is

$$\begin{pmatrix} \mathbf{1}_a \\ \mathbf{0}_{b,a} \end{pmatrix} \mapsto \begin{pmatrix} \alpha E_1 \\ \beta E_1 \end{pmatrix} .$$

Since $E_2 \in \mathrm{O}(b, a)$, we have $\alpha^t\alpha - \beta^t\beta = \mathbf{I}_a$.

This completes the proof of the proposition. As a final remark, and keeping the identification (4.13) in mind, we note that $(g \cdot N_0)^2 : V^{2,2}_{\mathbb{R}} \to V^{0,0}_{\mathbb{R}}$ is represented by

$$\begin{pmatrix} 0 & 0 & 0 \\ 0 & 0 & 0 \\ E_1^t(\alpha^t\alpha - \beta^t\beta)E_1 & 0 & 0 \end{pmatrix} = \begin{pmatrix} 0 & 0 & 0 \\ 0 & 0 & 0 \\ E_1^tE_1 & 0 & 0 \end{pmatrix} .$$

□

## 5. Example: period domain for $\mathbf{h} = (2, *, 2)$

The goal of this section is illustrate how the material of §§2–4 can be applied to study nilpotent orbits on the period domain $D$ for Hodge numbers

$$\mathbf{h} \;=\; (2, h^{1,1}, 2) .$$

In §5.1 we identify the $\mathbb{R}$–split polarized mixed Hodge structures on $D$ ("Step 1" of Remark 3.10). Section 5.2 describes representation theoretic constraints on the degenerations coming from the faces of a cone underlying a nilpotent orbit. In §5.3 we see that these are the only constraints: one may construct nilpotent cones, that underlie nilpotent orbits, from commuting $\mathfrak{sl}(2)$s exhibiting all remaining degenerations.

### 5.1. **The PMHS.** Set

$$m \;=\; 4 + h^{1,1} \quad \text{so that} \quad V_{\mathbb{R}} \;=\; \mathbb{R}^m .$$

We have

$$G_{\mathbb{R}} \;=\; \mathrm{Aut}(V_{\mathbb{R}}, Q) \;\simeq\; \mathrm{O}(4, m-4) \quad \text{and} \quad \mathfrak{g}_{\mathbb{R}} \;=\; \mathrm{End}(V_{\mathbb{R}}, Q) \;\simeq\; \mathfrak{so}(4, m-4) .$$

Modulo the action of $G_{\mathbb{R}}$ there are at most five polarized $\mathbb{R}$–split PMHS on $D$. The Hodge diamonds (HD) for these PMHS are depicted in Table 5.1, where they are denoted I, . . . , V. As you can see from the table, we need $m \geq 8$ to get all five degenerations.

TABLE 5.1. Hodge diamonds and signed Young diagrams

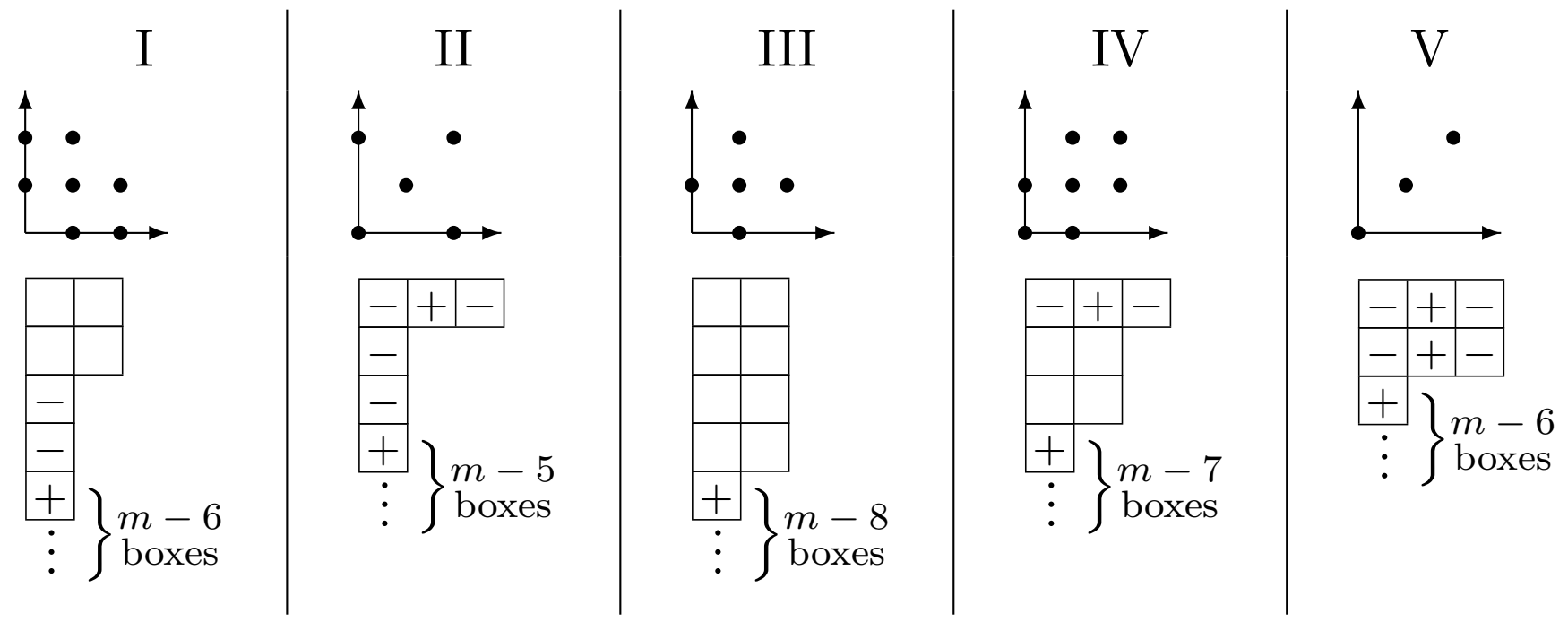

TABLE 5.2. The diagonal group $M_{\mathbb{R}}$

| | $M_{\mathbb{R}}^{\mathrm{ss}}$ | $\mathrm{rank}_{\mathbb{R}}$ | $\max \dim_{\mathbb{R}} \sigma$ |
|---|---|---|---|
| I | $\mathrm{O}(m-6) \times \mathrm{Sp}(2,\mathbb{R})$ | 1 | 1 |
| II | $\mathrm{O}(m-4,2)$ | 2 | $m-4$ |
| III | $\mathrm{O}(m-8) \times \mathrm{Sp}(4,\mathbb{R})$ | 2 | 3 |
| IV | $\mathrm{O}(m-6,2) \times \mathrm{Sp}(2,\mathbb{R})$ | 3 | $m-5$ |
| V | $\mathrm{O}(m-4,4)$ | 4 | $m-4$ |

By Lemma 3.11, the number of commuting $\mathfrak{sl}(2)$'s obtained from the Cattani, Kaplan and Schmid construction [4] is bounded by the real rank of $M_{\mathbb{R}}^{\mathrm{ss}}$. The subgroups, which are determined by Proposition 4.2, are listed in Table 5.2. In the table the maximal dimension of the nilpotent cones $\sigma \subset \mathcal{N}^0$ underlying a nilpotent orbit are taken from §4.2.

5.2. **Degenerations coming from the faces of the cone.** Together Theorem 2.21 and Corollary 3.6 provide representation theoretic constraints on the degenerations associated with the faces of a nilpotent cone underlying a nilpotent orbit on a period domain $D$. While illustrated with examples in weight two period domains, the discussion of this section is general and applies to arbitrary period domains.

Given a polarizing nilpotent $N \in \mathfrak{g}_{\mathbb{R}}$, let

$$\mathcal{N} \;=\; \mathrm{Ad}(G_{\mathbb{R}}) \cdot N$$

denote the conjugacy class. As discussed in §2.5, these orbits are enumerated by (partially) signed Young diagrams, a.k.a. Đoković's chromosomes, and the diagrams are determined by the Hodge diamond of $V_{\mathbb{C}}$. Recall that if $\sigma = \{\lambda^i N_i \mid \lambda^i > 0\}$ is a nilpotent cone underlying a nilpotent orbit, then Corollary 3.6 yields

$$N \in \sigma \implies \sigma \subset \mathcal{N} \quad \text{and} \quad N_i \subset \overline{\mathcal{N}}\,. \tag{5.1}$$

In particular, every nilpotent in $\sigma$ is of the same type as $N$. (E.g. if $N$ is of type II, then every element of the cone is of type II.) So we can speak of the "type of $\sigma$," and there are fives types of cones on the period domains for $\mathbf{h} = (2, h^{1,1}, 2)$.

Given two nilpotents $N_1$ and $N_2$ we write $\mathcal{N}_1 \le \mathcal{N}_2$ if $\mathcal{N}_1 \subset \overline{\mathcal{N}_2}$. Đoković's Theorem 2.21 characterizes this partial ordering; for the conjugacy classes $\mathcal{N}$ of the polarizing $N$ in Table 5.1 we have

$$\mathcal{N}_{\mathrm{I}} \;<\; \left\{ \begin{array}{c} \mathcal{N}_{\mathrm{II}} \\ \mathcal{N}_{\mathrm{III}} \end{array} \right\} \;<\; \mathcal{N}_{\mathrm{IV}} \;<\; \mathcal{N}_{\mathrm{V}}\,. \tag{5.2}$$

Given (5.1), and the dimension constraints listed in Table 5.2, this tells us something about the degenerations corresponding to the faces of the cone. For example,

(a) Any face of type I is necessarily one–dimensional.
(b) If $\sigma$ is of type II, then the faces of $\sigma$ are either of type I or type II — because only $\mathcal{N}_{\mathrm{I}}, \mathcal{N}_{\mathrm{II}} \le \mathcal{N}_{\mathrm{II}}$.
(c) Likewise, if $\sigma$ is of type III, the faces of $\sigma$ are either of type I or type III.
(d) Any face of type III is of dimension at most three.
(e) If $\sigma$ is of type IV, then no face is of type V.

5.3. **Commuting horizontal** $\mathrm{SL}(2)$**s.** The representation theoretic constraints (b), (c) and (e) above are the only restrictions on the degenerations coming from faces of the cone.[6] In fact, we may construct cones from commuting $\mathfrak{sl}(2)$s whose faces realize all *a priori* possible degenerations. Moreover, we may choose these $\mathfrak{sl}(2)$s to be of the simplest type – corresponding to a root $\alpha$ – and to be defined over $\mathbb{Q}$. (The resulting degenerations correspond to those of [17, §6] obtained from Cayley transforms $\mathbf{c}_\alpha$ in non-compact imaginary roots.) The details are as follows.

It will be convenient in the discussion that follows to use the representation theoretic language of roots and root spaces $\mathfrak{g}^\alpha \subset \mathfrak{g}_{\mathbb{C}}$. For the reader who does not regularly work with this nomenclature, we recall that $\mathfrak{g}^\alpha$ is 1–dimensional. Bases of the four root spaces $\mathfrak{g}^{\alpha_i}$ that we will consider below can be described as follows: Given a basis $\{e_1,\ldots,e_m\}$ of $V_\mathbb{C}$, let $\{e^1,\ldots,e^m\}$ denote the dual basis of $V_\mathbb{C}^*$, and $\{e^k_j := e_j\otimes e^k \mid 1\le j,k\le m\}$ the associated basis of $\mathrm{End}(V_\mathbb{C})$. It is possible to chose $\{e_1,\ldots,e_m\}$ so that

$$Q(e_j,e_k) \ = \ \delta^{m+1}_{j+k}$$

and

$$\begin{array}{lll}
\mathfrak{g}^{\alpha_1} \text{ and } \mathfrak{g}^{-\alpha_1} & \text{are spanned by} & e^3_2 - e^{m-1}_{m-2} \text{ and } e^2_3 - e^{m-2}_{m-1}\,;\\
\mathfrak{g}^{\alpha_2} \text{ and } \mathfrak{g}^{-\alpha_2} & \text{are spanned by} & e^4_1 - e^{m}_{m-3} \text{ and } e^1_4 - e^{m-3}_{m}\,;\\
\mathfrak{g}^{\alpha_3} \text{ and } \mathfrak{g}^{-\alpha_3} & \text{are spanned by} & e^{m-2}_2 - e^{m-1}_{3} \text{ and } e^2_{m-2} - e^{3}_{m-1}\,;\\
\mathfrak{g}^{\alpha_4} \text{ and } \mathfrak{g}^{-\alpha_4} & \text{are spanned by} & e^{m-3}_1 - e^{m}_{4} \text{ and } e^1_{m-3} - e^{4}_{m}\,,
\end{array}$$

respectively.[7]

[6]This fails for more general period domains: there are additional constraints, essentially imposed by horizontality [18].

Each root $\alpha$ determines a root subalgebra

$$\mathfrak{sl}^\alpha(2,\mathbb{C}) \ := \ \mathfrak{g}^\alpha \,\oplus\, [\mathfrak{g}^\alpha,\mathfrak{g}^{-\alpha}] \,\oplus\, \mathfrak{g}^{-\alpha} \ \simeq \mathfrak{sl}(2,\mathbb{C})\,.$$

Note that the four TDS $\mathfrak{sl}^{\alpha_i}(2,\mathbb{C})$ commute. Given $\{\beta_1,\ldots,\beta_s\}\subset\{\alpha_1,\ldots,\alpha_4\}$, let

$$\mathfrak{sl}(\beta_1,\ldots,\beta_s;\mathbb{C}) \ := \ \mathfrak{sl}^{\beta_1}(2,\mathbb{C})\times\cdots\times\mathfrak{sl}^{\beta_s}(2,\mathbb{C})$$

denote the subalgebra of commuting $\mathfrak{sl}(2,\mathbb{C})$'s. A *generic nilnegative of* $\mathfrak{sl}(\beta_1,\ldots,\beta_s;\mathbb{C})$ is an element $N\in\mathfrak{g}^{-\beta_1}\oplus\cdots\oplus\mathfrak{g}^{-\beta_s}$ whose projection to each negative root space $\mathfrak{g}^{-\beta_i}$ is nonzero.

Following the approach of [17, §6], we may construct commuting $\mathfrak{sl}(2)$'s, each of which is defined over $\mathbb{Q}$, so that

$$\left.\begin{array}{c}
\mathfrak{sl}^{\alpha_i}(2,\mathbb{Q})\\
\mathfrak{sl}(\alpha_i,\alpha_j;\mathbb{Q})\,,\ (i,j)=(1,3),(2,4)\\
\mathfrak{sl}(\alpha_i,\alpha_j;\mathbb{Q})\,,\ (i,j)\neq(1,3),(2,4)\\
\mathfrak{sl}(\alpha_i,\alpha_j,\alpha_k;\mathbb{Q})\,,\ i,j,k \text{ distinct}\\
\mathfrak{sl}(\alpha_1,\ldots,\alpha_4;\mathbb{Q})
\end{array}\right\}
\begin{array}{c}\text{admits a}\\ \text{generic}\\ \text{nilnegative}\\ \text{of type}\end{array}
\left\{\begin{array}{l}
\mathrm{I}\\ \mathrm{II}\\ \mathrm{III}\\ \mathrm{IV}\\ \mathrm{V}\,.
\end{array}\right.$$

[7] In terms of the simple roots $\{\sigma_1,\ldots,\sigma_r\}$ of $\mathfrak{g}_\mathbb{C}$, we have $\alpha_1=\sigma_2$, $\alpha_2=\sigma_1+\sigma_2+\sigma_3$, $\alpha_3=\sigma_2+2\sigma_3+\alpha'$ and $\alpha_4=\sigma_1+\sigma_2+\sigma_3+\alpha'$, where $\alpha'=\begin{cases}2(\sigma_4+\cdots+\sigma_r) & \text{if } m=2r+1\,,\\ 2(\sigma_4+\cdots+\sigma_{r-2})+\sigma_{r-1}+\sigma_r & \text{if } m=2r\,.\end{cases}$

These commuting $\mathfrak{sl}(2)$s yield 4–dimensional nilpotent cones $\sigma$ whose faces realize every combination of degeneration not ruled out by §5.2.

## 6. Deligne systems

In this section, we give a counterexample to the following assertion of Kato [16, Theorem 1.4], and provide a corrected statement:

**Theorem 6.1.** *Let $(V, W, N_1, \ldots, N_n, F)$ be a Deligne–Hodge system of $n$ variables. Then for $N'_j = \sum_{k=1}^{j} a_{j,k} N_k$ $(1 \leq j \leq n)$ with $a_{j,k} > 0$ $(1 \leq k \leq j \leq n)$ such that $a_{j,k}/a_{j,k+1} \gg 0$ $(1 \leq k < j \leq n)$, $(V, W, N'_1, \ldots, N'_n, F)$ is an* IMHM *of $n$ variables.*

In a nutshell, the problem is that one needs a polarizability condition on the original Deligne-Hodge system to guarantee that some modification is an IMHM. More precisely, the flaw in the proof, which starts on page 857 of [16] appears to be the following: The second sentence of the third paragraph of the proof says, "on $\mathrm{gr}^W_w$, put *the* bilinear form in Proposition 3.2.7." We've added the emphasis here because the key problem with the proof appears to be the word "the." To every object $\mathbf{V} = (V, W, N_1, \ldots, N_r, F)$ in $\mathrm{DH}_r$, one can associate an $\mathrm{SL}_2$-orbit or, equivalently, an object $\hat{\mathbf{V}}$ in $\widehat{\mathrm{DH}}_r$. Any $\hat{\mathbf{V}}$ can be polarized by a bilinear form $Q$. But $Q$ is not unique. If $\hat{\mathbf{V}}$ is irreducible, then $Q$ is unique up to non-zero scalar multiple. But, in general, there is no unique $Q$ even up to scalar multiple. This becomes a problem because the condition that the $N_i$ be infinitesimal isometries of $Q$ is non-trivial, and they impose (possibly contradictory) conditions on what $Q$ can be.

Before proceeding to give the counterexample in §6.2, we first provide a short account of Deligne systems in §6.1. In §6.3 we revisit the counterexample from a categorical point of view. In §6.4 we show that Kato's theorem holds in the presence of a suitable graded-polarization condition which can be stated in terms of the associated $\mathrm{SL}_2$-orbit. In §6.5df, we discuss the geometry of Deligne systems with a given underlying $\mathrm{SL}_2$-orbit.

The definition of infinitesimal mixed Hodge module appears in §4 of [15], and the notion of Deligne system will be defined in the next section. For completeness, we record the definition of an IMHM here:

*Definition* 6.2. An infinitesimal mixed Hodge module consists of

(1) A finite dimensional real vector space $V_{\mathbb{R}}$ equipped with an increasing filtration $W$ and a collection of non-degenerate bilinear forms $Q_k : Gr^W_k \otimes Gr^W_k \to \mathbb{R}$ of parity $(-1)^k$;
(2) A decreasing filtration $F$ of $V_{\mathbb{C}} = V_{\mathbb{R}} \otimes \mathbb{C}$;
(3) Nilpotent endomorphisms $N_1, \dots, N_r$ of $V_{\mathbb{R}}$ which preserve $W$ and act by infinitesimal isometries on $Gr^W$.

such that

(a) $N_j(F^p) \subset F^{p-1}$ for all $j$ and $p$;
(b) $e^{\sum_j z_j N_j} F$ induces a nilpotent orbit of pure Hodge structure of weight $k$ on $Gr^W_k$ which is polarized by $Q_k$;
(c) For any subset $J$ of $\{1, \dots, r\}$ there exists a relative weight filtration $M(J)$ such that (i) $N_j M_k(J) \subset M_{k-2}(J)$ for all $j \in J$ and (ii) $M(J)$ is the weight filtration of $\sum_{j \in J} N_j$ relative to $W$.

In particular, if $W$ is pure of weight $k$–i.e. $Gr^W_\ell = 0$ unless $\ell = k$ then a IMHM is the same thing as a nilpotent orbit of pure Hodge structure of weight $k$. [Condition $(c)$ follows from the results of Cattani and Kaplan].

6.1. **Preliminary Remarks.** Fix a field $K$ of characteristic zero, and let $W$ be an increasing filtration of a finite dimensional $K$-vector space $V$. Then, a grading of $W$ is a semisimple endomorphism $Y$ of $V$ with integral eigenvalues such that

$$W_k = \bigoplus_{j \le k} E_j(Y) \tag{6.3}$$

where $E_j(Y)$ is the $j$-eigenspace of $Y$.

Let $N$ be a nilpotent endomorphism of $V$ which preserve $W$, i.e. $N(W_k) \subseteq W_k$. Then (cf. [26]), there exists at most one relative weight filtration $M = M(N, W)$ such that

(a) $N(M_k) \subseteq M_{k-2}$ for all $k$;
(b) If $Gr^W_k$ is non-zero and $\ell \ge 0$ then the induced map

$$N^\ell : Gr^M_{k+\ell} Gr^W_k \to Gr^M_{k-\ell} Gr^W_k.$$

is an isomorphism for each non-negative integer $\ell$.

In the geometric case of interest, $W$ is the weight filtration of an admissible variation of mixed Hodge structure over the punctured disk and $N$ is the local monodromy logarithm. In this setting $N$ is a $(-1,-1)$-morphism of the limit mixed Hodge structure $(F, M)$ where $M = M(N, W)$. Let

$$V = \bigoplus_{p,q} I^{p,q} \tag{6.4}$$

be the associated Deligne bigrading of $(F, M)$ and $Y = Y_{(F,W)}$ be the grading of $M$ which acts as multiplication by $p+q$ on $I^{p,q}$. Then,

$$[Y, N] = -2N \tag{6.5}$$

since $N$ is a $(-1,-1)$-morphism. Moreover, $Y$ preserves $W$ since $(F, M)$ induces a mixed Hodge structure on each $W_k$ by (3.13) of [26].

*Definition* 6.6. A 1-variable Deligne system over $K$ consists of the following data:

- An increasing filtration $W$ of a finite dimensional $K$-vector space $V$;
- A nilpotent endomorphism $N$ of $V$ which preserves $W$ such that $M = M(N, W)$ exists;
- A grading $Y$ of $M$ which preserves $W$ and satisfies $[Y, N] = -2N$.

A morphism of Deligne systems $(W, N, Y) \to (\tilde{W}, \tilde{N}, \tilde{Y})$ is an endomorphism $T$ the underlying $K$-vector spaces such that $T(W_i) \subset \tilde{W}_i$ and

$$\tilde{Y} \circ T - T \circ Y = 0, \qquad \tilde{N} \circ T - T \circ N = 0\,.$$

*Example* 6.7. By the remarks of the previous paragraphs, if $(e^{zN}F, W)$ is an admissible nilpotent orbit then $(W, N, Y)$ is a Deligne system where $Y = Y_{(F,M)}$ and $M = M(N, W)$.

To continue, we recall the following: An $\mathrm{sl}_2$-pair consists of a nilpotent endomorphism $N$ of a finite dimensional $K$-vector space $V$ and grading $H$ of the monodromy weight filtration $W(N)$ such that $[H, N] = -2N$. Moreover, there is a 1-1 correspondence between $\mathrm{sl}_2$-pairs and representations $\rho$ of $\mathrm{sl}_2(K)$ on $V$ such that

$$N = \rho\begin{pmatrix} 0 & 0 \\ 1 & 0 \end{pmatrix}, \qquad H = \rho\begin{pmatrix} 1 & 0 \\ 0 & -1 \end{pmatrix}$$

That such a representation determines an $\mathrm{sl}_2$-pair follows from the structure of the irreducible representations of $\mathrm{sl}_2(K)$. Conversely, given an $\mathrm{sl}_2$-pair, the elements of the kernel of $N : E_{-k}(H) \to E_{-k-2}(H)$ are lowest weight vectors for $\rho$.

In particular, given a 1-variable Deligne system $(W, N, Y)$ let $Y'$ be a grading of $W$ which commutes with $Y$ and

$$N = N_0 + N_{-1} + N_{-2} + \cdots \tag{6.8}$$

be the decomposition of $N$ into eigencomponents relative to $\mathrm{ad}(Y)$. Then, $N_0$ and $H = Y - Y'$ induce the action of an $\mathrm{sl}_2$-pair on each $Gr^W_k$. Let $\rho_k$ be the corresponding representation of $\mathrm{sl}_2(K)$ on $Gr^W_k$ and $\rho$ be the representation of $\mathrm{sl}_2(K)$ which acts as $\rho_k$ on $E_k(Y')$ via isomorphism $E_k(Y') \cong Gr^W_k$. Let

$$N_0^+ = \rho \begin{pmatrix} 0 & 1 \\ 0 & 0 \end{pmatrix} \tag{6.9}$$

Accordingly, given a Deligne system $(W, N, Y)$, each choice of grading $Y'$ of $W$ which commutes with $Y$ determines a corresponding $\mathrm{sl}_2$-triple. Moreover, a short calculation shows that since both $N_0$ and $H = Y - Y'$ commute with $Y'$ so does $N_0^+$. Likewise, it is easy to see that each component $N_{-k}$ appearing in (6.8) is weight $-2$ for $\mathrm{ad}\, Y$ since $[Y, Y'] = 0$.

**Theorem 6.10** ([8]). *Let $(W, N, Y)$ be a Deligne system. Then, there exists an unique, functorial grading $Y' = Y'(N, Y)$ of $W$ which commutes with $Y$ such that*

$$[N - N_0, N_0^+] = 0\,. \tag{6.11}$$

A sketch of Deligne's proof is as follows (cf. Theorem (4.4) in [20]): Let

$$\mathrm{gl}_{-r}(W) = \{\, \alpha \in \mathrm{gl}(V) \mid \alpha(W_k) \subseteq W_{k-r} \,\}$$

and $\mathrm{gl}^Y_{-r}(W)$ be the subalgebra of elements of $\mathrm{gl}_{-r}(W)$ which commute with $Y$. Then, the set of all gradings of $W$ which commute with $Y$ is an affine space upon which the group $\exp(\mathrm{gl}^Y_{-1}(W))$ acts simply transitively via the adjoint action.

Deligne now claims by induction that it is possible construct a sequence of gradings $Y'_0, Y'_1, \dots,$ such that if $(N_0, Y - Y'_r, N_0^+)$ is $\mathrm{sl}_2$-triple associated to $Y'_r$ then

$$[N - N_0, N_0^+] \in \mathrm{gl}^Y_{-r-1}(W)\,. \tag{6.12}$$

Given the finite length of $W$, this process terminates in the desired grading $Y'$. The induction base $r = 0$ is trivial since any grading $Y'_0$ of $W$ which commutes with $Y$ will suffice. Suppose therefore that the required gradings $Y'_1, \dots, Y'_{k-1}$ have been constructed.

Let $N_0, N_{-1}, \dots$ be the components of $N$ relative to $\operatorname{ad} Y'_{k-1}$. Then, since $Y'_{k-1}$ commutes with $(N_0^+, Y - Y'_{k-1}, N_0)$ it follows that

$$N_{-k} = [N_0, \gamma_{-k}] + N'_{-k} \tag{6.13}$$

where $N'_{-k}$ is highest weight $k-2$ for $(N_0^+, Y - Y'_{k-1}, N_0)$ and

$$\gamma_{-k} \in E_k(\operatorname{ad}(Y - Y'_{k-1})) \cap E_{-k}(\operatorname{ad} Y'_{k-1})\,. \tag{6.14}$$

In particular, equation (6.14) implies that $\gamma_{-k} \in E_0(\operatorname{ad} Y)$. A short calculation shows that

$$Y'_k = \operatorname{Ad}(1 + \gamma_{-k}) Y'_{k-1}$$

is a grading of $W$ which commutes with $Y$ and satisfies (6.12) for $r = k$.

Uniqueness and the compatibility of Deligne's result with direct sums, tensor products and duals now follow by standard arguments. In particular, one consequence of Deligne's construction is that $N_{-1} = 0$ since for $k > 0$, equation (6.11) implies that $N_{-k}$ is either zero or of highest weight $k-2$. Likewise, if the data $(V, W, N, Y)$ can be defined over a subfield $k \subset K$ then $Y'$ is defined over $k$.

**Lemma 6.15** ([8])**.** *If $T$ is a morphism of Deligne systems $(W, N, Y) \to (\tilde{W}, \tilde{N}, \tilde{Y})$ then*

$$T \circ Y' = \tilde{Y}' \circ T \tag{6.16}$$

*Proof.* Let $V$ and $\tilde{V}$ be the underlying vector spaces of $(W, N, Y)$ and $(\tilde{W}, \tilde{N}, \tilde{Y})$ respectively. Given endomorphisms $A : V \to V$ and $B : \tilde{V} \to \tilde{V}$, define

$$\operatorname{ad}(B, A)L = B \circ L - L \circ A$$

for any linear map $L : V \to \tilde{V}$. Let $Y' = Y(N, Y)$ and $\tilde{Y}' = Y'(\tilde{N}, \tilde{Y})$. Then,

$$T = \sum_{i \le 0} T_i$$

where $\mathrm{ad}(\tilde{Y}', Y')T_i = iT_i$. Equation (6.16) is the assertion that $T = T_0$. Thus, we must show that $T_i = 0$ for $i < 0$.

The action $\mathrm{ad}(B, A)$ defined above is compatible with Lie brackets and hence the $\mathrm{sl}_2$ representations $(N_0, H, N_0^+)$ and $(\tilde{N}_0, \tilde{H}, \tilde{N}_0^+)$ induce a representation of $\mathrm{sl}_2$ on $\mathrm{Hom}(V, \tilde{V})$. To show that $T_i = 0$ for $i < 0$, we show that

(a) $\mathrm{ad}(\tilde{H}, H)T_i = -iT_i$
(b) $\mathrm{ad}(\tilde{N}_0^+, N_0^+) \circ \mathrm{ad}(\tilde{N}_0, N_0)T_i = 0$

By the representation theory $\mathrm{sl}_2$, condition $(b)$ implies that $\mathrm{ad}(\tilde{N}_0, N_0)T_i = 0$. Condition then implies $T_i = 0$ since $T_i$ has weight $-i > 0$.

Direct calculation shows that $\mathrm{ad}(\tilde{N}_0, N_0)T_0 = 0$ since this is the degree zero eigencomponent of $\mathrm{ad}(\tilde{N}, N)T$ with respect to $\mathrm{ad}(\tilde{Y}', Y')$ and $\mathrm{ad}(\tilde{N}, N)T = 0$ by hypothesis. A similar computation shows that

$$\mathrm{ad}(\tilde{Y}', Y') \circ \mathrm{ad}(\tilde{Y}, Y)T_i = i\mathrm{ad}(\tilde{Y}, Y)T_i$$

Therefore, since $\mathrm{ad}(\tilde{Y}, Y)T = 0$ by hypothesis and $T = \sum_i T_i$ it follows that

$$\mathrm{ad}(\tilde{Y}, Y)T_i = 0$$

for all $i$. Accordingly, $T_i$ is weight $-i$ for $\mathrm{ad}(\tilde{H}, H)$ where $\tilde{H} = \tilde{Y} - \tilde{Y}'$ and $H = Y - Y'$.

By the previous paragraph $\mathrm{ad}(\tilde{N}_0, N_0)T_0 = 0$ and $\mathrm{ad}(\tilde{H}, H)T_0 = 0$, and hence

$$\mathrm{ad}(\tilde{N}_0^+, N_0)T_0 = 0$$

To finish the proof, let $i$ be the largest integer such that $T_i = 0$ and $i < 0$. By the Jacobi identity,

$$\mathrm{ad}(\tilde{N}_0^+, N_0) \circ \mathrm{ad}(\tilde{N}_i, N_i)T_0 = 0$$

On the other hand,

$$\mathrm{ad}(\tilde{N}_0, N_0)T_i + \mathrm{ad}(\tilde{N}_i, N_i)T_0 = 0$$

since this equals the eigencomponent of $\mathrm{ad}(\tilde{N}, N)T$ of weight $i$ with respect to $\mathrm{ad}(\tilde{Y}', Y')$. Comparing these two equations, it follows that

$$\mathrm{ad}(\tilde{N}_0^+, N_0^+) \circ \mathrm{ad}(\tilde{N}_0, N_0)T_i = 0$$

□

*Remark* 6.17. The definition of Deligne system given above is due to Kato. The original (and equivalent) formulation in Christine Schwarz's paper [24] is that a Deligne system is given by the data $(W, N, Y', Y)$. Deligne's theorem (6.10) is then stated as the assertion that given $(W, N, Y)$ as above there is a unique choice of grading $Y'$ of $W$ which completes $(W, N, Y)$ to a Deligne system. As outlined in (6.19) below, Deligne also considered the several variable case. Axioms for the several variable Deligne systems are stated in [24].

The origin of Deligne's letter to Cattani and Kaplan is a question related to a mysterious splitting operation which arises in Schmid's $\mathrm{SL}_2$-orbit theorem. Namely (Prop (2.20),[4]), given any mixed Hodge structure $(F, W)$ on $V$ there exists a unique, real element

$$\delta \in \bigoplus_{p,q<0} \mathrm{gl}(V)^{p,q}$$

such that $(e^{-i\delta}F, W)$ is a mixed Hodge structure which is split over $\mathbb{R}$. On the other hand (Lemma (6.60),[4]), by the $\mathrm{SL}_2$-orbit theorem, if $e^{zN}F$ is a nilpotent orbit of pure Hodge structure then

$$e^{iyN}F = e^{\zeta}\left(1 + \sum_{k>0} g_k y^{-k}\right) e^{iyN} e^{-i\delta} F$$

where $\zeta$ is given by universal Lie polynomials in the Hodge components of $\delta$. In particular, since $\zeta$ is real and depends only on the Hodge components of $\delta$, it follows that

$$\hat{F} = e^{\zeta} e^{-i\delta} F$$

is another split mixed Hodge structure attached to an arbitrary mixed Hodge structure $(F, W)$. In [1] this operation $(F, W) \mapsto (\hat{F}, W)$ is called the $\mathrm{sl}_2$-splitting of $(F, W)$. In the work of Kato and Usui, this operation is called the canonical splitting. For future reference we let $\hat{Y}_{(F,W)} = Y_{(\hat{F},W)}$.

In [8], Deligne asserts that if $(e^{zN}F, W)$ is an admissible nilpotent orbit with limit mixed Hodge structure $(F, M)$ which is split over $\mathbb{R}$ then

$$\hat{Y}_{(e^{iN}F,W)} = Y'(N, Y_{(F,M)}) \tag{6.18}$$

This was proven by the first two authors in [1].

In the second part of his letter, Deligne focuses on applying his construction to a several variable system

$$\begin{pmatrix} W^0 & W^1 & \cdots & W^{r-1} & Y^r \\ & N_1 & \cdots & N_{r-1} & N_r \end{pmatrix} \tag{6.19}$$

where

— $W^0, \ldots, W^r$ are increasing filtrations a finite dimensional $K$-vector space $V$;
— $Y^r$ is a grading of $W^r$;

such that

(a) $N_1, \ldots, N_r$ are commuting nilpotent endomorphisms of $V$ which preserve $W^0$;
(b) $M(N_j, W^{j-1})$ exists and equals $W^j$ for $j = 1, \ldots, r$;
(c) Let $1 \le j \le r$, $0 \le k < j-1$, $\ell \in \mathbb{Z}$, and let $U = W^k_\ell$. Then the restriction $W^j|_U$ of $W^j$ to $U$ is the relative monodromy filtration of $N_j|_U$ with respect to $W^{j-1}|_U$;
(d) $N_j(W^k_\ell) \subseteq W^k_\ell$ for any $j, k, \ell$, and $N_j(W^k_\ell) \subseteq W^k_{\ell-2}$ if $k \ge j$;
(e) $Y^r$ preserves each $W^j$ and $[Y^r, N_j] = -2N_j$ for all $j$.

*Definition* 6.20. A morphism of Deligne systems

$$T : (W, N_1, \ldots, N_r; Y) \to (\tilde{W}, \tilde{N}_1, \ldots, \tilde{N}_r, \tilde{Y}) \tag{6.21}$$

is a homomorphism of the underlying vector spaces such that

$$T(W^j_i) \subset \tilde{W}^j_i, \qquad T \circ N_j = \tilde{N}_j \circ T$$

for all $j$ (and $i$), and $T \circ Y^r = \tilde{Y}^r \circ T$.

**Theorem 6.22.** *(Deligne [8]) Starting from $Y^r$, the iterative application of the construction $Y^{j-1} = Y'(N_j, Y^j)$ to a system (6.19) satisfying the conditions (a)–(e) yields a system of commuting gradings such that if $\hat{N}_j$ is the degree zero part of $N_j$ with respect to $\operatorname{ad} Y^{j-1}$ and $H_j = Y^j - Y^{j-1}$ then*

$$(\hat{N}_1, H_1), \ldots, (\hat{N}_r, H_r) \tag{6.23}$$

*are commuting $\mathrm{sl}_2$-pairs.*

In particular, by condition $(b)$ the weight filtrations $W^1,\dots,W^r$ are determined by $W^0$ and $N_1,\dots,N_r$. Following [24], we therefore define a *Deligne system* to consist of data $(W,N_1,\dots,N_r;Y^r)$ which generate a system (6.19) satisfying conditions $(a)$–$(e)$ with $W^0=W$. For future use, we define a *pre-Deligne system* to consist of data $(W,N_1,\dots,N_r)$ as above which satisfy conditions $(a)$–$(d)$.

**Lemma 6.24.** *If $(W^0,N_1,\dots,N_r;Y^r)$ is a Deligne system with associated gradings $Y^i$ then $(W^0,N_1,\dots,N_j;Y^j)$ is also a Deligne system.*

**Corollary 6.25.** *Let $T$ be a morphism of Deligne systems* (6.21) *with associated grading $Y^i$ and $\tilde{Y}^i$. Then,*

$$\tilde{Y}^i\circ T = T\circ Y^i \tag{6.26}$$

*Proof.* Equation (6.26) is true for $i=r$ by hypothesis. Likewise, we know that $T$ is a morphism of the Deligne systems

$$(W^{r-1},N_r,Y^r)\to(\tilde{W}^{r-1},\tilde{N}_r,\tilde{Y}^r)$$

By functoriality, this implies (6.26) for $i=r-1$. Accordingly, we have a morphism of Deligne systems

$$(W^0,N_1,\dots,N_{r-1};Y^{r-1})\to(\tilde{W}^0,\tilde{N}_1,\dots,\tilde{N}_{r-1};\tilde{Y}^{r-1})$$

and hence (6.26) holds by downward induction. □

**Lemma 6.27** (Deligne [8])**.** *The set of systems* (6.19) *satisfying conditions* $(a)$–$(e)$ *form an abelian category.*

*Sketch.* Let $T$ be a morphism of Deligne systems (6.21). Then, equation (6.26) holds for all $i$ by the previous Corollary. This forces $T$ to be compatible with all of the associated filtrations and representations of $\mathrm{sl}_2$. We leave the details to the reader. □

Returning to the splitting operation (6.18), suppose that $\theta(z)=\exp(\sum_j z_jN_j)F$ is a polarizable nilpotent orbit of pure Hodge structure of weight $k$ on $V$ with limit mixed Hodge structure $(F,W^r)$. Then, $(W,N_1,\dots,N_r;Y^r)$ is a Deligne system where $Y^r=Y_{(F,W^r)}$ and $W$ pure of weight $k$ on $V$. As shown in [1], in the case where $(F,W^r)$ is split over $\mathbb{R}$, the resulting $\mathrm{sl}_2$-pairs (6.23) generate the representation of $\mathrm{SL}_2^r(\mathbb{R})$ occurring in the $\mathrm{SL}_2$-orbit theorem of Cattani, Kaplan and Schmid [4].

More generally, for any IMHM $(W, N_1, \dots, N_r; F)$ one obtains an representation of $\mathrm{SL}_2^r(\mathbb{R})$ by applying Deligne's construction to the Deligne system

$$(W, N_1, \dots, N_r; Y^r)$$

where $Y^r = \hat{Y}_{(F,W^r)}$. This motivates the following:

*Definition* 6.28 (cf. [16]). A Deligne–Hodge system $(W, N_1, \dots, N_r; F)$ consists of a pre-Deligne system $(W, N_1, \dots, N_r)$ over $K = \mathbb{C}$ equipped with:

(i) A real structure $V = V_{\mathbb{R}} \otimes \mathbb{C}$ to which $W$ and $N_1, \dots, N_r$ descend;
(ii) A decreasing filtration $F$ of $V$;

such that

(f1) $N_j F^p \subset F^{p-1}$ for any $1 \leq j \leq n$ and $p \in \mathbb{Z}$;
(f2) $(F, W^r)$ is a mixed Hodge structure. Furthermore, for $1 \leq k < n$, $w \in \mathbb{Z}$ and for $U = W^k_w$, $(F|_U, W^r|_U)$ is a mixed Hodge structure.

Let $D_r$ denote the category of Deligne systems in $r$ variables over $K = \mathbb{C}$ and $DH_r$ denote the category of Deligne–Hodge systems in $r$-variables. Then, we have a forgetful functor

$$(V_{\mathbb{R}}, W, N_1, \dots, N_r; F) \mapsto (W, N_1, \dots, N_r; Y^r) \tag{6.29}$$

from $DH_r$ to $D_r$.

*Remark* 6.30. In section 2.2.1 of [16], Kato constructs a functor from $DH_r$ to the category $D_{r,\mathbb{R}}$ of Deligne systems over $\mathbb{R}$ by setting $Y^r = \hat{Y}_{(F,W^r)}$.

In Proposition 5.6.2 of [15], Kashiwara shows that the category of IMHM is abelian, and $W$, $M(N_1, W)$, $Gr^W$, $Gr^{M(N_1,W)}$, etc. are exact functors. In Proposition (1.8) of [16], Kato asserts that the category of Deligne–Hodge systems is abelian via the embedding of Theorem (6.1) of Deligne–Hodge systems into the category of IMHM. Since Theorem (6.1) is false, this invalidates Kato's proof.

In §3.2.1 of [16], Kato defines a category of $\mathrm{SL}_2$-orbits. An alternative description of this category is as follows: Let $DH_r \to DH_r$ be the functor defined by the rule

$$(V_{\mathbb{R}}, W, N_1, \dots, N_r; F) \mapsto (V_{\mathbb{R}}, W, \hat{N}_1, \dots, \hat{N}_r; \hat{F}) \tag{6.31}$$

where the $\hat{N}_j$ are the degree zero part of $N_j$ with respect to $\operatorname{ad} Y^{j-1}$ and $(\hat{F}, W^r)$ is the sl$_2$-splitting of $(F, W^r)$. An object of $DH_r$ is an SL$_2$-orbit if it is fixed by this functor. The set $\widehat{DH_r}$ of all SL$_2$-orbits in $DH_r$ is a full subcategory.

*Example* 6.32. $\widehat{DH_0}$ is the category of mixed Hodge structures which are split over $\mathbb{R}$. In the case where $W$ is pure, $\widehat{DH_1}$ is the set of nilpotent orbits with limit mixed Hodge structure split over $\mathbb{R}$. In the case where $W$ is mixed, $\widehat{DH_1}$ consists of admissible nilpotent orbits with limit mixed Hodge structures which are split over $\mathbb{R}$ and $N = N_0$.

**Lemma 6.33.** *If $\hat{\mathbf{V}}$ is an object of $\widehat{DH_r}$ then $\hat{\mathbf{V}}$ is an IMHM.*

*Proof.* This is Proposition 3.2.7 in Kato. □

*Definition* 6.34. Let $\mathbf{V} = (N_1, \ldots, N_r; W, F)$ be an object of $DH_r$ with underlying vector space $V$. If $\mathbf{V}$ is pure of weight $w$ then $Q : V \otimes V \to \mathbb{R}(w)$ polarizes $\mathbf{V}$ if

(a) Each $N_j$ is an infinitesimal isometry of $Q$;
(b) $Q$ polarizes the associated SL$_2$-orbit $\hat{\mathbf{V}}$ obtained by the application of (6.31) to $\mathbf{V}$, i.e. $\hat{\mathbf{V}}$ satisfies the axioms of an IMHM with $Q$ as the polarizing form.

If $\mathbf{V}$ is not pure then $\mathbf{V}$ is said to be grade-polarizable if there exists a polarization for each of the induced Deligne–Hodge systems on $Gr^W$.

6.2. **A Counterexample.** In this section we construct an explicit counterexample to Kato's theorem (6.1) in the case where $\mathbf{V}$ is not graded-polarizable.

Define

$$(V, W, N_1, N_2; F) \tag{6.35}$$

as follows:

— $V$ is the four dimensional real vector space with ordered basis $(e_1, f_1, e_2, f_2)$;

— $N_1$ and $N_2$ are the nilpotent endomorphisms

$$N_1(e_1) = f_1, \qquad N_1(f_1) = 0, \qquad N_1(e_2) = f_2, \qquad N_1(f_2) = 0$$

$$N_2(e_1) = f_2, \qquad N_2(f_1) = 0, \qquad N_2(e_2) = 0, \qquad N_2(f_2) = 0$$

— $W$ is the increasing filtration on $V$, with $\operatorname{gr}^W_k V = 0$ for $k \neq 1$;

— $F$ is the decreasing filtration on $V_{\mathbb{C}}$ with $\operatorname{gr}^k_F V = 0$ for $k \notin [0,1]$ and with $F^1 = \langle e_1, e_2 \rangle$.

**Proposition 6.36.** *The data* (6.35) *defines an object* $\mathbf{V}$ *of* $\mathrm{DH}_2$.

*Proof.* We need to check conditions $(a)$–$(d)$ and $(f)$.

(a) Clearly $N_1N_2 = N_2N_1 = N_1^2 = N_2^2 = 0$. So the the $N_i$ commute and are nilpotent. They also obviously respect $W$ because $W$ is trivial.

(b) Set $W^0 = W$. Then $W^1 := M(N_1, W^0)$ exists and is split by the endomorphism

$$Y^1(e_1) = 2e_1, \qquad Y^1(f_1) = 0, \qquad Y^1(e_2) = 2e_2, \qquad Y^1(f_2) = 0 \tag{6.37}$$

Since $N_2(W^1_k) \subset W^1_{k-2}$, the filtration $W^2 := M(N_2, W^1)$ exists and is equal to $W^1$. In fact, the mixed Hodge structure, $(V, W^2, F)$ is split over $\mathbb{R}$ with $Y^1 = Y^2 = Y_{(W^2,F)}$ the Deligne splitting. ($Y^2v = (p+q)v$ for $v \in I^{(p,q)}_{(W^2,F)}$.) For this Hodge structure, we have $I^{(1,1)} = \langle e_1, e_2\rangle$ and $I^{(0,0)} = \langle f_1, f_2\rangle$.

(c) This is the condition that, if $0 \le k < j-1 \le n-1$ and $U = W^{(k)}_w$ for some $w$, then $W^{(j)}_U$ is the relative monodromy filtration of $N_j|U$ with respect ot $W^{(j-1}$. That is trivial in this case, because we have to have $k = 0$ and $j = 2$, and $W^{(0)}$ is trivial.

(d) Also trivial in this case.

(f.1) This is the horizontality condition which is trivial for the given $F$.

(f.2) The requirement here is that, for $0 \le k < n, w \in \mathbb{Z}$ and $U := W^{(k)}_w$, $(W^{(n)}|_U, F|_U)$ is a mixed Hodge structure. When $k = 0$ we just have to check that $(V, W^{(n)}, F)$ is a mixed Hodge structure. That is fairly obvious. It is also clear when $k = 1$, because $W^1 = W^2$.

□

Suppose now that Theorem 6.1 holds for (6.35). Then there exists an $a \in \mathbb{R}_+$ such that upon setting $N'_1 = N_1$ and $N'_2 = aN_1 + N_2$ the data

$$(V, N'_1, N'_2, W, F) \tag{6.38}$$

underlies an IMHM.

For $z = (z_1, z_2) \in \mathbb{C}^2$, set $N'(z) = \sum_{i=1}^2 z_i N'_i$ and $F'(z) = e^{N'(z)}F$. For $i = 1, 2$, set $e_i(z) = \exp(N'(z))e_i$. Then $F'^1(z) = \langle e_1(z), e_2(z)\rangle$.

Since $(V, W, N_1', N_2', F)$ is an IMHM, there exists a skew-symmetric form

$$Q : V \otimes V \to \mathbb{R}(-1)$$

respecting the $N_i'$ and polarizing $(V, F'(y))$ for $y = (y_1, y_2) \in \mathbb{R}^2$ with $y_1, y_2 \gg 0$. Fix this $Q$. (Here, when we say that $Q$ respects the $N_i'$, we mean that the $N_i'$ are infinitesimal isometries of $Q$.)

**Lemma 6.39.** *Both $N_1$ and $N_2$ respect $Q$.*

*Proof.* Since $N_1 = N_1'$, and $N_1'$ respects $Q$, $N_1$ respects $Q$. So $N_2 = N_2' - aN_1'$ also respects $Q$. □

**Lemma 6.40.** *We have*

(a) $Q(f_1, f_2) = 0$.
(b) $Q(e_1, e_2) = 0$.
(c) $Q(e_2, f_2) = 0$.

*Proof.* (a) Since $N_1$ respects $Q$,

$$0 = Q(N_1 e_1, f_2) + Q(e_1, N_1 f_2) = Q(f_1, f_2) + 0 = Q(f_1, f_2).$$

(b) Set $q(z) := Q(e_1(z), e_2(z))$ and note that $q(z)$ is a polynomial in the variables $z_1, z_2$. Since $Q$ polarizes $F'(y)$ for $y = (y_1, y_2)$ with $y_1, y_2 \gg 0$, $q(y) = 0$ for $y_1, y_2 \gg 0$. But, since $q$ is a polynomial, this implies that $q = 0$ identically. So $q(0) = Q(e_1, e_2) = 0$.

(c) Since $N_2$ respects $Q$, we have $0 = Q(N_2 e_1, e_2) + Q(e_1, N_2 e_2) = Q(f_2, e_2)$. So $Q(e_2, f_2) = 0$ as well. □

**Proposition 6.41.** *Contrary to Lemma* (6.40)*(c), $Q(f_2, e_2) > 0$. Consequently, the Deligne system $(V, N_1', N_2', W, F)$ from* (6.38) *does not admit a form $Q$ making it into a* IMHM.

*Proof.* Suppose $(V, N_1', N_2', W, F)$ is an IMHM polarized by $Q$. Assume $y_1, y_2 \gg 0$ so that $(V, F(y))$ is a pure Hodge structure of weight 1 polarized by $Q$. Moreover, assume $y_1$ and $y_2$ are positive.

Let $C(y)$ denote the Weil operator on the pure Hodge structure $(V, F(y))$. It is given on $V^{pq}$ by multiplication by $i^{p-q}$. (We are using the sign conventions from §1.2 of Kashiwara [15].)

Since $Q$ is a polarization on $(V, F(y))$, the form $Q(Cu, \bar{v})$ is positive definite. It follows that $Q(Ce_2(y), \bar{e}_2(y)) > 0$. So we compute

$$\begin{aligned} 0 < Q(Ce_2(y), \bar{e}_2(y)) &= Q(C(e_2 + i(y_1 + ay_2)f_2), e_2 - i(y_1 + ay_2)f_2) \\ &= iQ(e_2 + i(y_1 + ay_2)f_2, e_2 - i(y_1 + ay_2)f_2) \\ &= 2(y_1 + ay_2)Q(e_2, f_2) \end{aligned}$$

Since $y_1, y_2$ and $a$ are positive, it follows that $Q(e_2, f_2) > 0$. □

6.3. **Categorical Comments.** An IMHM $\mathbf{V}$ is split if $\mathbf{V} \cong \oplus Gr_k^W \mathbf{V}$. The category of split infinitesimal Hodge modules is, more or less by definition, polarizable Tannakian in the sense of Saavedra-Rivano (see, for example, page 169 of [9]). It follows that the category IMHM$^s$ of split infinitesimal mixed Hodge modules is semi-simple (by Proposition 4.11 on page 169 of [9]).

However, we have the following.

**Proposition 6.42.** *The example in* (6.35) *is not a semi-simple object.*

*Proof.* Let $H = \langle e_2, f_2 \rangle$. By restriction of $(W, N_1, N_2, F)$ to $H$ we obtain a sub-object $\mathbf{H}$ of $\mathbf{V}$. (The restriction of $N_2$ to $H$ is zero). But this sub-object is clearly not a direct summand. □

The proposition gives another way to see that Theorem (6.1) fails for $\mathbf{V}$: Since split objects in IMHM are semi-simple, the theorem would imply that split objects in $DH_r$ are also semi-simple.

6.4. **Graded Polarizability.** Every graded-polarizable Deligne–Hodge system gives rise to an IMHM:

**Theorem 6.43.** *Let $(V, W, N_1, \ldots, N_n, F)$ be a graded polarizable* DH *system of $n$ variables. Then for $N'_j = \sum_{k=1}^{j} a_{j,k} N_k$ $(1 \le j \le n)$ with $a_{j,k} > 0$ $(1 \le k \le j \le n)$ such that $a_{j,k}/a_{j,k+1} \gg 0$ $(1 \le k < j \le n)$, $(V, W, N'_1, \ldots, N'_n, F)$ is an IMHM of $n$ variables.*

Essentially, Kato's proof goes though word for word upon addition of this polarizability hypothesis: Suppose $(Q_w)$ is a fixed polarization of the graded of $\mathbf{V}$. Write $D$ for the mixed period domain corresponding to $Q_w$. It sits in the so-called compact dual $\check{D}$. Write $G$ for the group of isometries of $Q$ preserving $W$. Then $G(\mathbb{C})$ acts on the algebraic variety $\check{D}$, while $G(\mathbb{R})$ acts on the open complex submanifold $D$. The $N_i$ and $H_i$ are all in the Lie algebra $\mathfrak{g}$ of $G$ and, most importantly, the function $\beta$ constructed by Kato lies in $G(\mathbb{R})$. On the other hand, note that, without $Q$, there really is no period domain $D$ (or $\check{D}$).

Now, Kato's proof shows that $\beta(y)F(y)$ converges to

$$I = \exp(\sum_j i\hat{N}_j)\hat{F}$$

as long as $y_i/y_{i+1} \gg 0$ for all $i$. In particular, $\beta(y)F(y)$ lies in $D$ for such $y$. Since $\beta(y) \in G(\mathbb{R})$, it follows that $F(y)$ lies in $D$ as well.

6.5. **Geometric Structure.** By Theorem (6.43), the question of when a given Deligne–Hodge system gives rise to an IMHM for an appropriate substitution $N_j \to N_j'$ reduces to a question about the polarizability of the underlying $\mathrm{SL}_2$-orbit. Accordingly, we make the following definition:

*Definition* 6.44. The $\mathrm{sl}_2$-type of a Deligne system $(W, N_1, \dots, N_r, Y^r)$ consists of the weight filtration $W$ and the associated $\mathrm{sl}_2$-pairs $(\hat{N}_1, H_1), \dots, (\hat{N}_r, H_r)$ of (6.23).

*Remark* 6.45. In the case where $W$ is pure of weight $k$ the sum

$$Y^r = k\,\mathrm{Id} + H_1 + \cdots + H_r \tag{6.46}$$

for any possible associated Deligne system. In particular, we can then recover the intermediate gradings $Y^j$ and hence the elements $H_j = Y^j - Y^{j-1}$ via the iterated application of Deligne construction $Y^j = Y(\hat{N}_j, Y^{j+1})$. Thus, in the pure case, an $\mathrm{sl}_2$-type is equivalent to $(W, \hat{N}_1, \dots, \hat{N}_r, Y^r)$. In [16] (cf. Prop. 3.3.2), Kato calls such Deligne systems an associated SL(2)-orbit.

The remainder of this section is devoted to proving that the set all Deligne systems with a given $\mathrm{sl}_2$-type forms an algebraic variety. We start with a series of lemmata:

**Lemma 6.47.** *Let $(W, \hat{N}, H)$ be an* $\mathrm{sl}_2$*-type. Let $\mathcal{Y}(W, \hat{N}, H)$ be the set of all gradings of $W$ which commute with $\hat{N}$ and $H$. Define* $\mathrm{gl}_{-1}(W, \hat{N}, H)$ *to be the subalgebra of* $\mathrm{gl}_{-1}(W)$ *consisting of elements which commute with $\hat{N}$ and $H$. Then, $\mathcal{Y}(W, \hat{N}, H)$ is an affine space upon which the subgroup* $\exp(\mathrm{gl}_{-1}(W, \hat{N}, H))$ *acts simply transitively.*

*Proof.* The set $\mathcal{Y}(W)$ of all gradings of $W$ is an affine space upon which the subalgebra $W_{-1}\mathrm{gl}(V)$ acts simply transitively by $Y \mapsto Y + \beta$. Accordingly, let $Y \in \mathcal{Y}(W, \hat{N}, H)$ and $\beta \in W_{-1}\mathrm{gl}(V)$. Then, clearly $Y + \beta \in \mathcal{Y}(W, \hat{N}, H)$ if and only if $\beta \in \mathrm{gl}_{-1}(W, \hat{N}, H)$. Thus, $\mathcal{Y}(W, \hat{N}, H)$ is an affine space upon which $\mathrm{gl}_{-1}(W, \hat{N}, H)$ acts simply transitively by $Y \mapsto Y + \beta$.

On the other hand, as discussed in Proposition (2.2) of [4] the group $\exp(W_{-1}\mathrm{gl}(V))$ also acts simply transitively on $\mathcal{Y}(W)$ via the adjoint action. Moreover, the relation between $\alpha$, $\beta \in W_{-1}\mathrm{gl}(V)$ defined by the equation

$$e^{\mathrm{ad}\,\alpha} Y = Y + \beta$$

is given by universal Lie polynomials in the eigencomponents of $\alpha$ and $\beta$ with respect to $\mathrm{ad}\, Y$. In particular, since $Y$ commutes with $\hat{N}$ and $H$, if $\beta$ also commutes with $\hat{N}$ and $H$ then so do all of its eigencomponents, and hence $\alpha$ also has this property. □

Let $(W, \hat{N}, H)$ be an $\mathrm{sl}_2$-type. Define $\mathcal{U}(W, \hat{N}, H)$ to be the set of pairs $(Y^0, N_-)$ where $Y^0 \in \mathcal{Y}(W, \hat{N}, H)$ and

$$N_- = \sum_{k \geq 2} N_{-k}$$

where $N_{-k}$ is either 0 or an element of $E_{-k}(\mathrm{ad}\, Y^0)$ which is of highest weight $k - 2$ for the representation of $\mathrm{sl}_2$ generated by $(\hat{N}, H)$. Let

$$\pi : \mathcal{U}(W, \hat{N}, H) \to \mathcal{Y}(W, \hat{N}, H)$$

denote projection $(Y^0, N_-) \to Y^0$.

**Lemma 6.48.** *$\pi : \mathcal{U}(W, \hat{N}, H) \to \mathcal{Y}(W, \hat{N}, H)$ is an equivariant vector bundle in the sense that for any $\alpha \in \mathrm{gl}_{-1}(W, \hat{N}, H)$:*

$$e^{\alpha} : \pi^{-1}(Y^0) \to \pi^{-1}(e^{\alpha}.Y^0)$$

*is a linear isomorphism on the fibers, and the group* $\exp(\mathrm{gl}_{-1}(W, \hat{N}, H))$ *acts transitively on the base.*

For future use, we record the following:

**Lemma 6.49.** *Let $W$ and $W'$ be increasing gradings of a finite dimensional vector space $V$ over a field of characteristic zero, with respective gradings $Y$ and $Y'$. If $[Y, Y'] = 0$ then $Y$ preserves $W'$.*

*Proof.* Mutually commuting semisimple endomorphisms can be simultaneously diagonalized over any field which contains all the eigenvalues of both endomorphisms. Since $Y$ and $Y'$ have integral eigenvalues

$$V = \oplus_{p,q} V(p, q), \qquad V(p, q) = E_p(Y) \cap E_q(Y')$$

Therefore, $Y$ preserves $W'_k = \oplus_{q \leq k} V(p, q)$. □

**Lemma 6.50.** *Let $\mathcal{S}$ denote the set of all Deligne systems with* $\mathrm{sl}_2$*-type $(W, \hat{N}, H)$. Then, the map $\Psi : \mathcal{S} \to \mathcal{U}(W, \hat{N}, H)$ which sends the Deligne system $(W, N, Y^1)$ to $(Y^0, N_-)$ where $Y^0 = Y(N, Y^1)$ and $N_- = \sum_{k \geq 2} N_{-k}$ is the sum of components of $N$ of negative weight with respect to* $\mathrm{ad}\, Y^0$ *is a bijection.*

*Proof.* A Deligne system with $\mathrm{sl}_2$-data $(W, N, H)$ produces an element of $\mathcal{U}(W, \hat{N}, H)$ via $\Psi$. Moreover, $\Psi$ is injective since a Deligne system with $\mathrm{sl}_2$-type $(W, \hat{N}, H)$ is recovered from its image under $\Psi$ by the rule:

$$N = \hat{N} + N_-, \qquad Y^1 = Y^0 + H \tag{6.51}$$

It remains to prove that $\Psi$ is surjective. Let $(Y^0, N_-)$ be an element of $\mathcal{U}(W, \hat{N}, H)$, and define $N$ and $Y^1$ by (6.51). It is sufficient to show that $(W, N, Y^1)$ is a Deligne system, since by construction $\Psi(W, N, Y^1) = (Y^0, N_-)$.

Accordingly, let $W^1$ be the weight filtration determined by $Y^1$. We need to check that $(W, N, Y^1)$ satisfies Deligne system axioms $(a)$–$(e)$.

(e) By construction, $[Y^1, Y^0] = 0$ and hence $Y^1$ preserves $W^1$ and $W^0$ by the previous Lemma. Likewise,

$$[Y^1, N] = [Y^1, \hat{N}] + \sum_{k \geq 2} [Y^0 + H, N_{-k}] = -2N \tag{6.52}$$

(a) By $(e)$, $N$ is nilpotent. Similarly, since $N$ has weights less than or equal to zero with respect to $\mathrm{ad}\, Y^0$, it preserves $W^0$.

(b) To verify that $W^1 = M(N, W^0)$, we begin with the assertion that $W^1 = M(\hat{N}, W^0)$: By construction,

$$[Y^1, \hat{N}] = [Y^0 + H, \hat{N}] = -2\hat{N} \tag{6.53}$$

so $\hat{N}$ lowers $W^1$ by 2. It remains to show that

$$\hat{N}^\ell : Gr^{W^1}_{k+\ell} Gr^{W^0}_k \to Gr^{W^1}_{k-\ell} Gr^{W^0}_k \tag{6.54}$$

is an isomorphism. However, since $H = Y^1 - Y^0$ and $[Y^1, Y^0] = 0$ it follows that

$$Gr^{W^1}_{k+\ell} Gr^{W^0}_k \cong E_{k+\ell}(Y^1) \cap E_k(Y^0) = E_\ell(H) \cap E_k(Y^0)$$

Accordingly, (6.54) is an isomorphism. Consequently, $M(N, W) = M(\hat{N}, W)$ since $[Y^1, N] = -2N$ and and changing $\hat{N} \to N$ does not change the induced action of $N$ on $Gr^W$.

(c) Since $[Y^1, Y^0] = 0$, it follows that the argument of $(b)$ above holds for the restriction of $N$ to $W^0_\ell$.

(d) By $(a)$+$(e)$ above, $N$ preserves $W^0$ and lowers the weights of $W^1$ by 2.

□

**Theorem 6.55.** *The set $\mathcal{S}$ of all $r$-variable Deligne systems $(W, N_1, \dots, N_r, Y^r)$ of* $\mathrm{sl}_2$*-type*

$$(W, \hat{N}_1, H_1, \dots, \hat{N}_r, H_r) \tag{6.56}$$

*is an algebraic variety.*

*Proof.* We induct on the number of variables. The case $r = 1$ is covered by the previous Lemma.

To continue, we observe that if $S = (W, N_1, \dots, N_r; Y^r)$ is a point of $\mathcal{S}$ then

$$W^j = M(\hat{N}_j, W^{j-1}), \qquad j = 1, \dots, r$$

and hence the weight filtrations $W^j$ attached to $S$ are determined by the $\mathrm{sl}_2$-data (6.56). Let $\mathcal{Y}$ denote the affine variety consisting of gradings of $W^{r-1}$.

Let $\mathcal{S}'$ denote the set of Deligne systems with $\mathrm{sl}_2$-data $(W, \hat{N}_1, H_1, \dots, \hat{N}_{r-1}, H_{r-1})$. Let $\mathcal{S}''$ denote the set of Deligne systems with $\mathrm{sl}_2$-data $(W^{r-1}, \hat{N}_r, H_r)$. By the induction hypothesis, both $\mathcal{S}'$ and $\mathcal{S}''$ are algebraic varieties, and hence so is the product $\mathcal{S}' \times \mathcal{S}''$. Let $\mathcal{P} \subset \mathcal{S}' \times \mathcal{S}''$ be the fiber product over $\mathcal{Y}$ defined by the maps

$$S' = (W, N_1, \dots, N_{r-1}, Y^{r-1}) \mapsto Y^{r-1}$$
$$S'' = (W^{r-1}, N_r, Y^r) \mapsto Y'(N_r, Y^r)$$

A point $S \in \mathcal{S}$ determines a point $(S', S'') \in \mathcal{P}$ by the rule

$$S' = (W, N_1, \dots, N_{r-1}, Y^{r-1}), \qquad S'' = (W^{r-1}, N_r, Y^r) \tag{6.57}$$

such that:

(I) $Y^r$ commutes with the gradings $Y^{r-1}, \dots, Y^0$ attached to $S'$;
(II) $N_r$ has only non-positive eigenvalues relative to $\operatorname{ad} Y^j$ for $0 \leq j \leq r-1$;
(III) $[N_r, N_j] = 0$ for $1 \leq j \leq r-1$;
(IV) $[Y^r, N_j] = -2N_j$ for $1 \leq j \leq r$;

Conversely, let $(S', S'') \in \mathcal{P}$ be a point of the form (6.57) which satisfies conditions $(I)$–$(IV)$. Then,

$$S = (W, N_1, \dots, N_r; Y^r) \tag{6.58}$$

is a Deligne system:

(a) Since each $N_j$ is part of a Deligne system, it is nilpotent. Similarly, since $(N_1, \dots, N_{r-1})$ are part of a Deligne system, they mutually commute. Condition $(III)$ implies the remaining commutativity conditions. Likewise, by hypothesis, $N_1, \dots, N_{r-1}$ preserve $W^0$. Condition $(II)$ implies that $N_r$ preserves $W^0$.
(b) The fact that $S'$ and $S''$ are Deligne systems implies the existence of all required relative weight filtrations.
(c) Since $S'$ and $S''$ are Deligne systems, this condition is automatic except for the extremal case: $W^r$ restricted to $U = W^k_\ell$ for $k < r-1$ is the relative weight filtration of $W^{r-1}|_U$ and $N_r|_U$. This follows from properties $(I)$ and $(II)$.

(d) Since $S'$ and $S''$ are Deligne systems, the only unresolved cases are $N_r(W^k_\ell) \subset W^k_\ell$ for $k < r-1$ (which follows from $(II)$) and $N_j(W^r_\ell) \subset W^r_{\ell-2}$ which follows from $(IV)$.

(e) $[Y^r, N_j] = -2N_j$ is property $(IV)$. Property $(I)$ implies that $Y^r$ preserves each $W^j$.

Let $\mathcal{A}$ be the algebraic subvariety of $\mathcal{P}$ defined by properties $(I)$–$(IV)$. Given $(S', S'') \in \mathcal{A}$, the corresponding Deligne system $S$ has the same $\mathrm{sl}_2$-data as $S'$ and $S''$, i.e. $S \in \mathcal{S}$. A simple check shows that the maps

$$\mathcal{S} \to \mathcal{A} \to \mathcal{S}, \qquad \mathcal{A} \to \mathcal{S} \to \mathcal{A}$$

are the identity, and hence $\mathcal{S}$ is isomorphic to the algebraic variety $\mathcal{A}$.

□

In the case where $(\hat{N}_1, H_1), \ldots, (\hat{N}_r, H_r)$ are all infinitesimal isometries of bilinear forms on $Gr^W$, we can ask that all $N_j$'s appearing above are also infinitesimal isometries. This is again an algebraic condition. More generally, if the initial $\mathrm{sl}_2$-data belongs a Mumford–Tate Lie algebra $\mathfrak{m}$ then it is an algebraic condition for all the $N_j$'s to belong to $\mathfrak{m}$. Likewise, the condition that each $N_j$ be horizontal with respect to a given Hodge filtration $F$ is also algebraic.

As independent check of the compatibility of the next few examples with the results of the earlier sections of this paper, we develop our examples starting from Lemma (6.24) in Schmid [23]. Namely, if $H_{\mathbb{C}} = H_{\mathbb{R}} \otimes \mathbb{C}$ is Hodge structure of weight $k$ equipped with a horizontal action of $\mathrm{sl}_2(\mathbb{C})$ then $H_{\mathbb{C}}$ is a direct sum of irreducible representations of the form $S(n) \otimes \mathbb{R}(m)$ and $S(n) \otimes E(p,q)$ where

— $S(n) = \mathrm{Sym}^n(\mathbb{C}^2)$ where $\mathbb{C}^2$ is equipped with the standard matrix action of $\mathrm{sl}_2$ with respect to the basis $e = (1,0)$ and $f = (0,1)$ of $\mathbb{C}^2$. Relative to the limit mixed Hodge structure $e$ is type $(1,1)$ and $f$ is type $(0,0)$. The polarization is defined by

$$Q(e^j f^{n-j}, e^{n-j} f^j) = (-1)^n (-1)^j j!(n-j)!$$

— $E(p,q) = \mathbb{C}^2$ equipped with the trivial action of $\mathrm{sl}_2$ and $e + if$ of type $(p,q)$.

— $\mathbb{R}(m)$ is $\mathbb{C}$ equipped with a pure Hodge structure of type $(-m,-m)$.

Before proceeding with the examples, we note that in the case where $W$ is pure of weight $k$, equation (6.46) applies and we shall omit $W$ from the data of the Deligne system. If $W$ pure of weight $k$, we also always have $N_1 = \hat{N}_1$.

*Example* 6.59. Let $(N_1, N_2; F)$ generate a pure nilpotent orbit of odd weight $2m+1$ of "vanishing cycle" type, i.e. there exist linearly independent elements $\alpha_1$ and $\alpha_2$ of the underlying real vector space $V_{\mathbb{R}}$ such that such that $Q(\alpha_1, \alpha_2) = 0$ and

$$N_j(\gamma) = Q(\gamma, \alpha_j)\alpha_j.$$

Let $W^j = W(\sum_{i \le j} N_i)[-(2m+1)]$ and assume without loss of generality that $(F, W)$ is split over $\mathbb{R}$. In Schmid's terminology, the corresponding horizontal $\mathrm{sl}_2$ action generated by $N = N_1 + N_2$ and the limit mixed Hodge structure $(F, W^2)$ decomposes the underlying vector space as

$$[S(1) \otimes \mathbb{R}(-m)] \oplus [S(1) \otimes \mathbb{R}(-m)] \oplus K$$

where $K$ is a pure Hodge structure of weight $2m+1$ with trivial $\mathrm{sl}_2$-action. The two $S(1)$ factors correspond to the isotypical components of highest weight 1 for $(N_1, H_1)$ and $(N_2, H_2)$.

More concretely, there exist dual elements $\alpha_1'$ and $\alpha_2'$ of type $(m+1, m+1)$ with respect to $(F, W^2)$ such that $Q(\alpha_j', \alpha_k) = \delta_{jk}$, and $(F, W^2)$ has Deligne bigrading

$$I^{m+1,m+1} = \mathrm{span}(\alpha_1', \alpha_2'), \qquad \bigoplus_p I^{p,2m+1-p} = K, \qquad I^{m,m} = \mathrm{span}(\alpha_1, \alpha_2).$$

The associated $\mathrm{sl}_2$-data is given by $(N_1, H_1)$ and $(N_2, H_2)$ where

$$\begin{array}{lll} E_1(H_1) = \ \mathrm{span}(\alpha_1'), & E_0(H_1) = K \oplus \mathrm{span}(\alpha_2, \alpha_2'), & E_{-1}(H_1) = \mathrm{span}(\alpha_1) \\ E_1(H_2) = \ \mathrm{span}(\alpha_2'), & E_0(H_2) = K \oplus \mathrm{span}(\alpha_1, \alpha_1'), & E_{-1}(H_2) = \mathrm{span}(\alpha_2) \end{array}.$$

To analyze the corresponding variety of Deligne systems via the process described above, we start with the set $\mathcal{S}$ attached to $W^0$ and $(N_1, H_1)$. Since $W^0$ is pure of weight $k$, $\mathcal{S}$ is a point. This becomes the set $\mathcal{S}'$ at the next step, and we have to consider the set $\mathcal{S}''$ attached to $W^1$ and $(N_2, H_2)$. By (6.46), $Y^1 = Y^0 + H_1$ which eliminates the freedom to pick a point in $\mathcal{Y}(W^1, N_2, H_2)$. The remaining freedom in $\mathcal{U}(W^1, N_2, H_2)$ is to select an element which is weight $-2$ for $\mathrm{ad}\, Y^1$ and highest weight 0 for $(N_2, H_2)$. This is exactly the space spanned by $N_1$, and so the possible set of Deligne systems consists of the triples

$$(N_1, N_2 + aN_1, Y^2) \tag{6.60}$$

where $a$ is an arbitrary scalar.

*Example* 6.61. More generally, let $(N_1, N_2; F)$ generate a two variable $\mathrm{sl}_2$-orbit of weight $k$ polarized by $Q$, i.e. the limit mixed Hodge structure $(F, W^2)$ split over $\mathbb{R}$ and $\hat{N}_j = N_j$. Let $Y^1 = k\,\mathrm{Id} + H_1$ and $Y^2 = Y^1 + H_2$. Then, the set of Deligne systems with the same $\mathrm{sl}_2$-data consists of all triples

$$(N_1, N_2 + \eta, Y^2)$$

where $\eta = \sum_{\ell \geq 2} \eta_{-\ell}$ and $\eta_{-\ell}$ is weight $-\ell$ with respect to $\mathrm{ad}(Y^1)$ and highest weight $\ell - 2$ for $(N_2, H_2)$.

In particular, let $\Omega$ denote the set of all $(-1,-1)$-morphisms of $(F, W^2)$ which are lowest weight $-2$ for $(N_1, H_1)$ and highest weight zero for $(N_2, H_2)$. Note that $\Omega \neq 0$ since $N_1 \in \Omega$. Fix an inner product on $\Omega$. Recall that $F_o = e^{iN_1 + iN_2} F$ belongs to the classifying space $D$. Moreover, since $D$ is an open subset of the compact dual $\check{D}$ in the analytic topology, there exists an real number $t > 0$ such that if $\omega \in \Omega$ has norm 1 and $|\tau| < t$ then $e^{i\tau\omega} F_o$ also belongs to $D$. The following argument shows that $(N_1, N_1 + t\omega, N_2; F)$ generates a nilpotent orbit: By assumption, the map

$$\theta(z_1, z_2, z_3) = \exp(z_1 N_1 + z_2(N_1 + t\omega) + z_3 N_2) F \tag{6.62}$$

is horizontal. Moreover $\omega$ commutes with $N_1$ and $N_2$ since it is lowest weight $-2$ for $(N_1, H_1)$ and highest weight zero for $(N_2, H_2)$.

To see that $\theta$ takes values in $D$ when the imaginary parts of $z_j = x_j + \sqrt{-1} y_j$ are sufficiently large, observe that

$$y_1 N_1 + y_2(N_1 + t\omega) + y_3 N_2 = (y_1 + y_2)\,(N_1 + \tau\omega) + y_3 N_2$$

where $\tau = (y_2 t)/(y_1 + y_2)$ is positive and less than $t$. Therefore, since $\omega$ is lowest weight $-2$ for $(N_1, H_1)$ and highest weight zero for $(N_2, H_2)$ it follows that

$$\begin{aligned} \exp(iy_1 N_1 + iy_2(N_1 + t\omega) + iy_3 N_2) F &= \exp(i(y_1 + y_2)(N_1 + \tau\omega) + iy_3 N_2) F \\ &= (y_1 + y_2)^{-H_1/2} y_3^{-H_2/2} e^{i\tau\omega} F_o \end{aligned}$$

Accordingly, since $(y_1 + y_2)^{-H_1/2} y_3^{-H_2/2}$ is a real automorphism of $Q$ and $e^{i\tau\omega} F_o \in D$ it follows that (6.62) takes values in $D$ provided the imaginary parts of $z_1$, $z_2$ and $z_3$ are positive. Likewise, if $\mathfrak{a}$ is an abelian subalgebra of $\Omega$ and $\omega_1, \dots, \omega_\ell$ are an orthonormal basis of $\mathfrak{a}$ then a similar type of argument shows that

$$(N_1, N_1 + t_1\omega_1, \dots, N_1 + t_\ell\omega_\ell, N_2; F)$$

generates a nilpotent orbit, provided $|t_j| < t$ for all $j$.

*Example* 6.63. Example (6.59) can be generalized to the weight $2m+1$ case where there exist pairwise orthogonal, linearly independent sets of vanishing cycles $\{\alpha_{11}, \dots, \alpha_{1p}\}$ and $\{\alpha_{21}, \dots, \alpha_{2q}\}$ such that $N_i = \frac{1}{2}\sum_\ell \omega_i^{\ell\ell}$ where

$$\omega_i^{jk}(u) = Q(u, \alpha_{ij})\alpha_{ik} + Q(u, \alpha_{ik})\alpha_{ij}$$

for $j \leq k$. Then, the associated set of Deligne systems is

$$\mathcal{S}(\mathbb{F}) = \{(N_1, N_2 + \sum_{j \leq k} c_{jk}\omega_1^{jk}, Y^2) \mid c_{jk} \in \mathbb{F}\}$$

where $\mathbb{F} \subseteq \mathbb{C}$ is the field of interest.

In Schmid's terminology, the associated representation of $\mathrm{sl}_2$ for $N = N_1 + N_2$ is a sum

$$[S(1) \otimes \mathbb{R}(-m)]^p \oplus [S(1) \otimes \mathbb{R}(-m)]^q \oplus K$$

where two groupings of $S(1)$ factors corresponding to the isotypical components of highest weight 1 for $(N_1, H_1)$ and $(N_2, H_2)$, and $K$ is a pure Hodge structure of weight $2m+1$ with trivial $\mathrm{sl}_2$-action. As in Example (6.59), there are dual elements

$$\alpha'_{11}, \dots, \alpha'_{1p}, \alpha'_{21}, \dots, \alpha'_{2q}$$

of type $(m+1, m+1)$ such that $Q(\alpha'_{ij}, \alpha_{k\ell}) = \delta_{ik}\delta_{j\ell}$. Let $A_i = \mathrm{span}(\alpha_{ij})$ and $A'_i = \mathrm{span}(\alpha'_{ij})$. Then, Deligne bigrading of $(F, W^2)$ is given by

$$I^{m+1,m+1} = A'_1 \oplus A'_2, \qquad \bigoplus_p I^{p,2m+1-p} = K, \qquad I^{m,m} = A_1 \oplus A_2.$$

The corresponding $\mathrm{sl}_2$-data is $(N_1, H_1)$ and $(N_2, H_2)$ where

$$E_1(H_1) = A'_1, \qquad E_0(H_1) = K \oplus A_2 \oplus A'_2, \qquad E_{-1}(H_1) = A_1$$

$$E_1(H_2) = A'_2, \qquad E_0(H_2) = K \oplus A_1 \oplus A'_1, \qquad E_{-1}(H_2) = A_2$$

Due to the short length of the monodromy weight filtrations, it follows by Example (6.61) we are looking for infinitesimal isometries which are lowest weight $-2$ for $(N_1, H_1)$ and highest weight zero for $(N_2, H_2)$. This subspace is spanned by the elements $\omega_1^{ij}$ for $i \leq j$.[8]

The short length of the monodromy weight filtrations also forces commutativity of all $\omega_1^{ij}$, so we can form a new monodromy cone using the techniques described at the end of Example (6.61) – just pick an inner product which makes all the $\omega_1^{ij}$ orthonormal and find the appropriate value of $t$.

*Example* 6.64. Turning now to §5, we consider a two variable example of the form

$$S(2) \oplus S(2) \oplus \mathbb{R}(-1)$$

where the two $S(2)$ factors are the isotypical components of highest weight 2 for $(N_1, H_1)$ and $(N_2, H_2)$, and $\mathbb{R}(-1)$ is a factor of Hodge type $(1,1)$ on which both copies of $\mathrm{sl}_2$ act trivially. The corresponding period domain has Hodge numbers $h^{2,0} = 2$ and $h^{1,1} = 3$.

More concretely, we start with a real vector space $V_{\mathbb{R}}$ with basis

$$\{\alpha_2, \alpha_1, \alpha_0, \beta_2, \beta_1, \beta_0, \gamma\}$$

in which we think of $\alpha_j = e^j f^{2-j}$ and $\beta_j = e^j f^{2-j}$ under the identification with $S(2) = \mathrm{Sym}^2(\mathbb{C}^2)$, and $\gamma$ is the generator of $\mathbb{R}(-1)$. Accordingly, $N_1$ annihilates $\{\beta_2, \beta_1, \beta_0, \gamma\}$ and acts on $\alpha_j = e^j f^{2-j}$ by the rule $N_1(\alpha_j) = j\alpha_{j-1}$. Likewise, $N_2$ annihilates $\{\alpha_2, \alpha_1, \alpha_0, \gamma\}$ and acts on $\beta_j = e^j f^{2-j}$ by the rule $N_2(\beta_j) = j\beta_{j-1}$. The polarizing form is given by

$$Q(\alpha_j, \alpha_{2-j}) = Q(\beta_j, \beta_{2-j}) = (-1)^j j!(2-j)!, \qquad Q(\gamma, \gamma) = 1$$

and all other pairings zero.

The limit Hodge filtration of $(F, W^2)$ is

$$I^{2,2} = \mathrm{span}(\alpha_2, \beta_2), \qquad I^{1,1} = \mathrm{span}(\alpha_1, \beta_1, \gamma), \qquad I^{0,0} = \mathrm{span}(\alpha_0, \beta_0)$$

[8]Sketch: One first checks that the stated conditions imply that such an element $\gamma$ vanishes on $A_2' \oplus K \oplus A_1 \oplus A_2$ and $\gamma(A_1') \subset A_1$. Let $\gamma(\alpha_{1\ell}') = \sum_k \gamma_\ell^k \alpha_{1k}$. Then, the infinitesimal isometry condition forces $\gamma_j^i = \gamma_i^j$.

which is type $(V)$ in the setting of §5. The nilpotent orbit $\theta(z_1) = e^{z_1 N_1} e^{iN_2} F$ is of type $(II)$ with limit Hodge numbers $h^{2,2} = h^{2,0} = h^{0,2} = 1$, $h^{1,1} = 3$ and $h^{0,0} = 1$.

The corresponding elements $H_1$ and $H_2$ are given by

$$E_2(H_2) = \text{span}(\beta_2), \qquad E_0(H_2) = \text{span}(\alpha_2, \alpha_1, \alpha_0, \beta_1, \gamma), \qquad E_{-2}(H_2) = \text{span}(\beta_0)$$

$$E_2(H_1) = \text{span}(\alpha_2), \qquad E_0(H_1) = \text{span}(\beta_2, \beta_1, \beta_0, \alpha_1, \gamma), \qquad E_{-2}(H_1) = \text{span}(\alpha_0)$$

To find the possible candidate deformations $N_2 \mapsto N_2 + \eta$ which preserve the underlying $\text{sl}_2$-orbit structure, we can start by first identifying all morphisms of type $(-1,-1)$ of $(F, W^2)$ which are infinitesimal isometries. This space has basis $\{\eta_1, \eta_2, \eta_3, \eta_4, N_1, N_2\}$ where

$$\eta_1(\alpha_2) = \beta_1, \qquad \eta_1(\beta_1) = \frac{1}{2}\alpha_0$$

$$\eta_2(\alpha_2) = \gamma, \qquad \eta_2(\gamma) = -\frac{1}{2}\alpha_0$$

$$\eta_3(\beta_2) = \alpha_1, \qquad \eta_3(\alpha_1) = \frac{1}{2}\beta_0$$

$$\eta_4(\beta_2) = \gamma, \qquad \eta_4(\gamma) = -\frac{1}{2}\beta_0$$

and annihilate all other basis elements. A short calculation shows that only $\eta_2$ and $\eta_4$ commute with $N_1$ and $N_2$ which is required for $(N_1, N_2 + \eta, Y^2)$ to be a Deligne system. However, $\eta_4$ turns out to be a lowest weight vector of weight $-2$ for $(N_2, H_2)$ whereas for a deformation of Deligne systems $(N_1, N_2 + \eta, Y^2)$, we would need $\eta$ to be a sum of highest weight vectors for $(N_2, H_2)$. On the other hand, $\eta_2$ is lowest weight $-2$ for $(N_1, H_1)$ and highest weight zero for $(N_2, H_2)$. Thus, the set of all Deligne systems with given $\text{sl}_2$-data in this case is

$$(N_1, N_2 + a\eta_2 + bN_1, Y^2), \qquad Y^2 = 2\text{Id} + H_1 + H_2$$

where $a$ and $b$ are arbitrary scalars.

As in Example (6.61), we obtain an associated 3 variable nilpotent orbit with data

$$(N_1, N_1 + t\eta_2, N_2; F)$$

Looking back to §5, it was predicted that for a degeneration of type $(V)$ with Hodge numbers $h^{2,0} = 2$ and $h^{1,1} = 3$ the maximum possible dimension of a nilpotent cone is 3, which is realized by this example.

*Example* 6.65. As a "degenerate" case of the above, observe that we can always augment a pure, 1-variable Deligne system $(N, Y^1)$ of weight $k$ to a 2-variable Deligne system $(N_1, N_2, Y^2)$ by setting $N_1 = N$, $N_2 = 0$ and $Y^2 = Y^1$. In this case, one finds that any Deligne system deformation

$$(N_1, N_2 + \eta, Y^2)$$

of $(N_1, N_2, Y^2)$ must have $\eta$ of lowest weight $-2$ for the representation $(N_1, H_1)$.

In the case where $(N, Y^1)$ arises from a pure nilpotent orbit $\theta(z) = e^{zN}F$ with limit mixed Hodge structure $(F, W^1)$ split over $\mathbb{R}$, the method of Example (6.61) shows that if $\eta$ is also $(-1,-1)$-morphism of $(F, W^1)$ then there exists a positive number $t$ such that $(N, N + t\eta, F)$ generates a pure nilpotent orbit.

A simple example of this type occurs when we start with a several variable pure nilpotent orbit generated by $(N_1, \dots, N_r; F)$ with limit mixed Hodge structure split over $\mathbb{R}$, and consider the 1-variable orbit $\theta$ generated by $N = \sum_j N_j$ and $F$. Then, all of the $N_j$'s are lowest weight $-2$ for the $\mathrm{sl}_2$-representation attached to $\theta$.

In particular, let $\theta(z) = e^{zN}F$ be a 1-variable pure nilpotent orbit with limit mixed Hodge structure $(F, W^1)$ split over $\mathbb{R}$. Let $G_{\mathbb{R}} = \mathrm{Aut}_{\mathbb{R}}(Q)$ and $G_{\mathbb{R}}^{0,0}$ be the subset of elements which preserve $(F, W^1)$. Let $\mathcal{N}^0$ be the orbit of $N$ under the adjoint action of $G_{\mathbb{R}}^{0,0}$. Then, any element of $\mathcal{N}$ which commutes with $N$ will be of lowest weight $-2$ for the representation attached to $\theta$. More generally, we can construct several variable nilpotent orbits in this way.

*Example* 6.66. Let $D$ a period domain with Hodge numbers $h^{2,0} = 1$ and $h^{1,1} = m+1$. Then, one possible type of 1-varaible $\mathrm{sl}_2$-orbit $\theta(z) = e^{zN}F$ corresponds to the Schmid form

$$S(2) \oplus [\mathbb{R}(-1)]^m$$

with $\alpha_j = e^j f^{2-j}$ for $j = 0, 1, 2$ spanning $S(2)$ and $[\mathbb{R}(-1)]^m$ generated by elements $\gamma_1, \dots, \gamma_m$. The polarization is

$$Q(\alpha_j, \alpha_{2-j}) = (-1)^j j!(2-j)!, \qquad Q(\gamma_i, \gamma_j) = \delta_{ij}$$

and all other pairings zero. The subspace of infinitesimal isometries which are $(-1,-1)$-morphisms and lowest weight $-2$ for $(N,H)$ is spanned by $\eta_1,\dots,\eta_m$ where

$$\eta_i(\alpha_2) = \gamma_i, \qquad \eta_i(\gamma_i) = -\frac{1}{2}\alpha_0$$

and $\eta_i$ annihilates all other basis elements. It is easy to see that $[\eta_i,\eta_j] = 0$ for all $i$ and $j$. Accordingly, there is a positive real number $t$ such that

$$(N, N + t_1\eta_1 \dots, N + t_m\eta_m; F)$$

generates a pure nilpotent orbit provided $|t_j| < t$ for all $j$.

*Example* 6.67. To obtain an example where the deformation of a pure $\mathrm{sl}_2$-orbit $(N_1,N_2;F)$ with a Deligne system deformation $(N_1, N_2+\eta, Y^2)$ with $\eta$ of weight $-3$ with respect to $Y^1$, we consider an $\mathrm{sl}_2$-orbit of Schmid type

$$S(3) \oplus [S(1) \otimes \mathbb{R}(-1)]$$

for $N = N_1 + N_2$ where $S(3)$ is the isotypical component of highest weight 3 for $N_1$ and $S(1)\otimes\mathbb{R}(-1)$ is the isotypical component of highest weight 1 for $N_2$.

More concretely, we have a 6 dimensional real vector space $V_{\mathbb{R}}$ with basis $\{\alpha_3,\alpha_2,\alpha_1,\alpha_0,\beta_1,\beta_0\}$ and $\mathrm{sl}_2$ action in which $N_1$ annihilates $\{\beta_1,\beta_0\}$ and acts on $\alpha_j = e^j f^{3-j}$ by the rule $N_1(\alpha_j) = j\alpha_{j-1}$. Similarly, $N_2$ annihilates $\{\alpha_3,\alpha_2,\alpha_1,\alpha_0\}$ and acts on $\beta_j = e^j f^{1-j}$ by the rule $N_2(\beta_j) = j\beta_{j-1}$. The polarization is given by

$$Q(\alpha_j,\alpha_{3-j}) = -(-1)^j j!(3-j)!, \qquad Q(\beta_1,\beta_0) = 1$$

and all other pairings equal to zero.

In this setting, the limit mixed Hodge structure of $(F,W^2)$ is of Hodge–Tate type with $\alpha_j$ of type $(j,j)$ and $\beta_j$ of type $(j+1,j+1)$. For the mixed Hodge structure $(e^{iN_2}F, W^1)$ the Hodge type of $\alpha_j$ remains unchanged but $\beta_1 + i\beta_0$ is now of Hodge type $(2,1)$. The linear map $\eta$ which annihilates $\{\alpha_2,\alpha_1,\alpha_0,\beta_1\}$ and acts by

$$\eta(\alpha_3) = \beta_1, \qquad \eta(\beta_0) = -\frac{1}{6}\alpha_0$$

is an infinitesimal isometry which is a morphism of type $(-1,-1)$ for $(F,W^2)$. It is weight $-3$ with respect to $\operatorname{ad} Y^1$ and commutes with $N_1$. A bit more calculation shows that it lowest weight $-3$ for $(N_1,H_1)$ and highest weight 1 for $(N_2,H_2)$. So $(N_1, N_2 + \eta, Y^2)$ is a deformation of the Deligne system $(N_1, N_2, Y^2)$ with the same $\mathrm{sl}_2$-type. Since it is also polarizable, it follows that there is an $a > 0$ such that

$$(N_1, N_2 + \eta + aN_1; F)$$

generates a polarizable nilpotent orbit. We also note that unlike the original orbit $(N_1, N_2, F)$, the new orbit has maximal unipotent monodromy since now

$$Gr_4^{W^2} = N_1(Gr_6^{W^2}) + (N_2 + \eta + aN_2)(Gr_6^{W^2})$$

Let $(N_1, \dots, N_r; F)$ generate a pure nilpotent of weight $k$ which is polarized by $Q$. Let $W^j = W(\sum_{i \le j} N_i)[-k]$ and assume that $(F, W^r)$ is split over $\mathbb{R}$. Let $(\hat{N}_j, H_j)$ be the associated $\mathrm{sl}_2$-data, and $G_{\mathbb{R}} = \mathrm{Aut}_{\mathbb{R}}(Q)$. Then, for any $g \in G_{\mathbb{R}}$ the data

$$(\mathrm{Ad}(g)N_1, \dots, \mathrm{Ad}(g)N_r; g(F)) \tag{6.68}$$

generates a nilpotent orbit of weight $k$ which is polarized by $Q$.

**Theorem 6.69.** *The nilpotent orbit generated by* (6.68) *has the same limit mixed Hodge structure and* $\mathrm{sl}_2$*-data as the orbit generated by* $(N_1, \dots, N_r; F)$ *if and only if* $g$ *preserves* $F$ *and each* $\mathrm{sl}_2$*-pair* $(\hat{N}_j, H_j)$.

*Proof.* Suppose that $g \in G_{\mathbb{R}}$ preserves $F$ and each $\mathrm{sl}_2$-pair. Then, $g$ preserves

$$Y^j = k\,\mathrm{Id} + H_1 + \dots + H_j \tag{6.70}$$

and hence $g$ preserves the weight filtrations $W^j$. In particular, $g$ preserves the mixed Hodge structure $(F, W^r)$ and hence $g$ preserves $Y^r = Y_{(F,W^r)}$. By the functoriality of Deligne's construction, it then follows that

$$Y(\mathrm{Ad}(g)N_j, \mathrm{Ad}(g)Y^j) = \mathrm{Ad}(g)Y(N_j, Y_j) = \mathrm{Ad}(g)Y^{j-1} = Y^{j-1}$$

and so the corresponding chain of gradings remains the same. Likewise, the component of $\mathrm{Ad}(g)N_j$ of degree 0 with respect to $Y^{j-1} = \mathrm{Ad}(g)Y^{j-1}$ is just $\mathrm{Ad}(g)\hat{N}_j$. By assumption $\mathrm{Ad}(g)\hat{N}_j = \hat{N}_j$. Therefore (6.68) is a nilpotent orbit with the same $\mathrm{sl}_2$-data as $(N_1, \dots, N_r; F)$.

Conversely, suppose that (6.68) defines a pure nilpotent orbit with the same $\mathrm{sl}_2$-data as the orbit generated by $(N_1, \dots, N_r; F)$. Then, we must have the same associated weight filtrations

$$\begin{aligned} W(\mathrm{Ad}(g)(N_1+\cdots+N_j))[-k] &= g(W(N_1+\cdots+N_j)[-k]) \\ &= W(N_1+\cdots+N_j)[-k] \end{aligned}$$

since we have the same associated set of gradings (6.70). As in the previous paragraph, $g$ must preserve the mixed Hodge structure $(F, W^r)$ and hence $\mathrm{Ad}(g)Y^r = Y^r$. By the functoriality of Deligne systems,

$$Y(\mathrm{Ad}(g)N_r, Y^r) = Y(\mathrm{Ad}(g)N_r, \mathrm{Ad}(g)(Y^r)) = \mathrm{Ad}(g)Y^{r-1}$$

Therefore

$$\mathrm{Ad}(g)Y^{r-1} - Y^r = H_r = Y^r - Y^{r-1}$$

which implies that $\mathrm{Ad}(g)Y^{r-1} = Y^{r-1}$. Repeating this argument down the chain of gradings shows that $g$ preserves each $H_j$. As in the first paragraph of the proof, the component of $\mathrm{Ad}(g)N_j$ of degree 0 with respect to $Y^{j-1} = \mathrm{Ad}(g)Y^{j-1}$ is just $\mathrm{Ad}(g)\hat{N}_j$. To obtain the same $\mathrm{sl}_2$-data, we must therefore have $\mathrm{Ad}(g)\hat{N}_j = \hat{N}_j$ □

In summary, the set of Deligne systems with given $\mathrm{sl}_2$-data forms an algebraic variety, and this remains so once we layer on the existence of a polarization and/or a filtration $F$ with respect to which all of the $N_j$'s are horizontal. Starting from a pure nilpotent orbit with limit mixed Hodge structure which is split over $\mathbb{R}$, the real points of the algebraic group $\mathcal{G}$ consisting of isometries which preserve the limit Hodge filtration and associated $\mathrm{sl}_2$-pairs acts upon the set of nilpotent orbits with these properties.

One related question of interest is how starting from a nilpotent orbit generated by $(N_1, \ldots, N_r; F, W)$ the chain of $\mathrm{sl}_2$ representations depends on the ordering of the $N_j$'s. In the special case of a pure nilpotent orbit, permuting the order of the $\hat{N}_j$'s permutes the order of the $\mathrm{sl}_2$-pairs $(\hat{N}_j, H_j)$. However, as Example (6.61) shows, the corresponding varieties of Deligne systems can have different dimensions.

## Appendix A. Proof of Lemma 3.5

First note that $\mathcal{N}^0 \subset \mathcal{W}_N^\circ$.

Note that $Y \in \mathfrak{m}_{\mathbb{R}}^{0,\mathrm{ss}}$ acts on $\mathfrak{g}_{\mathbb{R}}^{\ell,\ell} \subset \mathfrak{m}_{\mathbb{R}}$ by the scalar $2\ell$.

*Claim.* Each $N' \in \mathcal{W}_N$ may be realized as the nilnegative element of a standard triple in $\mathfrak{m}_\mathbb{R}$ containing $Y$ as the neutral element.

*Remark.* Malcev's Theorem implies that $N$ and $N'$ are conjugate under the action of $\mathrm{Ad}(M^0_\mathbb{C})$. Unfortunately, Malcev's Theorem does not hold over $\mathbb{R}$.[9]

To prove the claim it suffices to construct $N'_+ \in \mathfrak{g}^{1,1}_\mathbb{R}$ with the property that $[N'_+, N] = Y$. Given $\ell \geq 0$, let

$$P_{2\ell}(N') \;:=\; \ker\{(N')^{2\ell+1} : \mathfrak{g}^{\ell,\ell}_\mathbb{R} \to \mathfrak{g}^{-\ell-1,-\ell-1}_\mathbb{R}\}\,.$$

Fix a basis $\{v^1_\ell, \ldots, v^{d_\ell}_\ell\}$ of $P_{2\ell}(N')$. Then

$$\bigcup_{\ell\geq 0} \big\{ (N')^k v^i_\ell \mid 1 \leq i \leq d_\ell\,,\; 0 \leq k \leq 2\ell \big\}$$

is a basis of $\mathfrak{m}_\mathbb{R}$, and we may define $N'_+ \in \mathrm{End}(\mathfrak{m}_\mathbb{R})$ by $(N')^k v^i_\ell \mapsto 2k(N')^{k-1} v^i_\ell$. Then $\{N'_+, Y, N'\}$ is a standard triple in $\mathrm{End}(\mathfrak{m}_\mathbb{R})$. Since both $Y, N' \in \mathfrak{m}_\mathbb{R} \subset \mathrm{End}(\mathfrak{m}_\mathbb{R})$, we necessarily have $N'_+ \in \mathfrak{m}_\mathbb{R}$. This proves the claim. ◇

From the proof of the claim we see that

$$\mathfrak{m}^0_\mathbb{R} \;=\; \bigoplus_{\ell \geq 0} (N')^\ell P_{2\ell}(N')\,,$$

and we deduce that the orbit

$$\mathcal{N}' \;:=\; \mathrm{Ad}(M^0_\mathbb{R}) \cdot N' \;\subset\; \mathfrak{g}^{-1,-1}_\mathbb{R}$$

has

$$\dim \mathcal{N}' \;=\; \dim_\mathbb{R} \mathfrak{m}^0_\mathbb{R} \,-\, \dim_\mathbb{R} P_0 \;=\; \dim_\mathbb{R} \mathfrak{g}^{-1,-1}_\mathbb{R}\,.$$

In particular, $\mathcal{N}'$ is open in $\mathfrak{g}^{-1,-1}_\mathbb{R}$, and this statement is independent of our choice of $N' \in \mathcal{W}_N$. The lemma now follows from the connectedness of $\mathcal{W}^\circ_N$.

[9]As an example notice that $N = \begin{pmatrix} 0 & 1 \\ 0 & 0 \end{pmatrix}$ is not conjugate to $-N$ in $\mathrm{SL}_2\mathbb{R}$ although they can both be completed to a standard triple containing $Y = \begin{pmatrix} 1 & 0 \\ 0 & -1 \end{pmatrix}$ as the neutral element.

## References

[1] Patrick Brosnan and Gregory Pearlstein. On the algebraicity of the zero locus of an admissible normal function. *Compos. Math.*, 149(11):1913–1962, 2013.

[2] Eduardo Cattani and Aroldo Kaplan. Polarized mixed Hodge structures and the local monodromy of a variation of Hodge structure. *Invent. Math.*, 67(1):101–115, 1982.

[3] Eduardo Cattani and Aroldo Kaplan. Degenerating variations of Hodge structure. *Astérisque*, (179-180):9, 67–96, 1989. Actes du Colloque de Théorie de Hodge (Luminy, 1987).

[4] Eduardo Cattani, Aroldo Kaplan, and Wilfried Schmid. Degeneration of Hodge structures. *Ann. of Math. (2)*, 123(3):457–535, 1986.

[5] Eduardo H. Cattani and Aroldo G. Kaplan. Horizontal $\mathrm{SL}_2$-orbits in flag domains. *Math. Ann.*, 235(1):17–35, 1978.

[6] David H. Collingwood and William M. McGovern. *Nilpotent orbits in semisimple Lie algebras.* Van Nostrand Reinhold Mathematics Series. Van Nostrand Reinhold Co., New York, 1993.

[7] Pierre Deligne. La conjecture de Weil. I. *Inst. Hautes Études Sci. Publ. Math.*, (43):273–307, 1974.

[8] Pierre Deligne. Personal letter to E. Cattani and A. Kaplan. 1993.

[9] Pierre Deligne, James S. Milne, Arthur Ogus and Shih Kuang-Yen. Hodge cycles, motives, and Shimura varieties. Lecture Notes in Mathematics, bf 900 (1982).

[10] Dragomir Ž. Djoković. Closures of conjugacy classes in classical real linear Lie groups. In *Algebra, Carbondale 1980 (Proc. Conf., Southern Illinois Univ., Carbondale, Ill., 1980)*, volume 848 of *Lecture Notes in Math.*, pages 63–83. Springer, Berlin, 1981.

[11] Murray Gerstenhaber. Dominance over the classical groups. *Ann. of Math. (2)*, 74:532–569, 1961.

[12] Murray Gerstenhaber. On dominance and varieties of commuting matrices. *Ann. of Math. (2)*, 73:324–348, 1961.

[13] Mark Green, Phillip Griffiths, and Matt Kerr. *Mumford-Tate groups and domains: their geometry and arithmetic*, volume 183 of *Annals of Mathematics Studies.* Princeton University Press, Princeton, NJ, 2012.

[14] Aroldo Kaplan and Gregory Pearlstein. Singularities of variations of mixed Hodge structure. *Asian J. Math.*, 7(3):307–336, 2003.

[15] Masaki Kashiwara. A study of variation of mixed Hodge structure. *Publ. Res. Inst. Math. Sci.*, 22(5):991–1024, 1986.

[16] Kazuya Kato. On SL(2)-orbit theorems. *Kyoto J. Math.*, 54(4):841–861, 2014.

[17] M. Kerr and C. Robles. Hodge theory and real orbits in flag varieties. arXiv:1407.4507, 2014.

[18] Matt Kerr and Colleen Robles. Partial orders and polarized relations on limit mixed hodge structures. In prepration, 2015.

[19] Anthony W. Knapp. *Lie groups beyond an introduction*, volume 140 of *Progress in Mathematics*. Birkhäuser Boston Inc., Boston, MA, second edition, 2002.

[20] Gregory J. Pearlstein. Degenerations of mixed Hodge structure. *Duke Math. J.*, 110(2):217–251, 2001.

[21] Colleen Robles. Classification of horizontal $\mathrm{SL}_2$'s. To appear in *Compositio Math.*, arXiv:1405.3163, 2014.

[22] C. Robles. Schubert varieties as variations of Hodge structure. *Selecta Math. (N.S.)*, 20(3):719–768, 2014.

[23] Wilfried Schmid. Variation of Hodge structure: the singularities of the period mapping. *Invent. Math.*, 22:211–319, 1973.

[24] Christine Schwarz. Relative monodromy weight filtrations. *Math. Z.*, 236(1):11–21, 2001.

[25] T. A. Springer and R. Steinberg. Conjugacy classes. In *Seminar on Algebraic Groups and Related Finite Groups (The Institute for Advanced Study, Princeton, N.J., 1968/69)*, Lecture Notes in Mathematics, Vol. 131, pages 167–266. Springer, Berlin, 1970.

[26] Joseph Steenbrink and Steven Zucker. Variation of mixed Hodge structure. I. *Invent. Math.*, 80(3):489–542, 1985.

[27] Steven Zucker. Generalized intermediate Jacobians and the theorem on normal functions. *Invent. Math.*, 33(3):185–222, 1976.

[28] Steven Zucker. Hodge theory with degenerating coefficients. $L_2$ cohomology in the Poincaré metric. *Ann. of Math. (2)*, 109(3):415–476, 1979.

*E-mail address*: pbrosnan@umd.edu

Department of Mathematics, Mathematics Building, University of Maryland, College Park, MD 20742-4015

*E-mail address*: gpearl@math.tamu.edu

Mathematics Department, Mail-stop 3368, Texas A&M University, College Station, TX 77843-3368

*E-mail address*: robles@math.duke.edu

Mathematics Department, Duke University, Box 90320, Durham, NC 27708-0320